\newif\ifarxived
\ifx\notArxived\Kaba\arxivedtrue\fi
\arxivedtrue
\documentclass[12pt]{article}
\ifarxived
\usepackage[bookmarks=true,bookmarksnumbered=true,
  bookmarkstype=toc,
  pdftitle={Strong downward Löwenheim-Skolem theorems for stationary logics, III --- mixed 
    support iteration},%
  pdfauthor={Sakaé Fuchino (\quad\quad \quad), André Ottenbreit Maschio 
  Rodrigues and Hiroshi Sakai (\quad\quad \quad\quad)},%
  pdfsubject={},%
  pdfkeywords={Strong Downward Löwenheim Skolem Theorem, stationary logic, mixed support teration, 
    Reflection Principles, supercompact cardinals}]{hyperref} 
\else
\usepackage[dvipdfmx,bookmarks=true,bookmarksnumbered=true,
  bookmarkstype=toc,
  pdftitle={Strong downward Löwenheim-Skolem theorems for stationary logics, III --- mixed 
    support iteration},%
  pdfauthor={Sakaé Fuchino (渕野 昌), André Ottenbreit Maschio 
  Rodrigues and Hiroshi Sakai (酒井 拓史)},%
  pdfsubject={},%
  pdfkeywords={Strong Downward Löwenheim Skolem Theorem, stationary logic, mixed support teration, 
    Reflection Principles, supercompact cardinals}]{hyperref} 
\fi
\ifarxived
\else
\usepackage[dvipdfmx]{pxjahyper}
\hypersetup{
  colorlinks=true,
  linkcolor=blue,
  filecolor=blue,      
  citecolor=cyan,
  urlcolor=cyan
}
\fi
\ifarxived
\usepackage{graphicx, color}
\else
\usepackage[dvipdfmx]{graphicx, color}
\fi
\usepackage{amsmath, amssymb}
\usepackage{bbm}
\usepackage{dsfont}
\usepackage{bbold}
\usepackage{mathtools}
\usepackage[normalem]{ulem}
\usepackage{setspace}

\newif\ifextended
\ifarxived\extendedtrue\fi
\extendedtrue
\definecolor{darkelectricblue}{rgb}{0.33, 0.41, 0.47} 
\newif\ifJapanese
\newif\iftesting
\testingtrue

\usepackage{theorem}
\newcommand{\bbd}[1]{{\mathbb{#1}}}
\theorembodyfont{\ifJapanese\rm\else\it\fi}
\newcount\minute	
\newcount\hour		
\newcount\hourMins  
\ifJapanese
\def\today%
{
  \the\year 年\,\zeroPadTwo{\the\month}月\,\zeroPadTwo{\the\day}日%
}
\fi
\def\now%
{
  \minute=\time    
  \hour=\time \divide \hour by 60 
  \hourMins=\hour \multiply\hourMins by 60
  \advance\minute by -\hourMins 
  \zeroPadTwo{\the\hour}:\zeroPadTwo{\the\minute}%
}
\def\zeroPadTwo#1%
{
  \ifnum #1<10 0\fi    
  #1
}
\title{Strong downward L\"owenheim-Skolem theorems for stationary logics, I}
\author{Saka\'e Fuchino (渕野 昌), Andr\'e Ottenbreit Maschio 
  Rodrigues\\and Hiroshi Sakai (酒井 拓史)${}^\ast$}
\date{}

\renewcommand{\baselinestretch}{1.2}
\renewcommand{\thefootnote}{(\arabic{footnote})\,}
\iftesting
\fi
\iftesting
\newcommand{\Label}[1]{\label{#1}\marginpar{{\renewcommand{\baselinestretch}{0.6}\tiny 
		  #1}}}
\else
\newcommand{\Label}[1]{\label{#1}}
\fi
\def\memo#1{\iftesting\marginpar{{\normalsize\renewcommand{\baselinestretch}{0.4}\tiny%
			#1\par}}\else\fi}%
\newcounter{frml}[section]
\newcounter{frmla}[section]
\def\thefrml{{\arabic{section}.\arabic{frml}}}
\def\thefrmla{{$\aleph$\arabic{section}.\arabic{frmla}}}
\def\frmlabel#1{\refstepcounter{frml}{\def\baka{#1}\ifx\baka\empty\else\label{#1}\fi}%
{\rm({\thefrml})\hfill\hfill\hfill}}
\def\frmlabela#1{\refstepcounter{frmla}{\def\baka{#1}\ifx\baka\empty\else\label{#1}\fi}%
{\rm({\thefrmla})\hfill\hfill\hfill}}
\def\xitem[#1]{\item[\frmlabel{#1}]\mbox{}%
	\iftesting\marginpar{{\renewcommand{%
				\baselinestretch}{0.6}\tiny#1}}\fi\ignorespaces}
\def\xitemq[#1]{\item[\frmlabel{#1}]\mbox{}%
	\ignorespaces}
\def\xitema[#1]{\item[\frmlabela{#1}]\mbox{}%
	\iftesting\marginpar{{\renewcommand{%
				\baselinestretch}{0.6}\tiny#1}}\fi\ignorespaces}
\def\xitemsub[#1]#2{\item[\frmlabel{#1}$_{#2}$]\mbox{}%
	\iftesting\marginpar{{\renewcommand{%
				\baselinestretch}{0.6}\tiny#1}}\fi\ignorespaces}
\def\xitemcite[#1]{\item[\rlap{\rm(\ref{#1})}\hspace*{3em}]\mbox{}%
	\iftesting\marginpar{{\renewcommand{%
				\baselinestretch}{0.6}\tiny#1}}\fi\ignorespaces}
\def\xitemciteb[#1]#2{\item[\rlap{\rm(\ref{#1}#2)}\hspace*{3em}]\mbox{}%
	\iftesting\marginpar{{\renewcommand{%
				\baselinestretch}{0.6}\tiny#1}}\fi\ignorespaces}
\def\xxitem[#1][#2]{\item[(\ref{#1}{\makebox[1.4ex][c]{#2}})]\mbox{}%
	\iftesting\marginpar{{\renewcommand{%
				\baselinestretch}{0.6}\tiny\{#1\}\{#2\}}}\fi\ignorespaces}
\def\xitemof#1{{\rm({\ref{#1}})}}
\def\xitemaof#1{{\rm({\ref{#1}})}}

\def\xitembof#1#2{$\mbox{\rm(\ref{#1}#2)}$}

\newenvironment{xitemize}{\begin{list}{}{\parsep=0.5\smallskipamount%
			\itemindent=-0.4ex%
			\itemsep=0.5\smallskipamount\leftmargin=4em\labelwidth=3em\labelsep=0.7em}}%
							 {\end{list}}
\def\assert#1{\noindent\makebox[4.8ex][r]{\rm(\makebox[2.2ex][c]{#1})}\ \ \ignorespaces}
\def\assertof#1{\makebox[4ex][c]{\rm(\makebox[2.2ex][c]{#1})}}%
\def\wassert#1{\noindent\makebox[4.8ex][r]{\em{\rm(\makebox[2.3ex][c]{#1})}}\hspace{0.8em}\ignorespaces}
\def\wassertof#1{\makebox[4.2ex][c]{\rm(\makebox[2.3ex][c]{#1})}}%

\newcommand{\bysame}[1]{\underline{\phantom{#1}}}%

\newtheorem{Thm}{\ifJapanese{\bf 定理}\else {\bf Theorem}\fi}[section]

\newtheorem{ThmA}{\ifJapanese{\bf 定理\,A}\else{\bf Theorem\,A}\fi}[section]

\newtheorem{Prop}[Thm]{\ifJapanese{\bf 命題}\else{\bf Proposition}\fi}
\newtheorem{Problem}[Thm]{\ifJapanese{\bf 未解決問題}\else{\bf Problem}\fi}

\newtheorem{Lemma}[Thm]{\ifJapanese{\bf 補題}\else{\bf Lemma}\fi}
\newtheorem{LemmaA}[ThmA]{\ifJapanese{\bf 補題\,A}\else{\bf Lemma\,A}\fi}
\newtheorem{factA}[ThmA]{Fact A}
\newtheorem{Cor}[Thm]{\ifJapanese{\bf 系}\else{\bf Corollary}\fi}

\newtheorem{Claim}{{\bf Claim}}[Thm]

\newcommand{\prf}{\ifJapanese{\bf 証明．\ }\ignorespaces\else{\bf 
		Proof.\ \ }\ignorespaces\fi}

\newcommand{\prfofClaim}{\raisebox{-.4ex}{\Large $\vdash$\ \ }}

\newcommand{\Thmof}[1]{\ifJapanese{定理\,\ref{#1}}\else{Theorem~\ref{#1}}\fi}

\newcommand{\Lemmaof}[1]{\ifJapanese{補題\,\ref{#1}}\else{Lemma \ref{#1}}\fi}
\newcommand{\LemmaAof}[1]{\ifJapanese{補題\,A\,\ref{#1}}\else{Lemma\,A\,\ref{#1}}\fi}

\newcommand{\Factof}[1]{{Fact~\ref{#1}}}
\newcommand{\FactAof}[1]{{Fact A\,\ref{#1}}}
\newcommand{\Propof}[1]{\ifJapanese{命題\,\ref{#1}}\else{Proposition~\ref{#1}}\fi}

\newcommand{\Corof}[1]{\ifJapanese{系\,\ref{#1}}\else{Corollary~\ref{#1}}\fi}

\newcommand{\Claimof}[1]{{Claim \ref{#1}}}

\newcommand{\pageof}[1]{\ifJapanese\pageref{#1}ページ\else p.\pageref{#1}\fi}
\newcommand{\sectionof}[1]{\ifJapanese{第\ref{#1}節}\else{Section~\ref{#1}}\fi}

\newcommand{\Thmabove}{{\ifJapanese 定理\else Theorem\fi\ \number\theThm}}

\newsavebox{\qedbox}\sbox{\qedbox}{
{\unitlength=0.05mm \begin{picture}(40,60)
\put(0,0){\framebox(30,44)[cc]{}}
\put(30,-7){\rule{7\unitlength}{44\unitlength}}
\put(10,-7){\rule{27\unitlength}{7\unitlength}}
\end{picture}}}
\newcommand{\qed}{\mbox{}\hfill\usebox{\qedbox}}
\newcommand{\smallqed}%
{\mbox{}\smallskip\hfill\raisebox{-.4ex}{\Large $\dashv$}}
\newcommand{\qedof}[1]%
{\mbox{} \hspace*{\fill}{\usebox{\qedbox}{\tiny~(#1)}}}
\newcommand{\Qedof}[1]%
{\mbox{} \hspace*{\fill}{\usebox{\qedbox}%
{\tiny~(#1~\number\theThm)}}}
\newcommand{\QedAof}[1]%
{\mbox{} \hspace*{\fill}{\usebox{\qedbox}%
{\tiny~(#1~\number\theThmA)}}}
\newcommand{\qedofThm}{\Qedof{\ifJapanese 定理\else Theorem\fi}}

\newcommand{\qedofCor}{\Qedof{\ifJapanese 系\else Corollary\fi}}
\newcommand{\qedofProp}{\Qedof{\ifJapanese 命題\else Proposition\fi}}
\newcommand{\qedofLemma}{\Qedof{\ifJapanese 補題\else Lemma\fi}}
\newcommand{\qedofLemmaA}{\QedAof{\ifJapanese 補題A\else Lemma\,A\fi}}

\newcommand{\qedskip}{\medskip}

\newcommand{\qedofClaim}%
{\mbox{}\hfill\raisebox{-.4ex}{\Large $\dashv$ }\nolinebreak%
\mbox{\tiny~(Claim~\number\theClaim)}}
\newcommand{\qedofSubclaim}%
{\mbox{}\hfill\raisebox{-.4ex}{\Large $\dashv$ }\nolinebreak%
\mbox{\tiny~(Subclaim~\number\theSubclaim)}}
\newcommand{\cardof}[1]{\mathopen{|\,}#1\mathclose{\,|}}

\newcommand{\Card}{{\it Card\/}}
\newcommand{\setof}[2]{\{#1\,:\,#2\}}
\newcommand{\ssetof}[1]{\{#1\}}

\newcommand{\subseteqand}[1]{\mathrel{\mathop{\subseteq}%
		\limits_{\scriptscriptstyle\hbox to 14pt{$\scriptscriptstyle #1$\hss}}}}

\newcommand{\dotcup}{\mathrel{\dot{\cup}}}

\newcommand{\mapping}[3]{#1:#2\rightarrow #3}

\newcommand{\elembed}[3]{#1:#2\stackrel{\preccurlyeq\hspace{0.8ex}}{\rightarrow}#3}
\newcommand{\combed}[3]{#1:#2\stackrel{{\scriptstyle\mathrel{{\leqslant}%
		\hspace{-0.72ex}{\lower-0.27ex\hbox{\scalebox{0.7}{$\scriptstyle\circ$}}}}\hspace{0.7ex}}}{\rightarrow}#3}

\newcommand{\fnsp}[2]{\mbox{}^{{#1}\hspace{-0.02em}}#2}
\newcommand{\imageof}{{}^{\,{\prime}{\prime}}}
\newcommand{\seqof}[2]{\langle#1\,:\,#2\rangle}
\newcommand{\pairof}[1]{\langle#1\rangle}
\newcommand{\psof}[1]{{\mathcal P}\/(#1)}
\newcommand{\incmptbl}[3]{#2\perp_{#1\,}#3}
\newcommand{\cmptbl}[3]{#2\raisebox{-0.5ex}{$\,\top_{#1}\,$}#3}
\newcommand{\forces}[2]{\,\|\hspace{-.35ex}\mbox{\sf--}_{\,#1\,}%
\mbox{\rm``}\,#2\,\mbox{\rm''}}
\newcommand{\notforces}[2]{\rlap{\ 
 /}\|\hspace{-.35ex}\mbox{\sf--}_{\,#1\,}%
 \mbox{\rm``}\,#2\,\mbox{\rm''}}
\newcommand{\modelof}[1]{\models\!\mbox{\rm``\,}#1\mbox{\rm''}}

\newcommand{\crit}{\mbox{\it crit\/}}
\newcommand{\bbone}{{\bbd1}}

\newcommand{\circleq}{\mathrel{{\leqslant}%
		\hspace{-0.86ex}{\lower-0.53ex\hbox{$\scriptscriptstyle\circ$}}}}

\newcommand{\variables}[2]{{#1}_0\ctentenc {#1}_{#2}}

\newcommand{\restr}{\restriction}
\newcommand{\dhrpr}{\downharpoonright}

\newcommand{\cf}{\mathop{cf\/}}

\newcommand{\Col}{{\rm Col}}
\newcommand{\Fn}{{\rm Fn}}
\newcommand{\Add}{{\mathrm{Add}}}

\newcommand{\trcl}{{\it trcl\/}}
\newcommand{\supp}{\mathop{\rm supp}}

\newcommand{\reals}{{\bbd{R}}}

\newcommand{\rationals}{{\bbd{Q}}}

\newcommand{\poC}{\bbd{C}}
\newcommand{\poO}{\bbd{O}}
\newcommand{\poP}{\bbd{P}}
\newcommand{\poQ}{\bbd{Q}}
\newcommand{\poR}{\bbd{R}}
\newcommand{\poS}{\bbd{S}}
\newcommand{\poT}{\bbd{T}}
\newcommand{\BaA}{{\mathbb{A}}}
\newcommand{\BaB}{{\mathbb{B}}}
\newcommand{\On}{{\rm On}}

\newcommand{\genG}{\mathbb{G}}
\newcommand{\genH}{\mathbb{H}}
\newcommand{\genK}{\mathbb{K}}
\newcommand{\genL}{\mathbb{L}}
\newcommand{\geng}{\mathbbm{g}}
\newcommand{\genh}{\mathbbm{h}}
\newcommand{\utgenG}{\utilde{\mathbb{G}}}

\newcommand{\conda}{\mathbbm{a}}
\newcommand{\condb}{\mathbbm{b}}

\newcommand{\condo}{\mathbbm{o}}
\newcommand{\condp}{\mathbbm{p}}
\newcommand{\condq}{\mathbbm{q}}
\newcommand{\condr}{\mathbbm{r}}
\newcommand{\conds}{\mathbbm{s}}
\newcommand{\condt}{\mathbbm{t}}
\newcommand{\condu}{\mathbbm{u}}
\newcommand{\boundingno}{{\mathfrak b}}

\newcommand{\LT}{{<}\,}
\newcommand{\LE}{{\leq}\,}

\newcommand{\llor}{{\bigvee\hspace{-1.5ex}\bigvee}\rule[-0.8ex]{0cm}{1ex}}
\newcommand{\lland}{{\bigwedge\hspace{-1.5ex}\bigwedge}\rule[-0.8ex]{0cm}{1ex}}

\newcommand{\ctenten}{,\mbox{}\hspace{0.08ex}{.}{.}{.}\hspace{0.1ex}}
\newcommand{\ctentenc}{,{}\linebreak[0]\hspace{0.04ex}{{.}{.}{.}\hspace{0.1ex},\,}\linebreak[0]}

\newcommand{\xmbox}[1]{ $\relax{\rm #1}\relax$ }
\newcommand{\qbox}[1]{\mbox{``\hspace{0.1em}}#1\mbox{''}}
\newcommand{\gmA}{\mathfrak{A}}
\newcommand{\gmB}{\mathfrak{B}}

\newcommand{\gmP}{\mathfrak{P}}
\newcommand{\gmQ}{\mathfrak{Q}}

\newcommand{\calC}{{\mathcal C}}
\newcommand{\calD}{{\mathcal D}}

\newcommand{\calF}{{\mathcal F}}

\newcommand{\calH}{{\mathcal H}}
\newcommand{\calI}{{\mathcal I}}
\newcommand{\calL}{{\mathcal L}}

\newcommand{\calO}{{\mathcal O}}
\newcommand{\calP}{{\mathcal P}}
\newcommand{\calQ}{{\mathcal Q}}

\newcommand{\calS}{{\mathcal S}}

\newcommand{\calU}{{\mathcal U}}

\newcommand{\uta}{\utilde{a}}
\newcommand{\utb}{\utilde{b}}

\newcommand{\utildef}{\utilde{f}}
\newcommand{\utf}{\utilde{f}}

\newcommand{\utg}{\utilde{g}}
\newcommand{\uth}{\utilde{h}}

\newcommand{\utu}{\utilde{u}}

\newcommand{\utS}{\utilde{S}}

\newcommand{\utpoQ}{\utilde{\mathbb Q}}
\newcommand{\utpoR}{\utilde{\mathbb R}}
\newcommand{\utpoS}{\utilde{\mathbb S}}

\newcommand{\utcondp}{\utilde{\condp}}
\newcommand{\utcondq}{\utilde{\condq}}
\newcommand{\utcondr}{\utilde{\condr}}

\newcommand{\varin}{\mathrel{\varepsilon}}

\newcommand{\ZFC}{{\sf ZFC}}

\newcommand{\ADS}{{\sf ADS}}

\newcommand{\MA}{{\sf MA}}
\newcommand{\MM}{{\sf MM}}
\newcommand{\PFA}{{\sf PFA}}
\newcommand{\PKL}{{\sf PKL}}

\newcommand{\GRP}{{\sf GRP}}

\newcommand{\Refl}{{\mathfrak{R}\mathfrak{e}\mathfrak{f}\mathfrak{l}\,}}

\newcommand{\SDLS}{{\sf SDLS}}
\newcommand{\HH}{{\sf HH}}
\newcommand{\intnl}{int}
\newcommand{\stat}{{stat}}

\newcommand{\continuum}{2^{\aleph_0}}

\newcommand{\HP}{{\sf HP}}

\newcommand{\st}{such that}
\newcommand{\wrt}{with respect to}
\newcommand{\Wolog}{Without loss of generality}

\newcommand{\nbhd}{neighborhood}

\newcommand{\uniV}{{\sf V}}
\newcommand{\po}{poset}

\newcommand{\pos}{posets}

\newcommand{\Pkl}[2]{\ifx\bakakaba#1\bakakaba\ifx\bakakaba#2\bakakaba{\mathcal 
    P}_\kappa(\lambda)\else{\mathcal P}_\kappa(#2)\fi\else{\mathcal P}_{#1}(#2)\fi}

\newcommand{\utildeT}[1]{%
	\hbox to 0pt{$\mathop{#1}\limits_{\raise0.25ex\hbox{$\scriptstyle\sim$}}$\hss}%
		\relax\phantom{\underline{#1}}}
\newcommand{\utildeS}[1]{%
	\hbox to 0pt{\smash{$\mathop{\scriptstyle #1}\limits_{%
				\raisebox{0.6ex}[0pt]{$\scriptscriptstyle\sim$}}$}\hss}%
		\relax\phantom{\mathord{{#1}_{\rule[-0.6ex]{0pt}{1pt}}}}}
\newcommand{\utildeSS}[1]{%
	\hbox to 0pt{$\mathop{\scriptscriptstyle #1}%
		\limits_{\scriptscriptstyle\sim}$\hss}%
		\relax\phantom{\underline{#1}}}
\newcommand{\utilde}[1]{%
	\mathchoice{\utildeT{#1}}{\utildeT{#1}}{\utildeS{#1}}{\utildeSS{#1}}}

{\end{minipage}\end{trivlist}}

\begin{document}
\title{Strong downward L\"owenheim-Skolem theorems\\ for stationary logics, III\\--- 
  mixed support iteration}
\ifarxived
\author{Saka\'e Fuchino\ifextended\ (\quad\quad \quad)\fi, Andr\'e Ottenbreit Maschio 
  Rodrigues$^\dagger$ \ifextended\\\fi and 
  Hiroshi Sakai\ifextended\ (\quad\quad \quad\quad)\fi
  }
\else
\author{Saka\'e Fuchino\ifextended\ (渕野 昌)\fi, Andr\'e Ottenbreit Maschio 
  Rodrigues$^\dagger$ \ifextended\\\fi and 
  Hiroshi Sakai\ifextended\ (酒井 拓史)\fi
  }
\fi
\maketitle
\renewcommand{\thefootnote}{$\ast$\ }
  \footnotetext{Graduate School of System Informatics, Kobe University \\Rokko-dai 1-1, Nada, Kobe 657-8501 Japan
   \\
    \scalebox{0.95}[1]{\tt fuchino\@diamond.kobe-u.ac.jp,\,andreomr\@gmail.com,\,hsakai\@people.kobe-u.ac.jp}}


\begin{abstract}
Continuing [Fuchino, Ottenbreit and Sakai\cite{III:bib-I,III:bib-II}] and [Fuchino and\\ 
  Ottenbreit\cite{III:bib-refl-conti}], we further study reflection principles  
in connection with the L\"owenheim-Skolem Theorems of stationary logics. In 
this paper, we mainly analyze the situations in the models obtained by mixed support 
iteration of a supercompact length and then collapsing another supercompact 
cardinal to make it $(\continuum)^+$. We show, among other things, that the reflection down to $\LT\continuum$ of 
the non-metrizability of topological spaces with small character is independent from the 
reflection properties studied in [Fuchino, Ottenbreit and Sakai\cite{III:bib-I,III:bib-II}] and [Fuchino 
  and Ottenbreit\cite{III:bib-refl-conti}]. 
\end{abstract}


\ifextended
{\color{darkelectricblue}

\begin{quote}
	\footnotesize
	\noindent
	\centerline{\normalsize\tt Contents\hspace{6em}\mbox{}}\mbox{}\\
       {\mbox{}\hspace{-1.6em}\tt\makebox[3.4ex][l]{\ref{III:intro}.}%
         Introduction}\ \ \dotfill\ \ \pageref{III:intro}\\ 
       {\mbox{}\hspace{-1.6em}\tt\makebox[3.4ex][l]{\ref{III:hamburger}.}%
         Reflection number of Hamburger's Hypothesis}\ \ \dotfill\ \ \pageref{III:hamburger}\\ 
       {\mbox{}\hspace{-1.6em}\tt\makebox[3.4ex][l]{\ref{III:preserv}.}%
         Preservation and non-preservation of stationarity of subsets of \\$\Pkl{}{}$ }\ \ \dotfill\ \ \pageref{III:preserv}\\ 
       {\mbox{}\hspace{-1.6em}\tt\makebox[3.4ex][l]{\ref{III:two}.}%
         Two dimensional Laver-generic large cardinals}\ \ \dotfill\ \ \pageref{III:two}\\ 
       {\mbox{}\hspace{-1.6em}\tt\makebox[3.4ex][l]{\ref{III:msi}.}%
         Mixed support iteration}\ \ \dotfill\ \ \pageref{III:msi}\\ 
       {\mbox{}\hspace{-1.6em}\tt\makebox[3.4ex][l]{\ref{III:refl}.}%
         Models with strong reflection properties down to $\LT\continuum$ and 
         with\\even stronger reflection properties but down to $\LE\continuum$}\ \ \dotfill\ \ \pageref{III:refl}\\ 
       {\mbox{}\hspace{-1.6em}\tt
         References}\ \ \dotfill\ \ \pageref{III:ref}\\ 
\end{quote}}
\fi
\renewcommand{\thefootnote}{}
\footnotetext{{\it Date:} August 26, 2020
  \qquad {\it Last update:} 
  \today\ (\now\ JST)\vspace{-1\smallskipamount}
}
\footnotetext{\vspace{-1\smallskipamount}}
\footnotetext{The first author was partially supported by JSPS Kakenhi Grant No.\ 20K03717.
  The second author was partially supported by the Monbukagakusho (Ministry of Education, 
  Culture, Sports, Science and Technology) Scholarship, Japan. The third author is 
  supported by JSPS Kakenhi Grant No.\ 18K03397. The authors would like to thank the 
  anonymous referee for valuable comments. 
\ifextended\\\par
  {\color{darkelectricblue} This is an extended version of the paper with the same title. Some extra remarks and 
    details omitted in the published version, as well as further corrections, may be found 
    in this version. The additional 
    stuff is typeset in dark electric blue like this paragraph. The most recent edition of this 
    version is downloadable as:\smallskip\\ \qquad 
    \url{https://fuchino.ddo.jp/papers/SDLS-III-xx.pdf}}\else An extended version of the paper 
  with some details omitted in the published version, and with some possible further update is 
  downloadable as: \url{https://fuchino.ddo.jp/papers/SDLS-III-xx.pdf} 
\fi}

\renewcommand{\thefootnote}{\arabic{footnote})\,}
\section{Introduction}\Label{III:intro}
Reflection properties of the following type are considered in various mathematical contexts: 
\begin{xitemize}
\xitem[III:Refl-0]If a structure $\gmA$ in the class $\calC$ has the property $\calP$, then there is a 
  structure $\gmB$ in relation $\calQ$ to $\gmA$ \st\ $\gmB$ has 
  the cardinality $<\kappa$ and $\gmB$ also has the property $\calP$. 
\end{xitemize}
We shall call ``$\LT\kappa$'' above the 
{\it reflection point} of the reflection property \xitemof{III:Refl-0}. If $\kappa$ is a successor 
cardinal $\mu^+$, we shall also say that the reflection point of the reflection property is $\leq\mu$.

An instance of \xitemof{III:Refl-0} is when $\calC=$ ``first countable
topological spaces'', $\calP=$ ``non-metrizable'', $\calQ=$ ``subspace'' and
$\kappa=\aleph_2$, that is, with the reflection point $\leq\aleph_1$. In this 
setting, the obtained reflection statement is: 
\begin{xitemize}
\xitem[III:Refl-1] For any first countable topological space $X$, if $X$ is 
  non-metrizable, then there is a subspace $Y$ of $X$ of cardinality $<\aleph_2$ 
  \st\ $Y$ is also non-metrizable.
\end{xitemize}
The consistency of the statement above is still unknown. This persistently open 
problem about the consistency of the assertion \xitemof{III:Refl-1} is called Hamburger's Problem 
after Peter Hamburger who asked a related question (see [Hajnal-Juh\'asz\cite{III:bib-hajnal-juhasz}]).  

The naturalness of the question can be seen in the following known partial 
solutions: With ``first countable'' replaced by ``compact'', the 
assertion \xitemof{III:Refl-1} is a theorem in \ZFC\ [Dow\cite{III:bib-dow}].
With ``first countable'' replaced by ``locally-compact'', the assertion 
\xitemof{III:Refl-1} is independent from \ZFC\ (for the consistency we need some 
very large cardinal since $\square_\kappa$ for some 
$\kappa$ implies the 
negation of the statement, see 
[Fuchino-Juh\'asz-Szentmikl\'ossy-Usuba\cite{III:bib-erice}]).

We shall call the following principle ``Hamburger's Hypothesis'' (with the 
reflection point $\LT\kappa$): 
\begin{xitemize}
\item[$\HH(\LT\kappa)$\,: ] For any topological space $X$ with $\chi(x,X)<\kappa$ for 
  all $x\in X$, if $X$ is non-metrizable then there is a subspace $Y$ of $X$ of 
  cardinality $<\kappa$ which is also non-metrizable. 
\end{xitemize}
Recall that the character 
$\chi(x,X)$ of a point $x$ in a topological space $X$ is  the minimal possible cardinality of 
a neighborhood  base of $x$ in $X$. Without the 
condition on the character of points, we easily obtain a counter-example to 
the reflection of non-metrizability (see [Hajnal-Juh\'asz\cite{III:bib-hajnal-juhasz}]). 

Note that the original Hamburger's Problem \xitemof{III:Refl-1} is equivalent to
$\HH(\LT\aleph_2)$ ([Hajnal-Juh\'asz\cite{III:bib-hajnal-juhasz}]).
\ifextended{\small\color{darkelectricblue}{[\![} If $\HH(\LT\aleph_2)$ holds, then we clearly have
  \xitemof{III:Refl-1}.
  \par Assume that \xitemof{III:Refl-1} holds. To see that $\HH(\LT\aleph_2)$ holds, 
  suppose that $X$ is a non-metrizable space with $\chi(x, X)<\aleph_2$ for all $x\in X$. 
  If $\chi(x,X)=\aleph_1$ for some $x\in X$, then there is a subspace $Y$ of $X$ of 
  cardinality $\aleph_1$ with $x\in Y$ and $\chi(x,Y)=\aleph_1$ (an elementary submodel 
  argument proves this easily: Let $\theta$ be sufficiently large and let
  $M\prec\calH(\theta)$ be \st\ $\cardof{M}=\aleph_1$, $\omega_1\subseteq M$, and
  $\pairof{X,\tau}\in M$. Then, $Y=X\cap M$ is as desired). By this $x$, $Y$ is not metrizable. 
  If $\chi(x,X)<\aleph_1$ for all $x\in X$, then, by the assumption, there is a non-metrizable 
  subspace $Y$ of $X$ of cardinality $\LE\aleph_1$. 
{]\!]}}\fi

\ifextended{\small\color{darkelectricblue}
\mbox{$\HH(<\aleph_1)$} does not hold:
$\omega_1$ in order topology as well as 
$X_\calF$ in the proof of \Thmof{III:P-vDw-2} for an unbounded
  $\calF\subseteq\fnsp{\omega}{\omega}$ is  a counter\-example. 
}\fi 

The following fact will be used in the proofs of 
\Corof{III:refl-a-1} and \Propof{III:P-gen-large-3}:
\begin{Thm}[{[}Dow, Tall and Weiss\cite{III:bib-dtw2}{]}]
  \Label{III:P-refl-a} Suppose that $X$ is a non-metrizable space, 
  $\delta\in\Card$ and $\poP=\Fn(\delta,2)$, the \po\ with finite conditions 
  adding $\delta$ many Cohen reals. Then we have 
  \begin{xitemize}
  \xitem[III:refl-a-0] $\forces{\poP}{\check{X}\mbox{ is non-metrizable}}$. \qed
  \end{xitemize}
\end{Thm}

Topological space  $X$ is considered here as a pair $X=\pairof{X,\tau}$ where
$\tau$ is the open base of the topology. Note that the family $\calO$ of all 
open sets in the ground model need not to satisfy the axioms of open sets in a 
generic extension, while an open base remains to be an open base in the generic extension. 

Let us call the \pos\ of the form $\Fn(\delta,2)$ for some ordinal $\delta$ 
{\it generalized Cohen \pos}. 

For a class $\calP$ of \pos, a cardinal $\kappa$ is said to be {\it generically 
supercompact by $\calP$}, if, for any $\lambda\geq\kappa$, there is a \po\ $\poP\in\calP$
\st, for a $(\uniV,\poP)$-generic $\genG$, there are classes $j, M\subseteq\uniV[\genG]$ 
\st\ 
\begin{xitemize}
\xitem[III:refl-a-2] $M$ is transitive in $\uniV[\genG]$ and $\elembed{j}{V}{M}$,
\xitem[III:refl-a-3] $crit(j)=\kappa$, 
\xitem[III:refl-a-4] $j(\kappa)>\lambda$, and  
\xitem[III:refl-a-5] $j\imageof\lambda\in M$. 
\end{xitemize}

\begin{Cor}
  \Label{III:P-refl-a-0}
  If $\kappa$ is generically supercompact by generalized Cohen \pos, then
  $\HH(\LT\kappa)$ holds. 
\end{Cor}
\prf Suppose that $X$ is a non-metrizable space with 
\begin{xitemize}
\xitem[III:refl-a-1] 
  $\chi(x,X)<\kappa$ for all
  $x\in X$. 
\end{xitemize}

\Wolog, $X=\pairof{\theta,\tau}$ for some ordinal $\theta$ and an open base $\tau$ on $\theta$. Let 
$\lambda\geq\theta$ be sufficiently large and let $\poP=\Fn(\mu,2)$ for some 
  cardinal $\mu$ \st, for a $(\uniV,\poP)$-generic filer $\genG$, there are 
  classes $j$, $M\subseteq\uniV[\genG]$ satisfying 
  \xitemof{III:refl-a-2}, \xitemof{III:refl-a-3}, \xitemof{III:refl-a-4}, and
  \xitemof{III:refl-a-5} for this $\lambda$. 

Let $\tau''=\setof{j(O)\cap j\imageof\theta}{O\in\tau}$. Then we have
$\pairof{j\imageof\theta,\tau''}$, $\pairof{\theta,\tau}\in M$, and
$M\models\pairof{\theta,\tau}\cong\pairof{j\imageof\theta,\tau''}$ by 
\xitemof{III:refl-a-5} (see, e.g.\ Lemma 2.5 in [Fuchino, Ottenbreit and Sakai\cite{III:bib-II}]).

By \Thmof{III:P-refl-a},
$\uniV[\genG]\modelof{\pairof{j\imageof\theta,\tau''}\mbox{ is non-metrizable}}$. 
By \xitemof{III:refl-a-1}, $M\modelof{\pairof{j\imageof\theta,\tau''}\xmbox{ is a
    sub-space of }\pairof{j(\theta),j(\tau)}}$. 

Thus, $M\modelof{\xmbox{there is a non-metrizable subspace }Y\xmbox{ of }j(X)\xmbox{ of cardinality }\LT j(\kappa)}$.
By elementarity, it follows that 
$\uniV\modelof{\xmbox{there is a non-metrizable subspace }Y\xmbox{ of }X\xmbox{ of cardinality }\LT \kappa}$.
\qedofCor\qedskip

In a model obtained as the generic extension by $\Fn(\kappa,2)$ where $\kappa$ is 
a supercompact cardinal, we have $\continuum=\kappa$ and $\kappa$ is generically 
supercompact by generalized Cohen \pos. 
Thus, 

\begin{Cor}[{[}Dow, Tall and Weiss\cite{III:bib-dtw2}{]}]
  \Label{III:P-refl-a-1}
  If\/ $\ZFC$ $+$ ``there is a supercompact cardinal'' is consistent, then so is 
  $\ZFC$\ $+$ $\HH(\LT\continuum)$. \qed
\end{Cor}

The Strong Downward L\"owenheim-Skolem Theorem 
$\SDLS^-(\calL^{\aleph_0}_{stat},\LT\kappa)$ for the stationary logic 
$\calL^{\aleph_0}_{\stat}$ down to $\LT\kappa$ is another natural reflection 
property. Here, the stationary logic $\calL^{\aleph_0}_{stat}$ is a monadic second 
order logic whose second order variables run over countable subsets of the 
underlining set of the structure in question. The only second-order quantifier in the logic 
is 
`$\stat$' (as well as its dual `$aa$' where the quantification ``$aa\,X$'' is introduced as 
the abbreviation of 
``$\neg\,stat\,X\,\neg$''). The semantics of the logic is introduced by 
the following step in the recursion in  
addition to the 
usual recursive definition of the semantics for first order part of the logic: for a structure $\gmA=\pairof{A\ctenten}$ 
and $\calL^{\aleph_0}_{stat}$-formula 
$\varphi=\varphi(x_0\ctenten, X_0\ctenten,X)$ in the corresponding signature, where $X_0\ctenten$, $X$ are the 
second order variables in $\varphi$, as well as for $a_0\ctenten\in A$ and
$U_0\ctenten\in[A]^{\aleph_0}$, 
\begin{xitemize}
\xitem[III:refl-a-6] 
    $\gmA\models stat\,X\,\varphi(a_0\ctentenc U_0\ctentenc X)$\\ 
  $\Leftrightarrow$\ \
  $\setof{U\in[A]^{\aleph_0}}{\gmA\models\varphi(a_0\ctentenc U_0\ctentenc U)}$ 
  is stationary in $[A]^{\aleph_0}$. 
\end{xitemize}

For a substructure $\gmB=\pairof{B\ctenten}$ of $\gmA$, the weak variant of 
elementary submodel relation $\prec^-_{\calL^{\aleph_0}_{\stat}}$ between $\gmB$ and $\gmA$ is 
defined by
\begin{xitemize}
\xitem[III:refl-a-7] 
  $\gmB\prec^-_{\calL^{\aleph_0}_{\stat}}\gmA$\ \ $\Leftrightarrow$ \\[\jot]
    $\gmB\models\varphi(\variables{b}{n-1})$ holds if and only if
  $\gmA\models\varphi(\variables{b}{n-1})$ holds 
  for all $\calL^{\aleph_0}_{\stat}$-formulas $\varphi=\varphi(\variables{x}{n-1})$ 
  without free second-order variables,
  and for all $\variables{b}{n-1}\in B$. 
\end{xitemize}

The reflection principle $\SDLS^-(\calL^{\aleph_0}_{stat},\LT\kappa)$ for a 
cardinal $\kappa\geq\aleph_2$ is defined by:
\begin{xitemize}
\item[$\SDLS^-(\calL^{\aleph_0}_{stat},\LT\kappa)$:\ \ ] For any structure $\gmA$ 
  in a countable signature, there is a substructure $\gmB$ of $\gmA$ of 
  cardinality $\LT\kappa$ \st\ $\gmB\prec^-_{\calL^{\aleph_0}_{stat}}\gmA$. 
\end{xitemize}

In [Fuchino, Ottenbreit and Sakai\cite{III:bib-I}], we also considered the version of 
SDLS without `$-$' by allowing second order free variables and second order 
parameters in the formulas $\varphi$ in \xitemof{III:refl-a-7}. However, it is proved there that the 
principle $\SDLS(\calL^{\aleph_0}_{stat},\LT\kappa)$ obtained in this way for a 
regular $\kappa$ is 
simply the conjunction of $\SDLS^-(\calL^{\aleph_0}_{stat},\LT\kappa)$ and
$\mu^{\aleph_0}<\kappa$ for all $\mu<\kappa$.

In the standard 
model of \PFA\ or under strongly Laver-generically supercompactness of a cardinal $\kappa$ 
for proper \pos\ (for definition of Laver-generic supercompactness, see 
\pageof{III:gen-large-6-0}),  we have the reflection principle 
$\SDLS^-(\calL^{\aleph_0}_\stat,\LT\aleph_2)$. Actually, 
$\MA^{+\omega_1}(\sigma\mbox{-closed})$ already implies this principle, and 
strongly Laver-generically supercompactness for properness of $\kappa$ implies $\kappa=\aleph_2$ 
and $\PFA^{+\omega_1}$.

If
$\MA^{+\omega_1}(\sigma\mbox{-closed})$ (or
  $\PFA^{+\omega_1}$, or $\MM^{+\omega_1}$, resp.) 
holds and $\poP$ is $\LT\aleph_2$-directed 
closed, then we have $\forces{\poP}{\MA^{+\omega_1}(\sigma\mbox{-closed})}$\ \ (or 
$\forces{\poP}{\PFA^{+\omega_1}}$, or $\forces{\poP}{\MM^{+\omega_1}}$ resp.) 
(Proposition 15 in [Fuchino and Ottenbreit\cite{III:bib-refl-conti}]). 

Suppose that $\MA^{+\omega_1}(\sigma\mbox{-closed})$ holds and
$\continuum=2^{\aleph_1}=\aleph_2$. 
and there is a supercompact cardinal $\kappa_1$.  
Let $\poP=\Col(\continuum, \kappa_1)$. In a generic extension by $\poP$, we 
still have $\MA^{+\omega_1}(\sigma\mbox{-closed})$ by the result mentioned above, 
and hence also $\SDLS^-(\calL^{\aleph_0}_\stat,\LT\aleph_2)$. 
On the other hand, $\poP$ forces
$\kappa_1$ to be $(\continuum)^+$ and makes $\kappa_1$ generically supercompact by
$\LT\aleph_2$-closed \pos\ (see, e.g.\ Lemma 4.10 in [Fuchino, Sakai and 
  Ottenbreit\cite{III:bib-I}]). By Theorem 4.13 in [Fuchino, Sakai and 
  Ottenbreit\cite{III:bib-I}], the assertion that $\kappa_1=\kappa^+$ is generically supercompact 
by $\LT\kappa$-closed \pos\ is equivalent to the Game Reflection Principle
$\GRP^{\LT\kappa}(\leq\kappa)$ under $2^{\LT\kappa}=\kappa$.  Thus, in this way, 
we obtain a model of a very strong reflection property with the reflection point
$\LT\continuum$, together with an even stronger reflection property but with the 
reflection point $\LE\continuum$. 

$\SDLS^-(\calL^{\aleph_0}_{stat},\LT\continuum)$ implies $2^{\aleph_0}=\aleph_2$ (Corollary 
2.3 in in [Fuchino, Sakai and  
  Ottenbreit\cite{III:bib-II}]). This 
means in particular that, if the continuum should be lager than $\aleph_2$, this 
reflection statement is not available. 
In the model obtained by 
iterating ccc \pos\ supercompact times with finite support along with a 
book-keeping provided by a Laver-function, the continuum is 
extremely large (e.g. in terms of existence of a saturated ideal) 
but the Strong Downward L\"owenheim-Skolem Theorem
$\SDLS^{int}(\calL^{\aleph_0}_{\stat},\LT\continuum)$ of the stationary logic
$\calL^{\aleph_0}_{\stat}$ with internal interpretation (a weakening of
$\SDLS^-(\calL^{\aleph_0}_{stat},\LT\continuum)$) holds, 
as well as the Strong Downward 
L\"owenheim-Skolem Theorem $\SDLS^{int}_+(\calL^{\PKL}_{stat},\LT\continuum)$ of 
internal interpretation of the \PKL-logic with the reflection 
point $\LE\continuum$  
(Theorem 2.10 and Proposition 3.1 in [Fuchino, Sakai and  
  Ottenbreit\cite{III:bib-II}] for $\SDLS^-(\calL^{\aleph_0}_{stat},\LT\continuum)$; 
Proposition 4.1 and Theorem 4.5 in [Fuchino, Sakai and  
  Ottenbreit\cite{III:bib-II}] for $\SDLS^{int}_+(\calL^{\PKL}_{stat},\LT\continuum)$ ) together 
with $\MA^{+\mu}$ for all $\mu<\continuum$.  The significance of
$\SDLS^{int}_+(\calL^{\PKL}_{stat},\LT\continuum)$ in this connection is that it implies 
that the continuum is very large (e.g.\ it implies that the continuum is at least weakly Mahlo).

For this model,  
there seems to be no way to force further to obtain a stronger reflection 
but with the reflection point $\LE\continuum$ without destroying the reflection 
properties already existing in the model.

In the present paper, we show that the mixed support supercompact 
time iteration, roughly speaking, with Easton support mixed with the finite support, 
bookkept along with a Laver function together with a further collapse of the second supercompact 
cardinal creates a model in which ``down to $\LT\continuum$'' type of reflection 
principles as mentioned above together with $\GRP^{\LT\continuum}(\LE\continuum)$ hold. 

Modifying the finite support part of this iteration, we show the independence of
$\HH(\LT\continuum)$ from the other strong reflection properties. 

For the definition of some of the set-theoretic principles and basic facts around them 
remained unexplained in the present paper, the reader should 
consult [Fuchino, Sakai and Ottenbreit\cite{III:bib-I,III:bib-II}]. These papers in extended 
version uploaded at the URLs given in the References may be also helpful since they contain some more 
details which were omitted in the submitted version of the papers. 

In particular, we are going to drop the definition of
$\SDLS^{int}_+(\calL^{\PKL}_{stat},\LT\continuum)$ and ask readers to consult [Fuchino, 
  Sakai and Ottenbreit\cite{III:bib-II}] for details. However, we shall cite the following infinitary 
combinatorial characterization of this principle. This will be used in 
\Propof{III:P-gen-large-3},\,\assertof{2} to show that this principle holds under certain instance of 
the two-dimensional Laver-generic large cardinal considered in \sectionof{III:two}. 

Extending the standard notation, for sets $s$ and $t$, we denote with $\Pkl{s}{t}$ the set
\begin{xitemize}
\xitem[III:refl-a-8] 
  $[\,t\,]^{\cardof{s}}=\setof{a\in\psof{t}}{\cardof{a}<\cardof{s}}$. 
\end{xitemize}

\begin{Lemma}[Proposition\,4.1 in {[}Fuchino,\,Sakai\,and\,Ottenbreit\cite{III:bib-II}{]}]
  For a regular cardinal $\kappa>\aleph_1$ $\SDLS^{int}_+(\calL^{\PKL}_{stat},\LT\kappa)$ 
  is equivalent to the assertion that $(\ast)^{int+\PKL}_{\LT\kappa,\lambda}$ holds for all 
  regular $\lambda\geq\kappa$ where
\begin{xitemize}
\item[$(\ast)^{\intnl+\PKL}_{\LT\kappa,\lambda}$: ] For any countable expansion $\gmA$ 
  of the structure $\pairof{\calH(\lambda),\kappa,\in}$ and any family
  $\seqof{S_a}{a\in\calH(\lambda)}$ \st\ $S_a$ is a stationary subset of
  $\Pkl{\kappa}{\calH(\lambda)}$ for all $a\in\calH(\lambda)$, there are stationarily many 
  $M\in\Pkl{\kappa}{\calH(\lambda)}$ \st\ $\cardof{\kappa\cap M}$ is regular, $\gmA\restr M\prec\gmA$ and
  $S_a\cap \Pkl{\kappa\cap M}{M}\cap M$ is stationary in $\Pkl{\kappa\cap M}{M}$ 
  for all $a\in M$.\qed
\end{xitemize}
\end{Lemma}

We shall use freely the following ``bullet notation'' of names in forcing construction, introduced by Asaf 
Karagila\,\footnote{The authors learned this extremely helpful notation in a tutorial lectures by Asaf 
  Karagila in Kyoto at the RIMS Set Theory Workshop 2019.}.

If $t(x_0\ctenten)$ is a term in some conservative expansion of the 
  language and the axiom system of the set theory by definitions then for a \po\ $\poP$ 
  and $\poP$-names $\uta_0$\ctenten, $t(\uta_0\ctenten)^\bullet$ denotes the 
  standard $\poP$-name $\utu$ \st\
  \begin{xitemize}
  \xitem[III:refl-a-9] 
    $\utu[\genG]=t^{\uniV[\genG]}(\uta_0[\genG]\ctenten)$ for any $(\uniV,\poP)$-generic 
    filter $\genG$ \\[\jot]
(or, more syntactically, 
    $\forces{\poP}{\utu\equiv t(\uta_0\ctenten)}$). 
  \end{xitemize}
  For example, 
  $\pairof{\uta,\utb}^\bullet$ is denoted as ${\rm op}(\uta,\utb)$ in [Kunen\cite{III:bib-kunen}].
  $t(\uta_0\ctenten)$ may have infinitely many parameters. For  
  example, if $\uta_\xi$, $\xi<\delta$ is a sequence of $\poP$-names in the ground 
  model, $\setof{\uta_\xi}{\xi<\delta}^\bullet$ may be introduced as 
  the $\poP$-name $\setof{\pairof{\uta_\xi,\bbone_\poP}}{\xi<\delta}$, while 
  $\seqof{\uta_\xi}{\xi<\delta}^\bullet$ may be introduced as the $\poP$-name
  $\setof{\pairof{\pairof{\check{\xi},\uta_\xi}^\bullet,\bbone_\poP}}{\xi<\delta}$. The choice of the 
  exact definition of each bullet name is left to the reader. We only assume that the 
  choice is done in a consistent way. If we want to emphasize that the bullet name
  $t(a_0\ctenten)^\bullet$ is a $\poP$-name, we put the subscript $\poP$ and write 
  $t(a_0\ctenten)^\bullet_\poP$. 

  For a \po\ $\poP$, $\poP$-check names of a ground model set $a$ is represented either simply 
  by $a$ or with a check as $\check{a}$. If it is necessary to make clear which \po\ is 
  involved, we shall also write $(a)^\surd_{\poP}$. This representation is used, 
  in particular, if a ground model set is given by a term. Thus we write, e.g.
  $(\psof{a})^\surd_\poP$, $(a\cup b)^\surd_\poP$, 
  $(\setof{x\in a}{\varphi(x,\ldots)})^\surd_\poP$ etc.

  A part of the results in the following,
  most of the materials in \sectionof{III:msi} in 
  particular, have been presented in the PhD thesis [Ottenbreit 
    Maschio Rodrigues\cite{III:bib-andre}] 
  of the second author, although some arguments and details are treated 
  differently from those in the PhD thesis.  

\section{Reflection number of Hamburger's Hypothesis}
\Label{III:hamburger}
\ifextended{\small\color{darkelectricblue} 

In the following, we shall examine the details of the example of the topological 
space given on p.158 in [van Douwen\cite{III:bib-vdouwen}]. 

A topological space $X=\pairof{X,\tau}$ where $\tau$ is an open base of the 
topology is said to be a {\it Moore space}, if $X$ is a regular Hausdorff space \st\ 
\begin{xitemize}
\xitema[III:vDw-0] there is a sequence $\calO_n\subseteq\tau$ of open covers of $X$ 
  (i.e. $\bigcup\calO_n=X$ for all $n\in\omega$) with the property that, for any 
  closed $C\subseteq X$ and $x\in X\setminus C$, there is $n\in\omega$ \st\ all 
  $O\in\calO_n$ with $x\in O$ is disjoint with $C$.
\end{xitemize}
The property \xitemof{III:vDw-0} is called the {\it developability} of $X$. 
If \xitemof{III:vDw-0} holds, we say that $\seqof{\calO_n}{n\in\omega}$ is a {\it development\/} of $X$ and $X$ 
is {\it developable}. 

The following is a warm-up exercise:

\begin{LemmaA}\wassertof{1} If $X$ is a metrizable space, then $X$ is a Moore space.\smallskip

\wassert{2} If $X$ is a Moore space then it is first countable.\ifextended\else\qed\fi

\end{LemmaA}
\ifextended{\small\color{darkelectricblue}
\prf \assertof{1}: Suppose that $X$ is a metrizable space. Then $X$ is Hausdorff and normal. 

Let $d$ be a metric on $X$ which induces the topology of $X$. Then
$\calO_n=\setof{S_d(x,\frac{1}{n+1})}{x\in\omega}$, for $n\in\omega$ form 
a development of $X$.\smallskip

\assertof{2}: Suppose that $\calO_n$, $n\in\omega$ witness that $X$ is a Moore space.
Let $x\in X$. 
For each $n\in\omega$, let $O_{x,n}$ for $n\in\omega$ be \st\ $x\in O_{x,n}$ and
$O_{x,n}\in\calO_n$. By the property \xitemof{III:vDw-0} this sequence is well-defined and 
$\setof{O_{x,n}}{n\in\omega}$ is an open \nbhd\ basis for $x$. 
\qedofLemmaA
\qedskip}\fi 

A topological space $X=\pairof{X,\tau}$ is {\it collectionwise Hausdorff\/} if, for any 
discrete closed set $D\subseteq X$, there is a family $\calU=\setof{U_d}{d\in D}$ of 
pairwise disjoint open sets  \st\ the mapping $D\ni d\mapsto U_d\in\calU$ is 1-1 and $d\in U_d$ for all $d\in D$. 

A (pairwise disjoint) family $\calC$ of closed subsets of a space $X$ is said to be 
{\it discrete} if, for any $x\in X$, there is a \nbhd\ $U$ of $x$ \st\ $U$ 
intersects with at 
most one element of $\calC$. 

$X=\pairof{X,\tau}$ is {\it collectionwise normal\/} if, for any discrete family 
$\calC$ of closed sets, there is a family $\calU=\setof{U_C}{C\in\calC}$ of 
pairwise disjoint open sets \st\ $\calC\ni C\mapsto U_C\in\calU$ is 1-1 and 
$C\subseteq U_C$ for each $C\in\calC$. 

The following is immediate from the definitions above.
\begin{LemmaA}
  \Label{III:P-vDw-0} For a Hausdorff space $X=\pairof{X,\tau}$, if $X$ is 
  collectionwise normal then $X$ is collectionwise Hausdorff.\qed
\end{LemmaA}

The following well-known facts are also used in the proof of \Thmof{III:P-vDw-2} below. 
\begin{factA}
  \Label{III:P-vDw-1} \wassert{1}{\em\it Any metrizable space $X$ is collectionwise normal. 
  In particular,  by \Lemmaof{III:P-vDw-0}, any metrizable space $X$ is collectionwise Hausdorff.}\smallskip

  \wassert{2} {\bf(\textrm[Bin\cite{III:bib-bin}])} {\em\it A collectionwise normal Moore space is metrizable. }\qed
\end{factA}
}\fi 

As usual, $\boundingno$ denotes the bounding number which is defined as the 
minimal possible cardinality of a subset of $\fnsp{\omega}{\omega}$ which is 
unbounded \wrt\ $\leq^*$ (coordinate-wise comparison modulo finite). 

\ifextended\else Topological spaces constructed in Remark 12.6 in {[}van 
  Douwen\cite{III:bib-vdouwen}{]} witness the following theorem: 
\fi

\begin{Thm}[{[}van Douwen\cite{III:bib-vdouwen}{]}]
  \Label{III:P-vDw-2}
  There is a Moore space $X$ of cardinality $\boundingno$ \st\ $X$ is not 
  collectionwise Hausdorff (and hence non-metrizable\ifextended{\color{darkelectricblue}\ by {\em\rm\FactAof{III:P-vDw-1},\assertof{1}}}\fi) but all subspaces of $X$ of 
  cardinality $<\boundingno$ are metrizable. \ifextended\else\qed\fi
\end{Thm}
\ifextended{\small\color{darkelectricblue}
\prf Let $\calF\subseteq\fnsp{\omega}{\omega}$ and let
\begin{xitemize}
\xitema[III:vDw-1] $X_\calF=\calF\dotcup\omega\dotcup\calF\times\omega\times\omega$. 
\end{xitemize}

We define the topology on $X_\calF$ by declaring that 
\begin{xitemize}
\xitema[III:vDw-2] elements of $\calF\times\omega\times\omega$ are discrete;
\xitema[III:vDw-3] each $f\in\calF$ has a \nbhd\ basis consisting of sets of the form
  $O_{f,s}=\ssetof{f}\cup\ssetof{f}\times f\setminus s$ where $s$ is a finite 
  subset of $\omega\times\omega$; and  
\xitema[III:vDw-4] for $k\in\omega\ (\subseteq X_\calF)$,
  $U_{k,n}=\ssetof{k}\cup\calF\times\ssetof{k}\times(\omega\setminus n)$ for 
$n\in\omega$ form a \nbhd\ basis of $k\in\omega\subseteq X_\calF$. 
\end{xitemize}

\begin{Claim}\Label{Cl-III:vDw-0}
  $X_\calF$ is a normal Hausdorff space.
\end{Claim}
\prfofClaim
To show that $X_\calF$ is normal, one of the cases to be checked is that
any closed $F\subseteq X_\calF$ and 
$\pairof{f,m,n}\in\calF\times\omega\times\omega\setminus F$ can be separated by 
open sets. $\ssetof{\pairof{f,m,n}}$ is the minimal open \nbhd\ 
of $\pairof{f,m,n}$ by \xitemof{III:vDw-2}.

For $g\in F\cap\calF$, if $g\not=f$, then $O_{g,s}$ for any 
$s\in[\omega\times\omega]^{\LT\aleph_0}$ does not contain $\pairof{f,m,n}$ and 
hence disjoint from $\ssetof{\pairof{f,m,n}}$. 
If $g=f$, then letting $s=\ssetof{\pairof{m,n}}$, $O_{g,s}$ does not contain
$\pairof{f,m,n}$ and hence disjoint from the open set $\ssetof{\pairof{f,m,n}}$.

For $k\in F\cap\omega$, $U_{k, n+1}$  is disjoint from $\ssetof{\pairof{f,m,n}}$.

For $\pairof{f',m',n'}\in F\cap\calF\times\omega\times\omega$, Since
$\pairof{f',m',n'}\not=\pairof{f,m,n}$, the open 
\nbhd\ $\ssetof{\pairof{f',m',n'}}$ of $\pairof{f',m',n'}$ is disjoint from
$\ssetof{\pairof{f,m,n}}$.

Thus, we find an open superset of $F$ disjoint from $\ssetof{\pairof{f,m,n}}$ by 
taking union of all the open sets as above. 

The rest of the proof can be done similarly. 
\qedofClaim

\begin{Claim}
  \Label{Cl-III:vDw-1} $X_\calF$ is developable. 
\end{Claim}
\prfofClaim
Let $\seqof{s_n}{n\in\omega}$ be an increasing sequence of finite subsets of
$\omega\times\omega$ \st\ $\omega\times\omega=\bigcup_{n\in\omega}s_n$.

For each $n\in\omega$, let
\begin{xitemize}
\xitema[III:vDw-5] $\calO_n={}
  \begin{array}[t]{@{}l}
    \setof{\ssetof{\pairof{f,k,l}}}{\pairof{f,k,l}\in\calF\times\omega\times\omega}\\[\jot]
    \cup\ \ \setof{O_{f,s_n}}{f\in\calF}\ \ \cup\ \ \setof{U_{k,n}}{k\in\omega}.
  \end{array}
  $ 
\end{xitemize}

Then $\seqof{\calO_n}{n\in\omega}$ is a development of $X$.
\qedofClaim

\begin{Claim}\Label{Cl-III:vDw-2}
  If $\calF\subseteq\fnsp{\omega}{\omega}$ is unbounded (\wrt\ $\leq^*$), then 
$X_\calF$ is not collectionwise Hausdorff. In particular, $\calF$ is 
  non-metrizable.  
\end{Claim}
\prfofClaim
$D=\calF\cup\omega$ as a subset of $X_\calF$ is discrete and closed. We show 
that this set is a counter-example to the collectionwise Hausdorffness. Suppose, 
toward a contradiction, that $\calU$ is a family of pairwise disjoint open sets in
$X_\calF$ which separates elements of $D$. \Wolog, we may assume that elements of 
$\calU$ are of the form either $O_{f,s}$ or $U_{k,n}$. 

Let $\mapping{f^*}{\omega}{\omega}$ be defined by
\begin{xitemize}
\xitema[III:vDw-6] $f^*(k)=n$ if $U_{k,n}\in \calU$.
\end{xitemize}

$f^*$ is well-defined since $\calU$ is pairwise disjoint. Since $\calF$ is 
unbounded, there is $g^*\in\calF$ \st\ $g^*\not\leq^*f^*$. Thus $g^*(k)>f^*(k)$ 
for infinitely many $k\in\omega$. Let $s\in[\omega\times\omega]^{\LT\aleph_0}$ be 
\st\ $g^*\in O_{g^*,s}\in\calU$ and let
$k\in\omega\setminus\setof{m\in\omega}{\pairof{m,n}\in s\mbox{ for some }n\in\omega}$ 
be \st\ $g^*(k)>f^*(k)$. Then, since $U_{k,f^*(k)}\in\calU$, we have 
$\pairof{g^*,k, g^*(k)}\in O_{g^*,s}\cap U_{k,f^*(k)}\not=\emptyset$. This is a 
contradiction to the pairwise disjointness of $\calU$.

Thus $X_\calF$ is not collectionwise Hausdorff. 
$X_\calF$ is non-metrizable by \Factof{III:P-vDw-1},\,\assertof{1}. 
\qedofClaim

\begin{Claim}\Label{Cl-III:vDw-3}
  If $\calF\subseteq\fnsp{\omega}{\omega}$ is bounded, then $X_\calF$ is 
  collectionwise normal, and hence $X_\calF$ is metrizable by \Factof{III:P-vDw-1},\,\assertof{2}.
\end{Claim}
\prfofClaim
Suppose that $\calC$ is a discrete family of closed sets in $X_\calF$. Let 
$g^*\in\fnsp{\omega}{\omega}$ be \st\ $f<^*g^*$ for all $f\in\calF$. 
For each $f\in\calF$, let $s_f\in[\omega]^{\LT\aleph_0}$ be \st\
$f\restr\omega\setminus s_f<g^*\restr\omega\setminus s_f$ (point-wise).

Since $\calC$ is discrete, for each $x\in X$ there is a \nbhd\ $V_x$ of $x$ \st\ 
$V_x$ intersects at most one element of $\calC$.

For $C\in\calC$, let
\begin{xitemize}
\xitema[III:vDw-7] $O_C={} 
  \begin{array}[t]{@{}l}
    (C\cap\calF\times\omega\times\omega)\\[\jot]
    \cup\ \bigcup\setof{O_{f,f\restr s_f}\cap V_f}{f\in C\cap\calF}\\[\jot]
    \cup\ \bigcup\setof{U_{k,g^*(k)}\cap V_k}{k\in C\cap\omega}.
  \end{array}$
\end{xitemize}
Then $\calU=\setof{O_C}{C\in\calC}$ separates elements of $\calC$. 

This shows that $X_\calF$ is collectionwise normal. By 
\Factof{III:P-vDw-1},\,\assertof{2} and since $X_\calF$ is a Moore space 
by \Claimof{Cl-III:vDw-0} and \Claimof{Cl-III:vDw-1}, it follows that $X_\calF$ is 
metrizable.\qedofClaim
\qedskip

Now, suppose that $\calF\subseteq\fnsp{\omega}{\omega}$ is  unbounded  with
$\cardof{\calF}=\boundingno$.\\ By \Claimof{Cl-III:vDw-2}, $X_\calF$ is non-metrizable 
for any $\calF_0\subseteq\calF$ of cardinality $<\boundingno$, the subspace
$X_{\calF_0}$ of $X_\calF$ is metrizable by \Claimof{Cl-III:vDw-3}. Since subspaces 
of $X_\calF$ of the form $X_{\calF_0}$ for $\calF_0\in[\calF]^{\LT\boundingno}$ 
are cofinal in $[\calF]^{\LT\boundingno}$ and since any subspace of a metrizable 
space is metrizable, it follows that all subspaces of $X_\calF$ of cardinality
$<\boundingno$ are metrizable. 
\qedofThm
}\fi 

\ifextended\else
The proof of \Thmabove\ gives a construction of topological spaces $X_\calF$ for 
$\calF\subseteq\fnsp{\omega}{\omega}$. The construction is enough absolute and 
these spaces satisfy the properties that $X_\calF$ is a developable normal Hausdorff space,  $X_\calF$ is a subspace 
of $X_{\calF'}$ for $\calF\subseteq\calF'\subseteq\fnsp{\omega}{\omega}$, and 
$X_\calF$ is metrizable if and only if $\calF$ is bounded in $\fnsp{\omega}{\omega}$ (\wrt\ $\leq^*$). 
\fi

\begin{Cor}
  \Label{III:P-vDw-3} There is a non-metrizable Moore space $X=\pairof{X,\tau}$ \st\
  $\forces{\poP}{\check{X}\xmbox{ is metrizable}}$ for a $\sigma$-centered \po\ $\poP$. 
\end{Cor}
\prf
Let $X=X_\calF$ for an unbounded family $\subseteq\fnsp{\omega}{\omega}$. Let 
$\poP$ be the Hechler forcing then $\forces{\poP}{\check{\calF}\mbox{ is bounded}}$. 
Thus,  \ifextended{\small\color{darkelectricblue} by \Claimof{Cl-III:vDw-3},}\fi 
$\forces{\poP}{X_{\check{\calF}}\mbox{ is metrizable}}$. By the absoluteness of the definition of 
$X_\calF$, we have $\forces{\poP}{\check{X}=X_{\check{\calF}}}$. 
\qedofCor
\qedskip

The reflection number $\Refl_\HP$ of Hamburger's Hypothesis is defined by:
\begin{xitemize}
\xitem[III:vDw-8] 
  $\Refl_\HP={}\left\{{}
  \begin{array}{@{}ll}
    \mbox{the minimal cardinal }\kappa\mbox{ \st,}\\
    \mbox{for any first countable 
      non-metrizable}\\
    \mbox{topological space }X,\mbox{ there is a non-}\\
    \mbox{metrizable subspace }Y\mbox{ of }X\mbox{ of}\\
    \mbox{cardinality }\LT\kappa;&\hspace{-1.2em}\mbox{if such }\kappa\mbox{ exists,}\\[2\jot]
    \infty; &\hspace{-1.2em}\mbox{otherwise.}
  \end{array}\right.$
\end{xitemize}

\begin{Lemma}
  \Label{III:P-vDw-4}
  \wassertof{1} $\boundingno<\Refl_\HP\leq\infty$.\smallskip

  \wassert{2} $\Refl_\HP=\infty$ is consistent.\smallskip

  \wassert{3} {\em(\bf[Bagaria and Magidor\cite{III:bib-bm}])}\ \,$\Refl_\HP\leq$ the least $\omega_1$-strongly compact cardinal (if it 
  exists). 
\end{Lemma}
\prf
\assertof{1}: By \Thmof{III:P-vDw-2}. 
\ifextended\else\smallskip\fi

\ifextended{\small\color{darkelectricblue} For $\aleph_1<\Refl_\HP$, we have more direct examples:
$\omega_1$ with the order topology or $E^\kappa_\omega$ for any cardinal of 
uncountable cofinality (also with the order topology) are among the examples 
showing the inequality $\aleph_1<\Refl_\HP$.}\smallskip\fi

\assertof{2}: This holds if $\Box_\kappa$ holds for cofinally many $\kappa$ (in
$\Card$) --- 
actually $\ADS^-(\kappa)$ for class many regular uncountable $\kappa$ is enough 
(see Proposition 6.3 in [Fuchino, Juh\'asz et al.\cite{III:bib-erice}])). \smallskip

\assertof{3}: \ifextended{\small\color{darkelectricblue} Suppose that $(X,\calO)$ is a first countable topological space 
\st\ all subspaces $Y\in[X]^{<\kappa}$ are metrizable. For each $x\in X$, let 
$\setof{O_{x,n}}{n\in\omega}$ be an open \nbhd\ base of $x$. 

Let $T$ be the $\calL_{\omega_1,\omega}$ theory in the language with the binary 
relation symbols $O_n(x,y)$ for all $n\in\omega$ coding ``$y\in O_{x,n}$'' and the 
binary symbols $d_q(x,y)$ for all $q\in\rationals_{\geq0}$ which should code
``$d(x,y)\leq q$'':
\begin{xitemize}
\xitema[x-1] {$T=\begin{array}[t]{@{}l}
    \setof{O_n(c_a,c_b)}{a,b\in X,\,b\in O_{a,n}}\\[\jot]
    \cup\ \setof{\neg O_n(c_a,c_b)}{a,b\in X,\,b\not\in O_{a,n}}\\[\jot]
    \cup\ \setof{\forall x\forall 
      y\ (d_q(x,y)\rightarrow d_q(y,x))}{q\in\rationals_{\geq0}}\\[\jot]
    \cup\ \setof{\forall x\forall y\ (d_q(x,y)\rightarrow d_{q'}(x,y))}{
      q,q'\in\rationals_{\geq 0},\,q\leq q'}\\[\jot]
    \cup\ \ssetof{\forall x \forall y\ (d_0(x,y)\rightarrow x\equiv y)}\\[\jot]
    \rlap{\(\cup\ \setof{\forall x\forall y\forall z\ 
      (d_q(x,y)\land d_{q'}(y,z)\rightarrow d_{q+q'}(x,z))}{q,q'\in\rationals_{\geq0}}\)}
    \\[\jot]
    \cup\ \setof{\forall x\llor_{q\in\rationals_{>0}} 
    \forall y\ (d_q(x,y)\rightarrow O_n(x,y))}{n\in\omega}\\[\jot]
    \cup\ \ssetof{\forall x\lland_{q\in\rationals_{>0}}\llor_{n\in\omega}
      \forall y\ (O_n(x,y)\rightarrow d_q(x,y))}
  \end{array}$}
\end{xitemize}
Clearly all $T'\in[T]^{<\kappa}$ are satisfiable. 

Since $\kappa$ 
is $\omega_1$-strongly compact, it follows that $T$ is also satisfiable. Let $M$ 
be a model of $T$.  
Then $\mapping{d}{X^2}{\reals}$ defined by
\begin{xitemize}
\xitema[x-2] 
  $d(a,b)=\inf\setof{q\in\rationals}{M\models d_q(c_a,c_b)}$ for $a$,
  $b\in X$ 
\end{xitemize}
is a metric on $X$ generating the topology of $(X,\calO)$.
}\else 
See: [Bagaria and Magidor\cite{III:bib-bm}]. 
\fi 
\qedofLemma

\section{Preservation and non-preservation of stationarity of subsets of $\Pkl{}{}$ }
\Label{III:preserv}

In the following, we show that the closedness of \pos\ cannot be used to establish 
reflection principles concerning the stationarity of subsets 
of $\Pkl{}{}$ for $\kappa>\aleph_1$ in the generic extensions. At least, not in a 
straight-forward generalization of the usage of $\sigma$-closed \pos\ in a forcing argument 
to obtain reflection properties on stationarity of subsets of $\Pkl{\aleph_1}{\lambda}$ in 
the generic extensions.

Actually, the examples of  
preservation and non-preservation of 
stationarity of subsets of $\Pkl{}{}$ in this section explain, why we need a mixed support 
iteration plus one further step with chain 
condition in connection with the following \Lemmaof{III:L-preserv-a} to establish 
(some of the) results in \sectionof{III:refl} but not in a much simpler way. 

It is well-known that ccc \pos\ and $\sigma$-closed \pos\ are proper. This means 
that such \pos\ preserve stationarity of subsets of $\Pkl{\aleph_1}{\lambda}$ for 
any uncountable $\lambda$. For \pos\ with $\kappa$-cc for regular cardinal $\kappa>\aleph_1$ 
we still have a corresponding lemma:
\begin{Lemma}
  \Label{III:L-preserv-a}
  Suppose that $\kappa$ is a regular uncountable cardinal and $\lambda\geq\kappa$. If
  $S\subseteq\Pkl{}{}$ is stationary and $\poP$ is a $\kappa$-cc \po, then we have
  $\forces{\poP}{\check{S}\mbox{ is stationary}}$. 
\end{Lemma}
\ifextended
{\small\color{darkelectricblue}
  \prf
  Suppose that $\utilde{\calC}$ is a $\poP$-name with
  $\forces{\poP}{\utilde{\calC}\mbox{ is a club in }\Pkl{\kappa}{\lambda}}$.

  In $\uniV$, 
  let
  $\calC=\setof{C\in\Pkl{\kappa}{\lambda}}{\forces{\poP}{\check{C}\varin\utilde{\calC}}}$. Then
  $\calC$ is club by the $\kappa$-cc of $\poP$. 
  Hence $\calS\cap\calC\not=\emptyset$. Since
  $\forces{\poP}{\check{\calC}\subseteq\utilde{\calC}}$, it follows that
  $\forces{\poP}{\check{\calS}\cap\utilde{\calC}\not\equiv\emptyset}$. 
  \qedofLemma\qedskip
}
\else
The standard proof for the case $\kappa=\aleph_1$ also works for this general 
lemma.
\fi

In contrast, $\kappa$-closed \po\ can destroy stationarity of ground model 
stationary set $\subseteq\Pkl{}{}$ if $\kappa>\aleph_1$. This makes consistency 
proofs of stationary reflection of stationary subsets of $\Pkl{}{}$ for
$\kappa\geq\aleph_2$ more involved. In the following, we shall examine situations 
where the stationarity of some subset of $\Pkl{}{}$ for $\kappa\geq\aleph_2$ is not 
preserved by a standard $\kappa$-closed \po. 

For cardinal $\kappa$ and a regular cardinal $\nu<\kappa$ we denote
\begin{xitemize}
\xitem[III:refl-1-0] 
  $E^\kappa_\nu=\setof{\alpha\in\kappa}{\cf(\alpha)=\nu}$. 
\end{xitemize}

The following Lemma is used for 
our first example of non-preservation of stationarity in 
\Propof{III:P-preserv-0}.

\begin{Lemma}
  \Label{III:L-preserv-0}Suppose that $\kappa$ is a regular cardinal with 
  $\kappa\geq\aleph_2$ and $X\supseteq\kappa^+$. Then, for any distinct 
  regular $\nu$, $\mu<\kappa$,
  \begin{xitemize}
  \xitem[III:preserv-0]
    $S=\setof{x\in\Pkl{}{X}}{\kappa\cap x\in E^\kappa_\nu,\,\cf(\sup(\kappa^+\cap x))=\mu}$
  \end{xitemize}
  is stationary in $\Pkl{}{X}$. 
\end{Lemma}
\prf Suppose that $C\subseteq\Pkl{}{X}$ is a club. Let 
$\mapping{f}{[X]^{<\aleph_0}}{X}$ be \st\
$C\ell^*(f)=\setof{x\in\Pkl{}{X}}{x\cap\kappa\in\kappa,\,x\mbox{ is closed \wrt\ }f}\subseteq C$.

Let $X_0\subseteq X$ be \st\ $X_0$ is closed \wrt\ $f$ and $\kappa^+\cap X_0\in E^{\kappa^+}_\mu$.
Let $\delta=\kappa^+\cap X_0$.

Let $\seqof{x_\xi}{\xi<\nu}$ be a continuously increasing sequence 
in $\Pkl{}{X_0}$ \st\
\begin{xitemize}
\xitem[III:preserv-1] $x_\xi$ is closed \wrt\ $f$ for all $\xi<\nu$;
\xitem[III:preserv-2] $\sup(\kappa^+\cap x_0)=\delta$; and 
\xitem[III:preserv-3] $\sup(\kappa\cap x_\xi)+1\subseteq x_{\xi+1}$ for all $\xi<\nu$.
\end{xitemize}
Note that this construction is possible since $\kappa$ is regular 
and $\nu$, $\mu<\kappa$. 

Let $x=\bigcup_{\xi<\nu}x_\xi$. Then
\begin{xitemize}
\xitem[III:preserv-4] $x$ is closed \wrt\ $f$;\qquad (by \xitemof{III:preserv-1})
\xitem[III:preserv-5] $\sup(\kappa^+\cap x)=\delta$;\qquad (by \xitemof{III:preserv-2})
\xitem[III:preserv-6] $\kappa\cap x\in\kappa$ and $\cf(\kappa\cap x)=\nu$.\qquad 
  (by \xitemof{III:preserv-3})
\end{xitemize}
$x\in S$ by \xitemof{III:preserv-5} and \xitemof{III:preserv-6}.
$x\in C\ell^*(f)$ by \xitemof{III:preserv-4} and 
\xitemof{III:preserv-6}. Thus we have
$\emptyset\not=S\cap C\ell^*(f)\subseteq S\cap C$. 
\qedofLemma
\qedskip

For a regular cardinal, $\Add(\kappa)$ denotes the set 
$\fnsp{\kappa>}{2}$ with the reverse inclusion. We denote with
$\Col(\kappa,\kappa^+)$ the set $\fnsp{\kappa>}{\kappa^+}$ with the reverse 
inclusion. $\Add(\kappa)$ and $\Col(\kappa,\kappa^+)$ are forcing equivalent to
$\Fn(\kappa,2,\kappa)$ and $\Fn(\kappa,\kappa^+,\kappa)$ in Kunen's 
notation in [Kunen\cite{III:bib-kunen}], respectively.
The \pos\ isomorphic to latter two \pos\ are also denoted as $\Col(\kappa,\ssetof{\kappa})$ and
$\Col(\kappa,\ssetof{\kappa^+})$ respectively, 
in the notation of [Kanamori\cite{III:bib-kanamori}]. Both of the \pos\ are $\kappa$-closed. 
$\Add(\kappa)$ adds a new subset of $\kappa$ while $\kappa^+$ is preserved if
$\kappa^{<\kappa}=\kappa$. $\Col(\kappa,\kappa^+)$ collapses $\kappa^+$ and makes 
it of cardinality and cofinality $\kappa$.

\begin{Prop}
  \Label{III:P-preserv-0}
  Suppose that $\kappa$ is a regular cardinal $\geq\aleph_2$ and
  $X\supseteq\kappa^+$. Then, there is a stationary $S\subseteq\Pkl{}{X}$,
  \st\ $\forces{\Col(\kappa, \kappa^+)}{\check{S}\xmbox{ is not stationary in }\Pkl{}{X}}$. 
\end{Prop}
\prf In $\uniV$, let
\begin{xitemize}
\xitem[III:preserv-a-0] 
  $S=\setof{x\in\Pkl{}{X}}{{}
  \begin{array}[t]{@{}l}
    x\cap\kappa\mbox{ and }\sup(x\cap\kappa^+)\mbox{ are limit ordinals, and}\\
    \cf(x\cap\kappa)\not=\cf(\sup(x\cap\kappa^+))}.
  \end{array}
$
\end{xitemize}
$S$ is a stationary subset of $\Pkl{}{X}$ by \Lemmaof{III:L-preserv-0}. We show 
that $\Col(\kappa,\kappa^+)$ forces that $S$ is not stationary.

Suppose that $\genG$ is a $(\uniV,\Col(\kappa,\kappa^+))$-generic filter. 
Note that, by $\LT\kappa$-closedness, $\Col(\kappa,\kappa^+)$ does not add any 
new sets of size $\LT\kappa$. Thus 
$\Pkl{}{X}^\uniV=\Pkl{}{X}^{\uniV[\genG]}$, all cofinalities $\LT\kappa$ are 
preserved in the generic extension $\uniV[\genG]$, and $\cf(\mu)=\kappa$ in
$\uniV[\genG]$ for $\mu=(\kappa^+)^{\uniV}$.

In $\uniV[\genG]$, 
let $\seqof{\gamma_\alpha}{\alpha<\kappa}$ be a continuously increasing sequence 
of ordinals cofinal in the ordinal $\mu$. Let
\begin{xitemize}
\xitem[III:preserv-a-1] 
  $C=\setof{x\in\Pkl{}{X}}{{}
  \begin{array}[t]{@{}l}
    x\cap\kappa\mbox{ and }\sup(x\cap\mu)\mbox{ are limit ordinals, and}\\
    \sup(x\cap\mu)=\gamma_{x\cap\kappa}}.
  \end{array}$
\end{xitemize}

Then $C$ is a club in $\Pkl{}{X}$ and $C\cap S=\emptyset$.
\qedofProp
\qedskip

In \Propof{III:P-preserv-0}, the crucial fact which made the set 
$S$ non-stationary in the generic extension was that the cardinal $\kappa^+$ is collapsed 
to be an ordinal of cofinality $\kappa$. However, stationarity of $\Pkl{}{}$ can be also destroyed 
by a $\LT\kappa$-closed forcing without collapsing cardinals: 

\begin{Prop}
  \Label{III:P-preserv-1}
  Suppose that $\kappa$ is a supercompact and $\cardof{X}\geq2^\kappa$. Then 
  there is a stationary $S\subseteq\Pkl{}{X}$ \st\
  $\forces{\Add(\kappa)}{\check{S}\xmbox{ is not a stationary subset of }\Pkl{}{X}}$.
\end{Prop}

Note that $\cardof{\Add(\kappa)}=\kappa$ since $\kappa$ is inaccessible and hence 
$\Add(\kappa)$ is $\kappa^+$-cc. Thus $\Add(\kappa)$ here preserves cardinals and 
cofinality. \qedskip

\noindent
\prf Let $\lambda=\cardof{X}$. \Wolog, we may assume that $X=\lambda$. 
In $\uniV$, let $\vec{B}=\seqof{B_\alpha}{\alpha<\lambda}$ be an enumeration of
$\psof{\kappa}$ and let
\begin{xitemize}
\xitem[III:preserv-7] 
  $S=\setof{x\in\Pkl{}{}}{{}
  \begin{array}[t]{@{}l}
    \mbox{\assertof{a}\ \ } \kappa\cap x\in\kappa\mbox{, and}\\[\jot]
    \mbox{\assertof{b}\ \ } \setof{B_\alpha\cap(x\cap\kappa)}{\alpha\in x}=\psof{\kappa\cap x}}.
  \end{array}$
\end{xitemize}

\begin{Claim}
  \Label{III:Cl-preserv-0} $S\in\calU$ for any normal ultrafilter $\calU$ over $\Pkl{}{}$.
\end{Claim}
\prfofClaim
Suppose that $\calU$ is a normal ultrafilter over $\Pkl{}{}$. It is enough to 
show that $j_\calU\imageof\lambda\in j_\calU(S)$ where
$\elembed{j_\calU}{\uniV}{M}$ is the elementary embedding induced by $\calU$. 

We have 
\begin{xitemize}
\xitem[III:preserv-7-0] 
  $\left(j_\calU\imageof\lambda\right)\cap j_\calU(\kappa)
  =\setof{j_\calU(\alpha)}{\alpha<\lambda,\,j_\calU(\alpha)<j_\calU(\kappa)}$\\[\jot]
  $=\setof{j_\calU(\alpha)}{\alpha<\kappa}=\setof{\alpha}{\alpha\in\kappa}
  =\kappa\in j_\calU(\kappa)$.
\end{xitemize}
For $\beta\in j_\calU\imageof\lambda$ with $\beta=j_\calU(\alpha)$ for
$\alpha\in\lambda$,
$j_\calU(\vec{B})(\beta)\cap (j_\calU\imageof\lambda\cap j_\calU(\kappa))
=j_\calU(B_\alpha)\cap\kappa=B_\alpha$. 

Thus we have
\begin{xitemize}
\xitem[III:preserv-7-1] 
  $\setof{j_\calU(\vec{B})(\beta)\cap (j_\calU\imageof\lambda\cap j_\calU(\kappa))}{\beta\in j\imageof\lambda}
  =\setof{B_\alpha}{\alpha<\lambda}$\\[\jot]
  $=\psof{\kappa}=\psof{j_\calU\imageof\lambda\cap j(\kappa)}$. 
\end{xitemize}
By elementarity, \xitemof{III:preserv-7-0} and \xitemof{III:preserv-7-1} imply $j_\calU\imageof\lambda\in j(S)$. 
\qedofClaim
\qedskip

Since any normal filter over $\Pkl{}{}$ contains all club sets and hence it 
consists of stationary sets, and since there {\it are} normal ultrafilters over
$\Pkl{}{}$ because $\kappa$ is 
supercompact,
we conclude that $S$ is a stationary subset of $\Pkl{}{}$. 

Thus the next Claim shows that $S$ is as desired:
\begin{Claim}
  $\forces{\Add(\kappa)}{S\mbox{ is not stationary in }\Pkl{}{}}$. 
\end{Claim}
\prfofClaim
Suppose that $\genG$ is a $(\uniV,\Add(\kappa))$-generic filter. 
In $\uniV[\genG]$, we have $\mapping{\bigcup\genG}{\kappa}{2}$. 
Let $A=(\bigcup\genG)^{-1}\imageof\ssetof{1}$. By genericity, $A$ is a new subset 
of $\kappa$. Let $\mapping{F}{\lambda}{\kappa}$ be defined by
$F(\alpha)=\min(B_\alpha\bigtriangleup A)$ for $\alpha<\lambda$. $F$ is well-defined 
since $B_\alpha\not=A$ for all $\alpha<\lambda$. We show that
$S\cap C_F=\emptyset$ where $C_F$ is the club set defined by 
\begin{xitemize}
\xitem[III:preserv-8] 
  $C_F=\setof{x\in\Pkl{}{}}{x\xmbox{ is closed \wrt\ }F}$.
\end{xitemize}

Suppose that $x\in S$. 
By \assertof{b} in \xitemof{III:preserv-7} and since $A\cap x\in\psof{x\cap\kappa}^\uniV$, there is 
an $\alpha^*\in x$ \st\ $B_{\alpha^*}\cap x=A\cap x$. But this implies that
$F(\alpha^*)=min(B_\alpha\bigtriangleup A)\not\in x$. Thus, $x$ is not closed \wrt\ 
$F$ and $x\not\in C_F$. 
\qedofClaim\\
\qedofProp
\qedskip

The non-preservation of stationarity of subsets of $\Pkl{}{}$ along the line of the results 
above is further studied in [Sakai\cite{III:bib-sakai}].

\section{Two dimensional Laver-generic large cardinals}\Label{III:two}
For properties $\gmP$ and $\gmQ$ of \pos, a cardinal $\kappa$ is 
{\it Laver-generically supercompact for $(\gmP,\gmQ)$} if, for any \po\ $\poP$ with
$\poP\models\gmP$, $(\uniV,\poP)$-generic $\genG$, and a cardinal $\lambda$, 
there are a $\poP$-name $\utpoQ$\vspace{-0.2\smallskipamount} of a \po\ with $\forces{\poP}{\utpoQ\models\gmQ}$, and
a $(\uniV,\poP\ast\utpoQ)$-generic $\genH$ with $i\imageof\genG\subseteq\genH$ 
where $\mapping{i}{\poP}{\poP\ast\utpoQ}$ is the canonical complete embedding,  
\st\ there 
are $j,M\subseteq\uniV[\genH]$ with
\begin{xitemize}
\xitem[III:gen-large-1-0] $M$ is a transitive class in $\uniV[\genH]$; 
\xitem[III:gen-large-2] $\elembed{j}{\uniV}{M}$;
\xitem[III:gen-large-3] 
  $\crit(j)=\kappa$ and $j(\kappa)>\lambda$;
\xitem[III:gen-large-4] 
  $\poP$, $\genH\in M$, and 
\xitem[III:gen-large-6] 
  $j\imageof\lambda\in M$.
\end{xitemize}
$\kappa$ is {\it strongly Laver-generically supercompact for $(\gmP,\gmQ)$} if $M$ in the 
definition of the Laver-generic supercompactness for $(\gmP,\gmQ)$ additionally satisfies
\begin{xitemize}
\xitem[III:gen-large-6-0] 
  $([M]^{\aleph_0})^{\uniV[\genH]}\subseteq M$.
\end{xitemize}

\ifextended{\small\color{darkelectricblue}
  If $\setof{\poP}{\poP\models\gmP}$ contains only trivial \pos, then the (strongly) Laver-generic 
  supercompactness for $(\gmP,\gmQ)$ coinsides with the (strongly) generic supercompactness 
  by \pos\ satisfying $\gmQ$. If $\gmP$ and $\gmQ$ are equivalent, and $\gmP$ is iterable, that  
  is, if for every $\poP\models\gmP$ and $\poP$-name $\utpoQ$ with
  $\forces{\poP}{\utpoQ\models\gmP}$, we have $\poP\ast\utpoQ\models\gmP$, then the (strongly) Laver-generic 
  supercompactness for $(\gmP,\gmP)$ is closely related to the (strongly) Laver generic 
  supercompactness for $\poP$ in the sense of \cite{III:bib-II} but may not be exactly the 
  same notion. 
}\else 
If $\gmP$ solely consists of trivial forcing, we shall say ``{\it strongly 
  Laver-generically supercompact for 
$\gmQ$}'' instead of ``strongly Laver-generically supercompact for 
$(\gmP,\gmQ)$''. 
\fi 

In the following, both of the properties $\gmP$ and $\gmQ$ considered in connection with the 
Laver-generic supercompactness imply the properness of the \po. In such a case, the 
model of the Laver-generic supercompactness constructed by forcing starting from a 
supercompact cardinal usually satisfies 
this strong version of Laver-generic supercompactness as well. This is because of the 
following well-known fact:
\begin{Lemma}
  \Label{III:P-gen-large-a-0} Suppose that $M\subseteq\uniV$ is an inner model with 
  \begin{xitemize}
  \xitem[III:gen-large-6-1] 
    $[M]^{\aleph_0}\subseteq M$. 
  \end{xitemize}
  If $\poP\in M$ is proper and $\genG$ is a $(\uniV,\poP)$-generic filter, then we have
  \begin{xitemize}
  \xitem[III:gen-large-6-2] 
    $([M[\genG]]^{\aleph_0})^{\uniV[\genG]}\subseteq M[\genG]$. 
  \end{xitemize}
\end{Lemma}
\prf
Let $a\in([M[\genG]]^{\aleph_0})^{\uniV[\genG]}$ and $\uta$ be a $\poP$-name of $a$. In
$\uniV$, let $\theta$ be a sufficiently large regular cardinal, and 
let $N\prec\calH(\theta)$ and $\condp\in\genG$ be \st\
\begin{xitemize}
\xitem[III:gen-large-6-3] 
  $\cardof{N}=\aleph_0$;
\xitem[III:gen-large-6-4] $\uta, \poP\in N$;
\xitem[III:gen-large-6-5] for any maximal antichain $A\subseteq\poP$ with $A\in N$, 
  $A\cap N$ is predense below $\condp$.
\end{xitemize}

Let $\utf\in N$ be a $\poP$-name \st\\
$\forces{\poP}{\mapping{\utf}{\omega}{\uta}\xmbox{ is a surjection}}$. For each $n\in\omega$, 
let $A_n\in M$ be a maximal antichain \st, for each $\condr\in A_n$, there is 
a $\poP$-name $\utb_{\,n,\condr}\in M$ 
\st\ $\condr\forces{\poP}{\utf(n)\equiv\utb_{\,n,\condr}}$. Note that we have
$\utb_{\,n,\condr}\in M\cap N$ if $\condr\in N$, 
since $\utb_{\,n,\condr}$ is uniquely determined for each $\condr\in A_n$. 

Let $\uta^*=\setof{\pairof{\utb_{\,n,\condr},\condr}}{n\in\omega,\condr\in A_n\cap N}$. 
Then $\uta^*\in M$ by \xitemof{III:gen-large-6-1} and $\uta^*[\genG]=\uta[\genG]$. Thus,
$a=\uta[\genG]\in M[\genG]$. 
\qedofLemma
\qedskip

The condition \xitemof{III:gen-large-6-0} can be even replaced with 
\begin{xitemize}
\xitemciteb[III:gen-large-6-0]{$'$} 
    $([M]^{j(\kappa)})^{\uniV[\genH]}\subseteq M$, 
\end{xitemize}
if we consider Laver-generic superhugeness instead of Laver-generic supercompactness.

This 
can be seen by means of the following:
\begin{Lemma}
  \Label{III:P-gen-large-a-1}
  Suppose that $M$ is an inner model of\/ $\uniV$ with 
  \begin{xitemize}
  \xitem[III:gen-large-6-6] 
    $V\modelof{[M]^\mu\subseteq M}$ 
  \end{xitemize}
  for a regular $\mu$. 
  If\/
  $\poP\in M$ is $\mu^+$-cc, then, for any $(\uniV,\poP)$-generic $\genG$, we have
  \begin{xitemize}
  \xitem[III:gen-large-6-7] 
    $([M[\genG]]^{\mu})^{\uniV[\genG]}\subseteq M[\genG]$. 
  \end{xitemize}
\end{Lemma}
  \prf \ifextended{\color{darkelectricblue}
    Note that $\poP\subseteq M$ since $M$ is transitive. }
  \fi
  Suppose $g\in\left(\fnsp{\mu}{M[\genG]}\right)^{\uniV[\genG]}$. We show that
  $g\in M[\genG]$. Let $\utg$ be a $\poP$-name of $g$. For each $\xi<\mu$, there is a maximal 
  pairwise incompatible $A_\xi\subseteq\poP$ \st, for each
  $\condp\in A_\xi$, there is a $\poP$-name $\uta_{\xi,\condp}\in M$ \st\
  $\condp\forces{\poP}{\utg(\xi)\equiv\uta_{\xi,\condp}}$. 
  By the $\mu^+$-cc of $\poP$,
  we have $\cardof{A_\xi}\leq\mu$ and hence $A_\xi\in M$ by \xitemof{III:gen-large-6-6}.

  Let 
  \begin{xitemize}
  \xitem[III:gen-large-6-8]
    $\uta_\xi=\setof{\pairof{\utb,\condq}}{{}
    \begin{array}[t]{@{}l}
      \condq\leq_\poP\condp\mbox{ for some }\condp\in A_\xi,\,\utb\mbox{ is a canonical}\\
      \poP\mbox{-name with }\condq\forces{\poP}{\utb\varin\uta_{\xi,\condp}}}.
    \end{array}$
  \end{xitemize}
  $\uta_\xi\in M$ since it is definable from $\leq\mu$ many parameters from $M$ and by 
  \xitemof{III:gen-large-6-6}. It is also clear by the definition above that $\forces{\poP}{\utg(\xi)\equiv\uta_\xi}$. 
  Let 
  \begin{xitemize}
  \xitem[III:gen-large-6-9] 
    $\utg^*=\setof{\pairof{\pairof{\check{\xi},\uta_\xi}^\bullet_{\poP},\bbone_\poP}}{\xi<\mu}$.
  \end{xitemize}
  Then $\utg^*\in M$ by \xitemof{III:gen-large-6-6} and $M[\genG]\ni\utg^*[\genG]=\utg[\genG]=g$. 
  \qedofLemma\qedskip

\Lemmaof{III:P-gen-large-a-1} is used with a generic superhuge 
cardinal to produce the strong generic superhugeness. This can be seen as follows:

Suppose that $\kappa$ is a superhuge cardinal and
$\vec{\poP}=\seqof{\poP_\alpha,\utpoQ_\beta}{\alpha\leq\kappa,\beta<\kappa}$ be an iteration with 
\begin{xitemize}
\xitem[III:gen-large-6-10] 
  $\cardof{\poP_\alpha}<\kappa$ for all $\alpha<\kappa$. 
\end{xitemize}
Let $\genG_\kappa$ be a
$(\uniV,\poP_\kappa)$-generic filter. 

For a given cardinal $\lambda$, let 
$\elembed{j}{\uniV}{M}$  be an elementary embedding into an inner model $M$ of $\uniV$ \st\
$\crit(j)=\kappa$, $j(\kappa)>\lambda$ and $[M]^{j(\kappa)}\subseteq M$. 
Let $\poP^*=j(\poP_\kappa)$. By elementarity $\poP^*$ is the $j(\kappa)$th iterand of the 
iteration $j(\vec{\poP})$ in $M$. There is the canonical complete embedding
$\combed{i}{\poP_\kappa}{\poP^*}$. Let $\genG^*$ be $(\uniV,\poP^*)$-generic filter with
\begin{xitemize}
\xitem[III:gen-large-6-11] 
  $i\imageof\genG_\kappa\subseteq\genG^*$. 
\end{xitemize}
By \xitemof{III:gen-large-6-10} and elementarity, 
we have $M\modelof{\cardof{\poP^*}\leq j(\kappa)}$. Hence $\cardof{\poP^*}\leq j(\kappa)$ 
and thus $\poP^*$  has $j(\kappa)^+$-cc. By \xitemof{III:gen-large-6-11}, $j$ can be 
lifted to
\begin{xitemize}
\xitem[III:gen-large-6-12] 
  $\elembed{\tilde{j}}{V[\genG_\kappa]}{M[\genG^*]}$; $\uta[\genG_\kappa]\mapsto j(\uta)[\genG^*]$
\end{xitemize}
for $\poP_\kappa$-names $\uta$. Now by \Lemmaof{III:P-gen-large-a-1}, we have
$([M[\genG^*]]^{j(\kappa)})^{\uniV[\genG^*]}\subseteq M[\genG^*]$.

Based on the observations above, we define $\kappa$ to be {\it strongly Laver-generically superhuge for 
  the pair of properties $(\gmP,\gmQ)$} if, for any \po\ $\poP$ with
$\poP\models\gmP$ and $(\uniV,\poP)$-generic $\genG$, there are a $\poP$-name
$\utpoQ$\vspace{-0.2\smallskipamount} of a \po\ with $\forces{\poP}{\utpoQ\models\gmQ}$ and
a $(\uniV,\poP\ast\utpoQ)$-generic $\genH$ with $i\imageof\genG\subseteq\genH$ 
where $\mapping{i}{\poP}{\poP\ast\utpoQ}$ is the canonical complete embedding,  
\st\ there 
are $j,M\subseteq\uniV[\genH]$ with \xitemof{III:gen-large-1-0} $\sim$ 
\xitemof{III:gen-large-4} and \xitembof{III:gen-large-6-0}{$'$}. 

\ifextended{\color{darkelectricblue}
If $\setof{\poP}{\poP\models\gmP}$ contains only trivial \pos, }
\else
If $\gmP$ solely consists of trivial forcing, 
\fi 
we shall say ``{\it strongly superhuge for 
$\gmQ$}'' instead of ``strongly Laver-generically superhuge for 
$(\gmP,\gmQ)$''.

The following is trivial. 
\begin{Lemma}
  \Label{III:P-gen-large-0}Suppose that\/ $\poP\models\gmP_0$ implies $\poP\models\gmP_1$ and\/ 
  $\poP\models\gmQ_1$ implies $\poP\models\gmQ_0$ for all \pos\ $\poP$ (i.e. these 
  implications are theorems in \ZFC). If $\kappa$ is 
  (strongly) Laver-generically supercompact/superhuge for 
  $(\gmP_1,\gmQ_1)$, then $\kappa$ is  (strongly) Laver-generically supercompact/huge for
  $(\gmP_0,\gmQ_0)$. \qed 
\end{Lemma}

We call a property $\gmP$ of \pos\ {\it iterable} if we can prove in \ZFC\ that 
\begin{xitemize}
\xitem[] 
  $\poP\ast\utpoQ\models\gmP$ for any \po\ $\poP$ with 
  $\poP\models\gmP$ and $\poP$-name $\utpoQ$ of a \po\ with
  $\forces{\poP}{\utpoQ\models\gmP}$. 
\end{xitemize}

\begin{Prop}
  \Label{III:P-gen-large-1}Suppose that $\gmP$ is the property ``forcing equivalent to a 
  \po\ of the form $\Col(\kappa,\mu)$ for 
  some $\mu$''{\footnote{In this section, we are back to Kanamori's notation of collapsing 
    \pos. $\Col(\kappa,\lambda)$ for an inaccessible $\lambda$ is thus the \po\ 
    collapsing all cardinals strictly between $\kappa$ and $\lambda$ by conditions of size
    $\LT\kappa$.}}, $\gmQ$ is iterable and we can prove (\/in {\em\ZFC}) that 
  \begin{xitemize}
  \xitem[III:gen-large-6-13] 
  	$\forall P\,(P\mbox{ is a }\sigma\mbox{-directed closed \po}\rightarrow P\models\gmQ)$. 
  \end{xitemize}
  If $\kappa$ is strongly Laver-generically supercompact for $(\gmP,\gmQ)$, 
  then, for any $\mu\geq\kappa$, and for $\poP=\Col(\kappa,\mu)$, we have
  \begin{xitemize}
  \xitem[III:gen-large-7] 
    $\forces{\poP}{\kappa\mbox{ is generically supercompact by \pos\ satisfying }\gmQ}$.
  \end{xitemize}
\end{Prop}
\prf The proof is a typical application of the master condition argument. 

\ifextended{\color{darkelectricblue} Let $\poP=\Col(\kappa,\mu)$ for some $\mu\geq\kappa$. }\fi
Note that $\poP\models\gmP$. Let $\genG$ be an 
arbitrary $(\uniV,\poP)$-generic filter. We have to show that
\ifextended{\color{darkelectricblue}
  \begin{xitemize}
  \xitema[III:gen-large-7-a] 
    $\uniV[\genG]\models\qbox{\kappa\mbox{ is generic supercompact by \pos\ with }\gmQ}$.
  \end{xitemize}
}\else
$\uniV[\genG]\models\qbox{\kappa\mbox{ is generic supercompact by \pos\ with }\gmQ}$. \fi

Let \ifextended{\color{darkelectricblue} $\theta\geq\kappa$ be arbitrary and let
  $\lambda=\max\ssetof{\theta,\mu^{\LT\kappa}}$. Let }\else
$\lambda\geq\theta$, $\mu^{\LT\kappa}$ and let \fi$\utpoQ$ be a $\poP$-name with
$\forces{\poP}{\utpoQ\mbox{ is a \po\ with }\gmQ}$ \st\ there is a
$(\uniV,\poP\ast\utpoQ)$-generic filter $\genH$ \st\ $i\imageof\genG\subseteq\genH$ for the 
canonical complete embedding $\mapping{i}{\poP}{\poP\ast\utpoQ}$, 
with $j$, $M\subseteq\uniV[\genH]$ \st\ \xitemof{III:gen-large-1-0} $\sim$ 
\xitemof{III:gen-large-6} and \xitemof{III:gen-large-6-0} hold. 

By \xitemof{III:gen-large-4}, we have $\genG\in M$. Let $\poR=j(\poP)$. By 
elementarity, 
\begin{xitemize}
\xitem[III:gen-large-7-0] 
  $M\models\qbox{\poR\xmbox{ is }\LT j(\kappa)\mbox{-directed closed}}$. 
\end{xitemize}

By \xitemof{III:gen-large-7-0} and \xitemof{III:gen-large-6-0},
\ifextended{\color{darkelectricblue}
$\uniV[\genH]\modelof{\poR\mbox{ is }\sigma\mbox{-directed closed}}$. 
Thus, $\uniV[\genH]\modelof{\poR\models\gmQ}$ by \xitemof{III:gen-large-6-13}. }\else
$\uniV\modelof{\poR\mbox{ is }\sigma\mbox{-directed closed}}$. 
Thus, $\uniV\modelof{\poR\models\gmQ}$ by \xitemof{III:gen-large-7}. \fi

By \xitemof{III:gen-large-3}, $M\models\cardof{j\imageof\genG}<j(\kappa)$. 
\ifextended{\color{darkelectricblue}
Hence, by \xitemof{III:gen-large-7-0}, there is (a master condition) $\condr\in\poR$ \st\
$M\models\condr\leq_\poR j\imageof\genG$. 

}\else
Hence, there is (a master condition) $\condr\in\poR$ \st\
$M\models\condr\leq_\poR j\imageof\genG$ by \xitemof{III:gen-large-7-0}. \fi
Let $\genK$ be 
a $(\uniV[\genH],\poR)$-generic filter with $\condr\in\genK$.
Then
\begin{xitemize}
\xitem[III:gen-large-8] 
  $\elembed{\tilde{j}}{\uniV[\genG]}{M[\genK]}$; $\uta^\genG\mapsto j(\uta)^{\genK}$
\end{xitemize}
is well defined and $j\subseteq\tilde{j}$. In particular we have
$\kappa=\crit(\tilde{j})$, $\tilde{j}>\lambda$ and
$\tilde{j}\imageof\lambda\in M[\genK]$.
\ifextended{\color{darkelectricblue}
we have
$\uniV[\genH]\modelof{\poR\mbox{ is }\sigma\xmbox{-directed closed}}$ by 
\xitemof{III:gen-large-7-0}, \xitemof{III:gen-large-6-0}, and since
$M[\genK]\modelof{\poR\mbox{ is }\sigma\xmbox{-directed closed}}$ by elementarity. 
}\else
By \xitemof{III:gen-large-7-0} and \xitemof{III:gen-large-6}, we have
$\uniV[\genH]\modelof{\poR\mbox{ is }\sigma\xmbox{-directed closed}}$.\fi
Thus, in $\uniV[\genG]$, letting $\poQ=\utpoQ[\genG]$, 
with the $\poQ$-name $\utpoR$ corresponding to $\poR$ \st\ 
$\forces{\poQ}{\utpoR\mbox{ is }\LT\theta\xmbox{-directed closed}}$,
$\poQ\ast\utpoR$ satisfies $\gmQ$ and it induces a generic elementary 
embedding for generic $\lambda$-supercompactness. 
Since \ifextended{\color{darkelectricblue} $\theta$ }\else$\genG$ \fi was arbitrary, it 
follows that \ifextended{\color{darkelectricblue} \xitemaof{III:gen-large-7-a} }\else\xitemof{III:gen-large-7} \fi holds.  
\qedofProp

\begin{Cor}
  \Label{III:P-gen-large-2} Suppose that $\kappa$ is strongly 
  Laver-generically supercompact for $(\gmP,\gmQ)$ where $\gmP$ 
  is the property ``forcing equivalent to a \po\ of the form $\Col(\kappa,\mu)$ for some $\mu$'' and $\gmQ$ is ``proper''. Suppose further that 
$\kappa_1>\kappa$ is a supercompact and let $\poP=\Col(\kappa,\kappa_1)$. Then \medskip
\\
  \wassert{a}
    $\forces{\poP}{\kappa\mbox{ is generically supercompact by proper
      \pos}}$; \smallskip\\
  \wassert{b} 
    $\forces{\poP}{\kappa^+\mbox{ is generically supercompact by }\LT\kappa\mbox{-closed \pos}};$\smallskip\\
  \wassert{c} $\forces{\poP}{\SDLS^{int}_+(\calL^{\aleph_0}_{stat},\LT\kappa)}$; and\smallskip\\
  \wassert{d} $\forces{\poP}{\GRP^{\LT\kappa}(\LE\kappa)}$. 
\end{Cor}
\prf \assertof{a}: By \Propof{III:P-gen-large-1}.

\noindent
\assertof{b}: By Lemma 4.10 in [Fuchino, Sakai and Ottenbreit\cite{III:bib-I}].

\noindent
\assertof{c}: By \assertof{a} and, Theorem 2.10 
and Propositions 3.1 
in [Fuchino, Sakai and Ottenbreit\cite{III:bib-II}].

\noindent
\assertof{d}: By \assertof{b} 
above and Lemma 4.11 in [Fuchino, Sakai and Ottenbreit\cite{III:bib-I}]. 
\qedofCor

\memo{!!!!}
\begin{Prop}
  \Label{III:P-gen-large-2-0}
  Suppose that\/ $\gmP$ is the property ``forcing equivalent to a \po\ of the form
  $\Col(\kappa,\mu)$ for some $\mu$'' and $\gmQ$ the property of \pos\ \st\
  \begin{xitemize}
  \xitem[III:gen-large-8-a] 
    $\forall P\,\forall \utilde{Q}\,({}\,
    \begin{array}[t]{@{}l}
      P\models\gmP\,\land\,
      \forces{P}{\utilde{Q}\mbox{ is a }\LT\theta\mbox{-directed closed \po\,}}\\
      \rightarrow P\ast\utilde{Q}\models\gmQ)
    \end{array}$ 
  \end{xitemize}
  for a cardinal $\theta$.
  
  If $\kappa$ is Laver-generically strongly superhuge for $(\gmP,\gmQ)$, 
  then, for any cardinal $\mu>\kappa$ and $\poP=\Col(\kappa,\mu)$, we have
  \begin{xitemize}
  \xitem[III:gen-large-8-a-a] 
    $\forces{\poP}{\kappa\mbox{ is strongly generically superhuge for }\gmQ}$.
  \end{xitemize}
\end{Prop}
\prf
Let $\genG$ be a $(\uniV,\poP)$-generic filter. For cardinals $\lambda$, let
$\lambda'=\max\ssetof{(\mu^{\LT\kappa})^+, \lambda,\theta}$. In $\uniV[\genG]$, let $\poQ$ 
be a \po\  with
$\poQ\models\gmQ$ and $\genH$ a $(\uniV[\genG],\poQ)$-generic filter \st\ there are $j$,
$M\subseteq\uniV[\genG][\genH]$ satisfying: $M$ is a transitive class in
$\uniV[\genG][\genH]$; $\elembed{j}{\uniV}{M}$;   $\crit(j)=\kappa$; $j(\kappa)>\lambda'$;
$\genG$, $\genH\in M$; and 
\begin{xitemize}
\xitem[III:P-gen-large-2-1] 
  $([M]^{j(\kappa)})^{\uniV[\genG][\genH]}\subseteq M$.
\end{xitemize}

Let $\poP^*=j(\poP)$. By elementarity and \xitemof{III:P-gen-large-2-1}, we have
$\poP^*=\Col(j(\kappa,j(\mu)))$ in $\uniV[\genG][\genH]$. Since $\poP^*$ 
is $j(\kappa)$-directed closed (in $M$ or in $\uniV[\genG][\genH]$), there is 
(a master condition) $\condr\in\poP^*$ with $\condr\leq_{\poP^*}j(\condp)$ for all $\condp\in\genG$. 
Let $\genG^*$ be $(\uniV[\genG][\genH],\poP^*)$-generic filter with $\condr\in\genG^*$. 
Then $j$ is lifted to
\begin{xitemize}
\xitem[III:P-gen-large-2-2] 
  $\elembed{\tilde{j}}{\uniV[\genG]}{M[\genG^*]}$; $\uta[\genG]\mapsto j(\uta)[\genG^*]$
\end{xitemize}
for $\poP$-name $\uta$. $j\subseteq\tilde{j}$ and hence it is clear that $\tilde{j}$ satisfies:
$\crit(\tilde{j})=\kappa$; $\tilde{j}(\kappa)>\lambda$;
$\genG^*\in M[\genG^*]$ since $\utgenG^*\in M$ by \xitemof{III:P-gen-large-2-1} and
$M\models\ZFC$ where $\utgenG^*$ is the standard $\poP^*$-name of $\genG^*$. We also have
\begin{xitemize}
\xitem[III:P-gen-large-2-3] 
  $([M[\genG^*]]^{j(\kappa)})^{\uniV[\genG][\genH][\genG^*]}\subseteq M[\genG^*]$ 
\end{xitemize}
by \Lemmaof{III:P-gen-large-a-1}. This shows that \xitemof{III:gen-large-8-a-a} holds. 
\qedofProp

\begin{Prop}
  \Label{III:P-gen-large-3}\wassertof{1} Let $\gmP$ be the property ``forcing equivalent 
  to a \po\ of the form $\Col(\kappa,\mu)$ 
  for some $\mu$'' and $\gmQ_\theta$ is the property ``forcing equivalent to a regular sub-\po\ of the 
  completion of a \po\ of the form `generalized Cohen \po\ $\times$ $\LT\theta$-closed \po' ''. 
  If $\kappa$ is strongly Laver-generically supercompact for $(\gmP,\gmQ_\theta)$ for all
  $\theta\in\Card$, then, for any cardinal $\mu$ and $\poP=\Col(\kappa,\mu)$, we have
  \begin{xitemize}
  \xitem[] 
    $\forces{\poP}{\HH(\LT\kappa)}$.
  \end{xitemize}

  \wassert{2} Let $\gmP$ be the property ``forcing equivalent to a \po\ of the form $\Col(\kappa,\mu)$ 
  for some $\mu$'' and $\gmQ_\theta$ is the property ``forcing equivalent to a regular sub-\po\ of the 
  completion of the \po\ of the form `ccc \po\ $\times$ $\LT\theta$-closed \po' ''. 
  If $\kappa$ is strongly Laver-generically supercompact for $(\gmP,\gmQ_\theta)$ for all
  $\theta\in\Card$, then, for any cardinal $\mu$ and $\poP=\Col(\kappa,\mu)$, we have
  \begin{xitemize}
  \xitem[] 
    $\forces{\poP}{\SDLS^{int}_+(\calL^\PKL_\stat,\LT\kappa)}$.
  \end{xitemize}

\end{Prop}
\memo{Scan\_2020-02-03..., p.19,20}
\prf \assertof{1}: Suppose that $\poP=\Col(\kappa,\mu)$ and $\genG$ is 
a $(\uniV,\poP)$-generic filter. In $\uniV[\genG]$, let $X=\pairof{X,\tau}$ be a 
non-metrizable topological space \st\ 
\begin{xitemize}
\xitem[III:gen-large-8-a-0] 
  $\chi(a,X)<\kappa$ for all $a\in X$. 
\end{xitemize}
Let
$\lambda_0=\cardof{X}$, 
$\theta=\lambda=\max\ssetof{(2^{\lambda_0})^+, \mu^{\LT\kappa}}$ 
and let $\elembed{j}{V}{M\subseteq\uniV[\genG][\genH]}$ be \st\ $\crit(j)=\kappa$,
\begin{xitemize}
\xitem[III:gen-large-8-a-1] $j(\kappa)>\lambda$, $j\imageof\lambda\in M$, and 
\xitem[III:gen-large-8-0] 
  $\left([M]^{\aleph_0}\right)^{\uniV[\genG][\genH]}\subseteq M$, 
\end{xitemize}
where 
$\genH=\tilde{\genH}\cap\poQ$ for a $(\uniV[\genG],\tilde{\poQ})$-generic filter
$\tilde{\genH}$ for a \po\ $\tilde{\poQ}$ in $\uniV[\genG]$ of the form
\begin{xitemize}
\xitem[III:gen-large-9] $\tilde{\poQ}\ \sim\  $generalized Cohen \po\ $\times$ $\LT\theta$-closed \po
\end{xitemize}
and $\poQ\circleq \tilde{\poQ}$. 

Let $\poP^*=j(\poP)$. By elementarity
$M\modelof{\poP^*\xmbox{ is }j(\kappa)\xmbox{-directed closed}}$. 
Since $\cardof{\poP}\leq\lambda$, $j\imageof\genG\in M[\genG]$ by \xitemof{III:gen-large-8-a-1} 
(see Lemma 2.5 in [Fuchino, Ottenbreit and Sakai\cite{III:bib-II}]). 
Since $\lambda<j(\kappa)$ there is $\condr\in\poP^*$ \st\ $\condr\leq_{\poP^*}j(\condp)$ 
for all $\condp\in\genG$. Let $\genG^*$ be a $(\uniV[\genG][\tilde{\genH}],\poP^*)$-generic filter 
with $\condr\in\genG^*$. $j$ is then lifted to
\begin{xitemize}
\xitem[III:gen-large-10] 
  $\elembed{\tilde{j}}{\uniV[\genG]}{M[\genG^*]
  \subseteq\uniV[\genG][\genH][\genG^*]\subseteq\uniV[\genG][\tilde{\genH}][\genG^*]}\,;$\\[\jot]
  \phantom{$\elembed{\tilde{j}}{\uniV[\genG]}{M[\genG^*]\subseteq{}}$}$\uta[\genG]\mapsto j(\uta)[\genG^*]$. 
\end{xitemize}

Let $\tau_0=\setof{j(O)\cap j\imageof X}{O\in\tau}$. By \xitemof{III:gen-large-8-a-0} and 
since $\kappa$ is the critical point of $j$, we have 
\begin{xitemize}
\xitem[III:gen-large-11] 
  $M\modelof{\pairof{j\imageof X,\tau_0}\mbox{ is a subspace of }\pairof{j(X),j(\tau)}}$.
\end{xitemize}
We also have 
\begin{xitemize}
\xitem[III:gen-large-12] 
  $M[\genG^*]\modelof{\pairof{j\imageof X,\tau_0}
  \mbox{ is homeomorphic to }\pairof{X,\tau}}$. 
\end{xitemize}
Hence the same property holds in $\uniV[\genG][\tilde{\genH}][\genG^*]$. 

Now generalized Cohen \po\ part of $\tilde{\genH}$ preserve the non-metrizability of 
$\pairof{X,\tau}$ by \Thmof{III:P-refl-a}. By the $\LT\theta$-closed part of $\tilde{\genH}$ 
no new metric on $X$ is added. Hence
$\uniV[\genG][\tilde{\genH}]\modelof{\pairof{X,\tau}\xmbox{ is non-metrizable}}$. It 
follows that $M\modelof{\pairof{X,\tau}\xmbox{ is non-metrizable}}$ and hence 
by $\LT\lambda$-closedness of $\poP^*$, it follows that
$M[\genG^*]\modelof{\pairof{X,\tau}\xmbox{ is non-metrizable}}$. Thus by 
\xitemof{III:gen-large-12},
$M[\genG^*]\modelof{\pairof{\tilde{j}\imageof X,\tau_0}\xmbox{ is non-metrizable}}$. Thus
\begin{xitemize}
\xitem[] 
  $M[\genG^*]\modelof{\tilde{j}(X)\xmbox{ has a non-metrizable subspace }Y\xmbox{ \\
  \phantom{\(M[\genG^*]\models''\)}of cardinality }\LT\tilde{j}(\kappa)}$.
\end{xitemize}
By elementarity of $\tilde{j}$ it follows that
\begin{xitemize}
\xitem[] 
  $\uniV[\genG]\modelof{X\mbox{ has a non-metrizable subspace }Y\mbox{ of cardinality }\LT \kappa}$.
\end{xitemize}
\smallskip

\assertof{2}: The proof is done similarly to \assertof{1}, by using 
\Lemmaof{III:L-preserv-a} in place of \Thmof{III:P-refl-a}.
\qedofProp

\begin{Lemma}
  \Label{III:P-gen-large-3-0}
  Let $\gmQ_\theta$ for a cardinal $\theta$ be as in {\em\rm\Propof{III:P-gen-large-3},\,\assertof{2}} and assume that 
  $\kappa$ is strongly generically superhuge for $\gmQ_\theta$ for all
  $\theta\in\Card$. Then, for any
  $\lambda\geq\kappa$, $\Pkl{\kappa}{\lambda}$ carries a $\sigma$-saturated normal 
  ideal. 
\end{Lemma}
\prf
Let $\lambda\geq\kappa$ and let $\poQ\circleq RO(\poS\times\poT)$ be \st\  
$\poS$ is $\LT\big(2^{2^{(\lambda^{<\kappa})}}\big)^+$-closed \po, $\poT$ is ccc, and that there are 
a $(\uniV,\poQ)$-generic filter $\genH$ and $j$, $M\subseteq\uniV[\genH]$ \st\
\ifextended{\color{darkelectricblue}$M$ is an inner model in $\uniV[\genH]$, }\fi $\elembed{j}{\uniV}{M}$,
$\crit(j)=\kappa$, $j(\kappa)>\lambda$,  
($j\imageof\lambda\in M$) and $([M]^{j(\kappa)})^{\uniV[\genH]}\subseteq M$. 
Let $\genK$ be a $(\uniV,\poS)$-generic filter and $\genL$ 
a $(\uniV[\genK],\poT)$-generic filter \st\ $\uniV[\genH]\subseteq\uniV[\genK][\genL]$. In
$\uniV[\genK][\genL]$,
\begin{xitemize}
\xitem[] 
  $\calI=\setof{X\in\left(\psof{\Pkl{}{}^\uniV}\right)^\uniV}{j\imageof\lambda\not\in j(X)}$
\end{xitemize}
is 
a $V$-normal ideal. Since $\uniV[\genK]\modelof{\poT\mbox{ is ccc}}$, it follows 
that, in $\uniV[\genK]$, 
$\calI'=\setof{X\in\left(\psof{\Pkl{}{}^\uniV}\right)^\uniV}{\forces{\poT}{\check{X}\in\utilde{\calI}}}$ is a
$\sigma$-saturated $\uniV$-normal ideal for $\poT$-name $\utilde{\calI}$ of $\calI$. Now,  
by the closedness of $\poS$, $\calI'\in\uniV$  
and $\calI'$ is a $\sigma$-saturated normal ideal in $\uniV$. 
\qedofLemma 

\memo{!!!!}
\begin{Prop}
  \Label{III:P-gen-large-4} Let $\gmP$ and $\gmQ_\theta$ for a cardinal $\theta$ be as in 
  {\em\rm\Propof{III:P-gen-large-3},\,\assertof{2}}. If $\kappa$ is strongly Laver-generically 
  superhuge for $(\gmP,\gmQ_\theta)$ for all cardinal $\theta$, then, for any $\lambda\geq\kappa$, and
  $\poP=\Col(\kappa,\lambda)$, 
  \begin{xitemize}
  \xitem[] 
    $\forces{\poP}{\Pkl{}{}\mbox{ carries a }\sigma\mbox{-saturated normal ideal}}$. 
  \end{xitemize}
\end{Prop}
\prf By \Propof{III:P-gen-large-2-0} and \Lemmaof{III:P-gen-large-3-0}. \qedofCor

\section{Mixed support iteration}\Label{III:msi}

The construction of the mixed support iteration we give here is similar to the one given 
in [Krueger\cite{III:bib-krueger1,III:bib-krueger2}]. Nevertheless, we will examine the  
details of our construction in the following, since there are 
a couple of points organized differently from [Krueger\cite{III:bib-krueger1,III:bib-krueger2}]. 

In this section, $\kappa$ is always a fixed supercompact cardinal and 
$\mapping{f}{\kappa}{\uniV_\kappa}$ is a Laver function, i.e.\ a function satisfying:
\begin{xitemize}
\xitem[III:msi-a-0] 
  for any set 
  $a$ and any $\lambda\geq\kappa$, there is $\elembed{j}{V}{M}$ \st\
  $crit(j)=\kappa$, $j(\kappa)>\lambda$,  
  $[M]^{\lambda}\subseteq M$ and $j(f)(\kappa)=a$
\end{xitemize}
(see e.g.\ Theorem 20.21 in [Jech\cite{III:bib-millennium-book}]).

Let $\mapping{\overline{f}}{\kappa}{\kappa}$ be defined by
\begin{xitemize}
\xitem[III:msi-a-1] $\overline{f}(\alpha)=\cardof{\trcl(f(\alpha))}$\ 
  \  for $\alpha<\kappa$. 
\end{xitemize}
Let
\begin{xitemize}
\xitem[III:msi-0]%
  $S=\setof{\alpha<\kappa}{{}
  \begin{array}[t]{@{}l}
    \alpha\mbox{ is a strongly Mahlo cardinal}\\[\jot]
    \mbox{closed \wrt\ }\overline{f}\ },\mbox{\ \  and let}
  \end{array}
  $
\xitem[III:msi-1-a] 
  $T=\kappa\setminus S$.
\end{xitemize}

Let $\mapping{\nu}{\kappa}{\kappa}$ be the mapping defined by 
\begin{xitemize}
\xitem[III:msi-1-0] 
  $\nu(\alpha)=\min(S\setminus(\alpha+1))$ for $\alpha\in\kappa$. 
\end{xitemize}

We treat iterations here as in [Jech\cite{III:bib-millennium-book}] \st\ elements 
of $\alpha$\/th step $\poP_\alpha$ of an iteration
$\seqof{\poP_\alpha,\utpoQ_\beta}{\alpha\leq\kappa,\beta<\kappa}$ are sequences 
of length $\alpha$.

Let $\seqof{\poO_\alpha,\utpoR_\beta}{\alpha\leq\kappa,\beta<\kappa}$ be a finite 
support iteration of ccc \pos\ which will be further specified later. This 
preparatory iteration should satisfy the following conditions:
\begin{xitemize}
\xitem[III:msi-2]
  $\utpoR_\alpha\in V_{\nu(\alpha)}$, and
\xitem[III:msi-3] 
  $\forces{\poO_\alpha}{\utpoR_\alpha\equiv\ssetof{\bbone_{\utpoR_\alpha}}}$, for 
  all $\alpha\in S$.
\end{xitemize}
We denote the canonical embeddings of $\poO_\alpha$ into $\poO_\beta$ for
$\alpha\leq\beta\leq\kappa$ by $i^*_{\alpha,\beta}$. 
Thus, $i^*_{\alpha,\beta}$ is the mapping defined by 
$i^*_{\alpha,\beta}(\condp)=\condp\cup\vec{\bbone}_{\alpha,\beta}$, where
$\vec{\bbone}_{\alpha,\beta}$ is the function $g$ on $\beta\setminus\alpha$ with
$g(\xi)=\bbone_{\utpoR_\xi}$ for all $\xi\in\beta\setminus\alpha$.

$\seqof{\poP_\alpha,\utpoQ_\beta}{\alpha\leq\kappa,\beta<\kappa}$\vspace{-0.8\smallskipamount} is our final 
iteration which is specified once the preparatory iteration
$\seqof{\poO_\alpha,\utpoR_\beta}{\alpha\leq\kappa,\beta<\kappa}$
\vspace{-0.8\smallskipamount} is fixed. The iteration
$\seqof{\poP_\alpha,\utpoQ_\beta}{\alpha\leq\kappa,\beta<\kappa}$ is 
defined recursively in \assertof{A} and \assertof{B} below, together with the 
commutative systems of 
complete embeddings $\mapping{\iota_\alpha}{\poO_\alpha}{\poP_\alpha}$, and
$\mapping{i_{\alpha,\beta}}{\poP_\alpha}{\poP_\beta}$ for 
$\alpha\leq\beta\leq\kappa$ which should satisfy
\begin{xitemize}
\xitem[III:msi-5] $i_{\alpha,\alpha}=id_{\poP_\alpha}$ for $\alpha\leq\kappa$;
\xitem[III:msi-5-0] $\iota_\beta\circ i^*_{\alpha,\beta}=i_{\alpha,\beta}\circ\iota_\alpha$, and
\xitem[III:msi-5-1] $i_{\beta,\gamma}\circ i_{\alpha,\beta}=i_{\alpha,\gamma}$ for
  $\alpha\leq\beta\leq\gamma\leq\kappa$;
\xitem[III:msi-8-0] 
  $\supp(\iota_\beta(\condo))=\supp(\condo)$, where $\supp(\cdot)$ is defined as 
  in \xitemof{III:msi-6} below, and
\xitem[III:msi-8-1] 
  $\iota_\alpha(\condo\restr\alpha)=\iota_\beta(\condo)\restr\alpha$ for 
$\alpha<\beta\leq\kappa$ and $\condo\in\poO_\beta$.
\end{xitemize}

We define now the {\it Easton-type mixed support} of the iteration as a sequence 
$\seqof{\calI_\alpha}{\alpha<\kappa}$ of ideals where each $\calI_\alpha$ for
$\alpha\leq\kappa$ is an ideal over $\alpha$.

\begin{xitemize}
\xitem[III:msi-8-2]
  $\calI_{\alpha+1}=\calI_\alpha\cup\setof{s\cup\ssetof{\alpha}}{s\in\calI_\alpha}$ for all $\alpha<\kappa$;
\xitem[III:msi-8-3] If $\gamma<\kappa$ is a limit ordinal but not a regular cardinal, then
  $\calI_{\gamma}=\setof{s\subseteq\gamma}{s\cap\alpha\in\calI_\alpha
  \xmbox{ for all }\alpha<\gamma\mbox{ and }s\cap T\xmbox{ is bounded in }\gamma}$;
\xitem[III:msi-8-4] If $\gamma\leq\kappa$ is a regular cardinal, then 
  $\calI_{\gamma}=\setof{s\subseteq\gamma}{s\cap\alpha\in\calI_\alpha
  \xmbox{ for all }\alpha<\gamma\mbox{ and }s\xmbox{ is bounded in }\gamma}$. 
\end{xitemize}

The following is easy to prove by induction on $\alpha\leq\kappa$:
\begin{Lemma}
  \Label{III:P-msi-0}\wassert{1} $\calI_\alpha$ is an ideal over $\alpha$ with
  $\setof{\ssetof{\beta}}{\beta<\alpha}\subseteq\calI_\alpha$ for all $\alpha\leq\kappa$.\smallskip

  \wassert{2} For 
  all $\alpha\leq\kappa$, $s\in\calI_\alpha$\ \ $\Leftrightarrow$\ \
  $s\cap T$ is finite and $\cardof{s\cap\mu}<\mu$ for all 
  regular infinite cardinal $\mu\leq\alpha$.\qed
\end{Lemma}

Now we are ready to define the iteration
$\seqof{\poP_\alpha,\utpoQ_\beta}{\alpha\leq\kappa,\beta<\kappa}$ in the 
following \assertof{A} and \assertof{B}:
\medskip

\assertof{A}
If $\seqof{\poP_\alpha,\utpoQ_\beta}{\alpha<\gamma,\beta<\gamma}$ has been 
defined for a limit $\gamma\leq\kappa$, let
\begin{xitemize}
\xitem[III:msi-8-5] 
  $\poP_\gamma=\setof{\condp}{{}
  \begin{array}[t]{@{}l}
    \condp\mbox{ is a sequence of length }\gamma,\\
    \condp\restr\alpha\in\poP_\alpha\mbox{ for all }\alpha<\gamma, \mbox{ and }
    \supp(\condp)\in\calI_\gamma\ }
  \end{array}$ 
\end{xitemize}
where
\begin{xitemize}
\xitem[III:msi-6] 
  $\supp(\condp)=\setof{\alpha<\gamma}{\condp(\alpha)\not=\bbone_{\utpoQ_\alpha}}$. 
\end{xitemize}

For $\condp$, $\condq\in\poP_\gamma$, 
\begin{xitemize}
\xitem[III:msi-8] 
  $\condp\leq_{\poP_\gamma}\condq$\ \ $\Leftrightarrow$\ \
  $\condp\restr\delta\leq_{\poP_\delta}\condq\restr\delta$ for 
  all $\delta<\gamma$. 
\end{xitemize}

We assume that the complete embeddings
$\mapping{\iota_\alpha}{\poO_\alpha}{\poP_\alpha}$, $\alpha<\gamma$  have been 
defined \st\   
\xitemof{III:msi-8-0} and \xitemof{III:msi-8-1} hold for all $\alpha\leq\beta<\gamma$. 
The mapping 
$\iota_\gamma$ is then defined by: 
\begin{xitemize}
\xitem[III:msi-9] 
  $\mapping{\iota_\gamma}{\poO_\gamma}{\poP_\gamma}$;\ \ 
  $\condo\mapsto\bigcup_{\delta<\gamma}\iota_\delta(\condo\restr\delta)$.
\end{xitemize}

For $\delta\leq\gamma$, let
$\vec{\bbone}_{\delta,\gamma}=\setof{\pairof{\alpha,\bbone_{\utpoQ_\alpha}}}{\delta\leq\alpha<\gamma}$ 
as before and let $i_{\delta,\gamma}$ be defined by
\begin{xitemize}
\xitem[III:msi-9-0] 
  $\mapping{i_{\delta,\gamma}}{\poP_\delta}{\poP_\gamma};\ \ 
  \condp\mapsto\condp\cup\vec{\bbone}_{\delta,\gamma}$.
\end{xitemize}

It is easy to check that $\iota_\gamma$ and $i_{\delta,\gamma}$, 
$\gamma\leq\delta$ are complete embeddings and 
\xitemof{III:msi-5} $\sim$ \xitemof{III:msi-8-0} hold for all indices
$\LE\gamma$. \smallskip

\assertof{B}
Now suppose 
that $\seqof{\poP_\alpha,\utpoQ_\beta}{\alpha\leq\gamma,\beta<\gamma}$,
$\seqof{\iota_\beta}{\beta<\gamma}$ and $\seqof{i_{\alpha,\beta}}{\alpha\leq\beta\leq\gamma}$ have been 
defined for some $\gamma<\kappa$.\smallskip

\assert{a} If $\gamma\in S$ and
\begin{xitemize}
\xitem[III:msi-9-1] 
  $f(\gamma)=\pairof{\mu,\theta,R}$ for some 
  cardinals $\mu$, $\theta>\gamma$ and a set $R$,\footnote{At the moment, $R$ does not play any 
    role. This component is added here so that we can later modify the construction.}
\end{xitemize}
then let
\begin{xitemize}
\xitem[III:msi-10]
  $\utpoQ_\gamma=(\Col(\gamma,\mu))^\bullet_{\poP_\gamma}$, and
\xitem[III:msi-11] 
  $\poP_{\gamma+1}=\setof{\condp\cup\ssetof{\pairof{\gamma,\utcondq}}}{
  \begin{array}[t]{@{}l}
    \condp\in\poP_\gamma,\utcondq\mbox{ is a 
      canonical }\poP_\gamma\mbox{-name\,\footnotemark}\\[-1\jot]
    \mbox{\st\ }
    \forces{\poP_\gamma}{\utcondq\varin\utpoQ_\gamma}}.
  \end{array}$
\end{xitemize}
\footnotetext{Adopting the terminology of [Cummings\cite{III:bib-cummings}], we call 
  a $\poP$-name $\uta$ is a {\it canonical\/ $\poP$-name} if, for any $\poP$-name
  $\utb$ with $\forces{\poP}{\uta\equiv\utb}$, we have $\cardof{\trcl(\uta)}\leq\cardof{\trcl(\utb)}$. }

For $\condp_0\cup\ssetof{\pairof{\gamma,\utcondq_0}}$,
$\condp_1\cup\ssetof{\pairof{\gamma,\utcondq_1}}\in\poP_{\gamma+1}$, 
\begin{xitemize}
\xitem[III:msi-12] 
  $\condp_0\cup\ssetof{\pairof{\gamma,\utcondq_0}}\ 
  \leq_{\poP_{\gamma+1}}\ \condp_1\cup\ssetof{\pairof{\gamma,\utcondq_1}}$
  \ $\Leftrightarrow$\ \ 
  $\begin{array}[t]{@{}ll}
    \condp_0\leq_{\poP_\gamma}\condp_1\mbox{ and }\\[\jot]
    \condp_0\forces{\poP_\gamma}{\utcondq_0\leq_{\utpoQ_\gamma}\utcondq_1}.
  \end{array}$
\end{xitemize}

For $\condo\in\poO_{\gamma+1}$, let
\begin{xitemize}
\xitem[III:msi-12-0] 
  $\iota_{\gamma+1}(\condo)=\iota_\gamma(\condo\restr\gamma)\cup\ssetof{\pairof{\gamma,\bbone_{\utpoQ_\gamma}}}$, 
\end{xitemize}
and, for $\alpha\leq\gamma$ and $\condp\in\poP_\alpha$, let
\begin{xitemize}
\xitem[III:msi-12-1] 
  $i_{\alpha,\gamma+1}(\condp)=i_{\alpha,\gamma}(\condp)\cup\ssetof{\pairof{\gamma,\bbone_{\utpoQ_\gamma}}}$.
\end{xitemize}
\smallskip

\assert{b} If $\gamma\in S$ but \assertof{a} does not hold, then let 
$\utpoQ_\gamma$ be a $\poP_\gamma$-name of trivial forcing and the rest is 
treated just as in the case \assertof{a}.

In both of the cases \assertof{a} and \assertof{b}, it is clear that the defined mappings 
are complete embeddings and satisfy \xitemof{III:msi-5} $\sim$ \xitemof{III:msi-8-0}.
\smallskip

\assert{c} If $\gamma\not\in S$, then let $\utpoQ_\gamma$ be 
the $\poP_\gamma$-name $\iota_\gamma(\utpoR_\gamma)$ and\,\footnote{With
  $\iota_\gamma$, we also denote the embedding 
  of $V^{\poO_\gamma}$ into $V^{\poP_\gamma}$ canonically induced by $\iota_\gamma$. }
\begin{xitemize}
\xitem[III:msi-13] 
  $\poP_{\gamma+1}=\setof{\condp\cup\ssetof{\pairof{\gamma,\iota_\gamma(\utcondr)}}}{
  \begin{array}[t]{@{}l}
    \condp\in\poP_\gamma,\ 
    \utcondr\mbox{ is a canonical }\poO_\gamma\mbox{-name}\\
    \mbox{\st\ }
    \forces{\poO_\gamma}{\utcondr\varin\utpoR_\gamma}}.
  \end{array}$
\end{xitemize}

For $\condp_0\cup\ssetof{\pairof{\gamma,\iota_\gamma(\utcondr_0)}}$,
$\condp_1\cup\ssetof{\pairof{\gamma,\iota_\gamma(\utcondr_1)}}\in\poP_{\gamma+1}$, 
\begin{xitemize}
\xitem[III:msi-14] 
  $\condp_0\cup\ssetof{\pairof{\gamma,\iota_\gamma(\utcondr_0)}}\ 
  \leq_{\poP_{\gamma+1}}\ \condp_1\cup\ssetof{\pairof{\gamma,\iota_\gamma(\utcondr_1)}}$
  \\
  $\Leftrightarrow$\ \ 
  $\begin{array}[t]{@{}l}
    \condp_0\leq_{\poP_\gamma}\condp_1\mbox{ and there is }\condo\in\poO_\gamma\mbox{ \st\ }
    \condp_0\leq_{\poP_\gamma}\iota_\gamma(\condo)\\
    \mbox{and }\condo\forces{\poO_\gamma}{\utcondr_0\leq_{\utpoR_\gamma}\utcondr_1}.
  \end{array}$ 
\end{xitemize}

For $\condo\in\poO_{\gamma+1}$, let
\begin{xitemize}
\xitem[III:msi-15] 
  $\iota_{\gamma+1}(\condo)=\iota_\gamma(\condo\restr\gamma)\cup\ssetof{\pairof{\gamma,\iota_\gamma(\condo(\gamma))}}$, 
\end{xitemize}
and, for $\alpha\leq\gamma$ and $\condp\in\poP_\alpha$, let
\begin{xitemize}
\xitem[III:msi-16] 
  $i_{\alpha,\gamma+1}(\condp)=i_{\alpha,\gamma}(\condp)\cup\ssetof{\pairof{\gamma,\bbone_{\utpoQ_\gamma}}}$.
\end{xitemize}

Also in this case, the mapping introduced are complete embeddings and 
\xitemof{III:msi-5} $\sim$ \xitemof{III:msi-8-0} are satisfied.

This finishes the 
construction of our Easton-type mixed support iteration.\medskip

The following three Lemmas can be proved easily with the standard argument in the order as we present them here. 

\begin{Lemma}
  \Label{III:P-msi-0-a-0}
  For an ordinal $\gamma\leq\kappa$ and a $\gamma$-sequence $\condp$, 
  \begin{xitemize}
  \item[] $\condp\in\poP_\gamma$\ \ $\Leftrightarrow$\ \ 
    $\begin{array}[t]{@{}l}
    \forces{\poP_\xi}{\condp(\xi)\varin\utpoQ_\xi}\mbox{ for all }\xi\in S\cap\gamma,\\
    \condp(\xi)=\iota_\xi(\utcondr)\mbox{ for a canonical }\poO_\gamma\mbox{-name }\utcondr\mbox{ with}\\
    \forces{\poO_\xi}{\utcondr\varin\utpoR_\xi}\mbox{ for all }\xi\in T\cap\gamma\mbox{, and}\\
    \supp(\condp)=\setof{\xi<\gamma}{\condp(\xi)\not\equiv\bbone_{\utpoQ_\xi}}\in\calI_\gamma.
    \end{array}$
  \end{xitemize}\vspace{-1.4\baselineskip}

\qed
\end{Lemma}

\begin{Lemma}
  \Label{III:P-msi-0-a-1}
  For $\delta\leq\gamma\leq\kappa$, $\condp_0$, $\condp_1\in\poP_\gamma$ 
  and $s\subseteq\gamma$, if $\supp(\condp_0)\subseteq\delta\cup s$,
  $\supp(\condp_1)\subseteq\delta\cup(\gamma\setminus s)$ and
  $\condp_0\restr\delta\leq_{\poP_\delta}\condp_1\restr\delta$, then
  \begin{xitemize}
  \xitem[] 
    $\condp_2=(\condp_0\restr(\delta\cup s)\,)\cup(\condp_1\restr(\gamma\setminus(\delta\cup s))\,)$\\[\jot]
    $\big(=(\condp_0\restr\supp(\condp_0))\cup(\condp_1\restr(\gamma\setminus\supp(\condp_0)))\ \big)\in\poP_\gamma$ 
  \end{xitemize}
  and
  $\condp_2$ is a maximal element 
  of\/ $\poP_\gamma$ below $\condp_0$ and $\condp_1$ \wrt\ $\leq_{\poP_\gamma}$. 
\end{Lemma}
\prf
\ifextended{\small\color{darkelectricblue}
We prove the assertion of the Lemma by induction on $\gamma$ with
$\delta\leq\gamma\leq\kappa$.

The rest will be written later. (Scan\_2020-02-03... p.18)
}
\else 
By induction on $\gamma$ with $\delta\leq\gamma\leq\kappa$. 
\fi 
\qedofLemma\qedskip

Note that, in the Lemma above, we are talking about ``a'' maximal element since 
$\leq_{\poP_\gamma}$ is merely a pre-ordering in general. 

The following can be proved applying the Pressing-down Lemma and \Lemmaof{III:P-msi-0-a-1} above. 
Note that, for $\alpha\in S\cup\ssetof{\kappa}$,
$R=\setof{\beta<\alpha}{\poP_\beta\xmbox{ is a direct limit of }\seqof{\poP_\xi}{\xi<\alpha}}$  
is a stationary subset of $\alpha$ by \xitemof{III:msi-0} and \xitemof{III:msi-8-4}. 
\begin{Lemma}
  \Label{III:P-msi-0-0}  For $\nu\in S\cup\ssetof{\kappa}$, we have $\cardof{\poP_\mu}<\nu$ 
  for all $\mu<\nu$, $\poP_\nu\subseteq V_\nu$ 
  and\/ $\poP_\nu$ has the $\nu$-cc. \qed
\end{Lemma}

For \pos\ $\poP$, $\poQ$, a mapping $\mapping{p}{\poQ}{\poP}$ is said to be a 
{\it projection} if 
\begin{xitemize}
\xitem[III:proj-0] $p(\bbone_\poQ)=\bbone_\poP$;
\xitem[III:proj-1] $p$ is order-preserving; and
\xitem[III:proj-2] for any $\condp\in\poP$ and $\condq\in\poQ$, if
  $\condp\leq_\poP p(\condq)$, then there is $\condq'\in\poQ$ \st\ $\condq'\leq_\poQ\condq$ and
  $p(\condq')\leq_\poP\condp$. 
\end{xitemize}
  
\ifextended{\small\color{darkelectricblue}
  Note that we do not assume that a projection is a surjection. However:
  \begin{LemmaA}
    \Label{III:LA-proj-0} If $\mapping{p}{\poQ}{\poP}$ is a projection then $p\imageof\poQ$ 
    is a dense subset of $\poP$.
  \end{LemmaA}
  \prf For $\condp\in\poP$, we have $\condp\leq_\poP\bbone_{\poP}=q(\bbone_\poQ)$. Thus by 
  \xitemof{III:proj-2}, there is $\condq'\in\poQ$ \st\ $p(\condq')\leq_\poP\condp$. \qedofLemmaA
  \qedskip
}\fi 

The following is standard and also easy to check:
\begin{Lemma}
  \Label{III:P-proj-0}
  Suppose that\/ $\poP$, $\poQ$ are \pos\ and $\mapping{p}{\poQ}{\poP}$ is a 
  projection.\smallskip

  \wassert{1} If\/ $\genH$ is a $(\uniV,\poQ)$-generic filter, then 
  $p\imageof\genH$ generates a $(\uniV,\poP)$-generic filter.\smallskip

  \wassert{2} If\/ $\genG$ is a $(V,\poP)$-generic filter, then letting
  \begin{xitemize}
  \xitem[III:proj-3] $\poQ/\genG=\setof{\condq\in\poQ}{p(\condq)\in\genG}$
  \end{xitemize}
  be \po\ with the pre-ordering $\leq_\poQ$ restricted to it, any
  $(\uniV[\genG],\poQ/\genG)$-generic filter $\genH$ is a $(\uniV,\poQ)$-generic 
  filer with $p\imageof\genH\subseteq\genG$. \qed
\end{Lemma}

Suppose that $\seqof{\poP_\alpha,\utpoQ_\beta}{\alpha\leq\kappa,\beta<\kappa}$ is 
an Easton-type mixed support iteration with the Laver-function 
$\mapping{f}{\kappa}{V_\kappa}$ and $S$ as above over a finite support iteration
$\seqof{\poO_\alpha,\utpoR_\beta}{\alpha\leq\kappa,\beta<\kappa}$. 

Note that, for $\alpha\leq\beta\leq\kappa$,
\begin{xitemize}
\xitem[III:proj-4] 
  $\mapping{p_{\beta,\alpha}}{\poP_\beta}{\poP_\alpha}$; $\condq\mapsto\condq\restr\alpha$
  is a projection, and 
\xitem[III:proj-5] $p_{\beta,\alpha}\circ i_{\alpha,\beta}=id_{\poP_\alpha}$. 
\end{xitemize}

For $\delta_0<\kappa$, let $\genG_{\delta_0}$ be 
a $(\uniV,\poP_{\delta_0})$-generic filter. 
Working in $\uniV[\genG_{\delta_0}]$, 
let $\delta_0\leq\gamma\leq\kappa$, and let
\begin{xitemize}
\xitem[III:msi-16-0] 
  $\poP_\gamma/\genG_{\delta_0}=\setof{\condp\in\poP_\gamma}{\condp\restr\delta_0\in\genG_{\delta_0}}$
\end{xitemize}
be the \po\ with the pre-ordering $\leq_{\poP_\gamma}$ restricted to
$\poP_\gamma/\genG_{\delta_0}$ and with the designated maximal element
$\bbone_{\poP_\gamma/\genG_{\delta_0}}=\bbone_{\poP_\gamma}$.

\begin{Lemma}
  \Label{III:P-msi-1}\wassert{1} A $(\uniV[\genG_{\delta_0}],\poP_\gamma/\genG_{\delta_0})$-generic 
  filter $\genH$ is also a $(\uniV,\poP_\gamma)$-generic filter with 
  $i_{\delta_0,\gamma}\imageof\genG_{\delta_0}\subseteq\genH$. \smallskip

  \wassert{2} If\/ $\genH$ is a $(\uniV,\poP_\gamma)$-generic filter with
  $i_{\delta_0,\gamma}\imageof\genG\subseteq\genH$, then $\genH$ is a
  $(\uniV[\genG_{\delta_0}],\poP_\gamma/\genG_{\delta_0})$-generic filter.
\end{Lemma}
\prf By \xitemof{III:proj-4}, \xitemof{III:proj-5} and \Lemmaof{III:P-proj-0}.\qedofLemma
\qedskip

It is well-known that
projections and complete embeddings are two interchangable notions for cBa 
\,\footnote{We call a \po\ $\poP=\pairof{\poP,\leq_\poP}$ a {\it cBa \po}\/ 
  if (the underlying set) $\poP$ of the \po\ coincides with the positive elements of a 
  complete Boolean algebra and $\leq_\poP$ coincides with the ordering of the complete 
  Boolean algebra.}:

\begin{Lemma}
  \Label{III:P-msi-1-a}
   For cBa \pos\ $\poP$  and $\poQ$, there is a complete embedding 
   $\mapping{i}{\poP}{\poQ}$ if and only if there is a projection $\mapping{p}{\poQ}{\poP}$.

   For cBa \pos\ complete embeddings are injections and projections are surjections. 
   \ifextended
   \else
   \qed\fi
\end{Lemma}
\ifextended{\small\color{darkelectricblue}
\prf 
Suppose that $\poP=\BaA^+$ and $\poQ=\BaB^+$ for complete Boolean algebras $\BaA$ and $\BaB$.

If $\mapping{i}{\poP}{\poQ}$ is a complete embedding, then $\mapping{p}{\poQ}{\poP}$ 
defined by $p(\condb)=\prod^\BaA\setof{\conda\in\poP}{i(\conda)\geq_\poQ\condb}$ for
$\condb\in\poQ$ is a projection.

If $\mapping{p}{\poQ}{\poP}$ is a projection, then the mapping $\mapping{i}{\poP}{\poQ}$ 
defined by $i(\conda)=\sum^{\BaB}\setof{\condb\in\poQ}{p(\condb)\leq_\poP\conda}$
for $\conda\in\poP$ is a complete embedding. Note that $i(\conda)\in\poQ$ by \LemmaAof{III:LA-proj-0}.
\qedofLemma\qedskip}\fi 

For the analysis of the structure of the iteration
\vspace*{-0.4ex}$\seqof{\poP_\alpha,\utpoQ_\beta}{\alpha\leq\kappa,\beta<\kappa}$, the 
following alternative treatment of 
the quotient $\poP_\gamma/\genG_{\delta_0}$ proves often to be more appropriate. 

For $\condp_0$, $\condp_1\in\poP_\gamma$, 
with $\supp(\condp_0)\cap\supp(\condp_1)=\emptyset$, we denote with\\
$\condp_0\land_{\poP_\gamma}\condp_1$, the element $\poP_\gamma$ defined by:
\begin{xitemize}
\xitem[III:msi-+18-a] 
  $\condp_0\land_{\poP_\gamma}\condp_1
  =\condp_0\restr\supp(\condp_0)\cup\condp_1\restr(\gamma\setminus\supp(\condp_0))$. 
\end{xitemize}

If it causes no confusion, we drop the subscript $\poP_\gamma$ in this notation and simply 
write $\condp_0\land\condp_1$ in place of $\condp_0\land_{\poP_\gamma}\condp_1$.

Suppose $\delta_0<\gamma\leq\kappa$ and $\genG_{\delta_0}$ is 
a $(\uniV,\poP_{\delta_0})$-generic filter. 
In $\uniV[\genG_{\delta_0}]$, 
let 
$\poP_\gamma\,|\,\genG_{\delta_0}=\setof{\condp\in\poP_\gamma}{\supp(\condp)\subseteq\gamma\setminus\delta_0}$ be 
the \po\ with the pre-ordering $\leq_{\poP_\gamma\,|\,\genG_{\delta_0}}$ defined by
\begin{xitemize}
\xitem[III:msi-28] 
  $\condq_0\leq_{\poP_\gamma\,|\,\genG_{\delta_0}}\condq_1$\ \ $\Leftrightarrow$\\[\jot]
  $i_{\delta_0,\gamma}(\condp)\land_{\poP_\gamma}\condq_0
  \leq_{\poP_\gamma}i_{\delta_0,\gamma}(\condp)\land_{\poP_\gamma}\condq_1$ for some
  $\condp\in\genG_{\delta_0}$ 
\end{xitemize}
for $\condq_0$, $\condq_1\in\poP_\gamma\,|\,\genG_{\delta_0}$, and with the designated 
maximal element $\bbone_{\poP_\gamma\,|\,\genG_{\delta_0}}=\bbone_{\poP_\gamma}$. 

Note that, for $\condq_0$, $\condq_1\in\poP_\gamma\,|\,\genG_{\delta_0}$,
\begin{xitemize}
\xitem[III:msi-28-0] 
    $\condq_0\leq_{\poP_\gamma}\condq_1$ implies
  $\condq_0\leq_{\poP_\gamma\,|\,\genG_{\delta_0}}\condq_1$, 
\end{xitemize}
since
$\bbone_{\poP_{\delta_0}}\in\genG_{\delta_0}$.

In the following, just for convenience, we shall often misuse the notation and write instead of
$i_{\delta_0,\gamma}(\condp)\land_{\poP_\gamma}\condq_0$ etc. simply $\condp\land\condq_0$ 
etc. The following Lemma is also formulated in this sloppy handling of the notation.

\begin{Lemma}
\Label{III:P-msi-1-0} \wassertof{1} For $\condp\in\poP_{\delta_0}$ and $\condq_0$,
$\condq_1\in\poP_\gamma\,|\,\genG_{\delta_0}$ If
$\condp\land\condq_0\leq_{\poP_\gamma}\condp\land\condq_1$, then for
$\condp'\leq_{\poP_{\delta_0}}\condp$, we have
$\condp'\land\condq_0\leq_{\poP_\gamma}\condp'\land\condq_1\leq_{\poP_\gamma}\condp\land\condq_1$. \smallskip

\wassert{2} For $\condq_0$, $\condq_1\in\poP_\gamma\,|\,\genG_{\delta_0}$ with
$\supp(\condq_0)\cap\supp(\condq_1)=\emptyset$, $\condq_0\land\condq_1$ is a join of 
$\condq_0$ and $\condq_1$  both \wrt\ 
$\leq_{\poP_\gamma}$ and \wrt\ $\leq_{\poP\,|\,\genG_{\delta_0}}$. 
\end{Lemma}
\prf \assertof{1}: By \Lemmaof{III:P-msi-0-a-1}. 
\smallskip

\assertof{2}: 
$\condq_0\land\condq_1$ is a join of 
$\condq_0$ and $\condq_1$  \wrt\ $\leq_{\poP_\gamma}$ by \Lemmaof{III:P-msi-0-a-1}. By 
\xitemof{III:msi-28-0}, it follows that
$\condq_0\land\condq_1\leq_{\poP_\gamma\,|\,\genG_{\delta_0}}\condq_0$, $\condq_1$. 
Suppose now that $\condr\leq_{\poP_\gamma\,|\,\genG_{\delta_0}}\condq_0$, $\condq_1$. Then 
there are $\conds_0$, $\conds_1\in\genG_{\delta_0}$ \st\
$\conds_0\land\condr\leq_{\poP_\gamma}\conds_0\land\condq_0$ and
$\conds_1\land\condr\leq_{\poP_\gamma}\conds_1\land\condq_1$. Let 
$\conds_2\in\genG_{\delta_0}$ be \st\ $\conds_2\leq_{\poP_{\delta_0}}\conds_0$, $\conds_1$. 
Then, by \assertof{1}, we have $\conds_2\land\condr\leq_{\poP_\gamma}\conds_2\land\condq_0$,
$\conds_2\land\condq_1$. By \Lemmaof{III:P-msi-0-a-1}, it follows that
$\conds_2\land\condr\leq_{\poP_\gamma}\conds_2\land(\condq_0\land\condq_1)$. Thus
$\condr\leq_{\poP_\gamma\,|\,\genG_{\delta_0}}\condq_0\land\condq_1$. 
\qedofLemma\qedskip

For $\gamma\leq\kappa$, $\condp\in\poP_\gamma$ and $X\subseteq\kappa$, let 
$\condp\dhrpr X$ be the condition $\condr\in\poP_\gamma$ defined by
\begin{xitemize}
\xitem[III:msi-18-0] 
  $\condr(\alpha)=\left\{
    \begin{array}{@{}ll}
      \condp(\alpha), &\mbox{if }\alpha\in X;\\[\jot]
      \bbone_{\utpoQ_\alpha}, &\mbox{otherwise}
    \end{array}
    \right.$
\end{xitemize}
for all $\alpha<\gamma$.

Since  
$\supp(\condp\dhrpr X)\subseteq\supp(\condp)$, we have 
$\condp\dhrpr X\in\poP_\gamma$ by \Lemmaof{III:P-msi-0-a-0}. 
By definition, it is also clear that $\condp\leq_{\poP_\gamma}\condp\dhrpr X$. 

For $X\subseteq\gamma$ and $P\subseteq\poP_\gamma$, let us write
\begin{xitemize}
\xitem[III:msi-28-1] 
  $P\dhrpr X=\setof{\condp\dhrpr X}{\condp\in P}$. 
\end{xitemize}
Note that the underlying set of $\poP_\gamma\,|\,\genG_{\delta_0}$ could be also 
described as\\ $\poP_\gamma\dhrpr(\gamma\setminus\delta_0)$ 
with this notation.

The \po\ $\poP_\gamma\,|\,\genG_{\delta_0}$ is forcing equivalent to
$\poP_\gamma/\genG_{\delta_0}$.

\begin{Lemma}
  \Label{III:P-msi-3}
  The mapping
  \begin{xitemize}
  \xitem[III:msi-29]
  $\mapping{i_\dhrpr}{\poP_\gamma/\genG_{\delta_0}}{\poP_\gamma\,|\,\genG_{\delta_0}}\ ;$
    \ \ $\condq\mapsto\condq\dhrpr(\gamma\setminus\delta_0)$
  \end{xitemize}
is a dense embedding. 
\end{Lemma}
\prf $i_\dhrpr$ is surjective: If $\condp\in\poP_\gamma\,|\,\genG_{\delta_0}$ then 
$\condp\in\poP_\gamma/\genG_{\delta_0}$ and $i_\dhrpr(\condp)=\condp$.

$i_\dhrpr(\bbone_{\poP_\gamma/\genG_{\delta_0}})=i_\dhrpr(\bbone_{\poP_\gamma})
=\bbone_{\poP_\gamma}\dhrpr(\gamma\setminus\delta_0)=\bbone_{\poP_\gamma}
=\bbone_{\poP_\gamma\,|\,\genG_{\delta_0}}$.\smallskip

$i_\dhrpr$ is order preserving: Suppose that $\condq_0$,
$\condq_1\in\poP_\gamma/\genG_{\delta_0}$ and $\condq_0\leq_{\poP_\gamma}\condq_1$. 
Then $\condq_0\restr\delta_0$,
$\condq_1\restr\delta_0\in\genG_{\delta_0}$ and
$\condq_0\restr\delta_0\leq_{\poP_{\delta_0}}\condq_1\restr\delta_0$. It follows that 
\begin{xitemize}
\xitem[] 
  $\condq_0\restr\delta_0\land(\underbrace{\condq_0\dhrpr(\gamma\setminus\delta_0)}_{=i_\dhrpr(\condq_0)})
  \leq_{\poP_\gamma}
  \condq_0\restr\delta_0\land(\underbrace{\condq_1\dhrpr(\gamma\setminus\delta_0)}_{=i_\dhrpr(\condq_1)})$. 
\end{xitemize}
Thus, $i_\dhrpr(\condq_0)\leq_{\poP_\gamma\,|\,\genG_{\delta_0}}i_\dhrpr(\condq_1)$. \smallskip

$i_\dhrpr$ is incompatibility preserving: Suppose that $\condq_0$,
$\condq_1\in\poP_\gamma/\genG_{\delta_0}$ and, $i_\dhrpr(\condq_0)$ and $i_\dhrpr(\condq_1)$ 
are compatible in $\poP_\gamma\,|\,\genG_{\delta_0}$. Then, there is
$\condr\in\poP_\gamma\,|\,\genG_{\delta_0}$ \st\
$\condr\leq_{\poP_\gamma\,|\,\genG_{\delta_0}}i_\dhrpr(\condq_0)$, $i_\dhrpr(\condq_1)$. By 
the definition of $\leq_{\poP_\gamma\,|\,\genG_{\delta_0}}$, this means that there are $\conds_0$,
$\conds_1\in\genG_{\delta_0}$ \st\
$\conds_0\land\condr\leq_{\poP_\gamma}\conds_0\land i_\dhrpr(\condq_0)$ and
$\conds_1\land\condr\leq_{\poP_\gamma}\conds_1\land i_\dhrpr(\condq_1)$. 

Let $\conds_2\in\genG_{\delta_0}$ be \st\
$\conds_2\leq_{\poP_{\delta_0}}\conds_0$, $\conds_1$, $\condq_0\restr\delta_0$,
$\condq_1\restr\delta_0$. Then
$\conds_2\land\condr\in\poP_\gamma/\genG_{\delta_0}$ and
$\conds_2\land\condr\leq_{\poP_\gamma}\condq_0$, $\condq_1$ by \Lemmaof{III:P-msi-1-0},\,\assertof{1}.
\qedofLemma
\qedskip

Working further in $\uniV[\genG_{\delta_0}]$, 
let 
\begin{xitemize}
\xitem[III:msi-17] 
  $\poS_{\delta_0,\gamma}=(\poP_\gamma\,|\,\genG_{\delta_0})\dhrpr S
  =\setof{\condp\in\poP_\gamma}{\supp(\condp)\subseteq S\setminus\delta_0}$
\end{xitemize}
be the \po\ with the pre-ordering $\leq_{\poP_\gamma\,|\,\genG_{\delta_0}}$ restricted to it and with the 
designated maximal element $\bbone_{\poS_{\delta_0,\gamma}}=\bbone_{\poP_\gamma}$. 
Let  
\begin{xitemize}
\xitem[III:msi-18] 
  $\poT_{\delta_0,\gamma}=(\poP_\gamma\,|\,\genG_{\delta_0})\dhrpr T
  =\setof{\condp\in\poP_\gamma}{\supp(\condp)\subseteq T\setminus\delta_0}$
\end{xitemize}
be the \po\ with the pre-ordering $\leq_{\poP_\gamma\,|\,\genG_{\delta_0}}$ restricted to it and with the 
designated maximal element $\bbone_{\poT_{\delta_0,\gamma}}=\bbone_{\poP_\gamma}$.

\memo{3.14 in Ottenbreit}
\begin{Lemma}
  \Label{III:P-msi-2}
  In $\uniV[\genG_{\delta_0}]$, the mapping
  \begin{xitemize}
  \xitem[III:msi-19] 
  $\mapping{\pi_{\delta_0,\gamma}}{\poS_{\delta_0,\gamma}\times\poT_{\delta_0,\gamma}}{
    \poP_\gamma\,|\,\genG_{\delta_0}}\,;$\ 
  	\ $\pairof{\conds,\condt}\mapsto \conds\,\land\condt$
  \end{xitemize}
  is a projection.


\end{Lemma}
\prf 
$\pi_{\delta_0,\gamma}\models\mbox{\xitemof{III:proj-0}}$ is clear by the definition of $\pi_{\delta_0,\gamma}$.
\smallskip

To show that
$\pi_{\delta_0,\gamma}$ is order-preserving, suppose that 
$\conds'\leq_{\poP_\gamma\,|\,\genG_{\delta_0}}\conds$ and
$\condt'\leq_{\poP_\gamma\,|\,\genG_{\delta_0}}\condt$. Then, there are $\condu_0$, 
$\condu_1\in\genG_{\delta_0}$ \st\
$\condu_0\land\conds'\leq_{\poP_\gamma}\condu_0\land\conds$ and
$\condu_1\land\condt'\leq_{\poP_\gamma}\condu_1\land\condt$.

Let $\condu_2\in\genG_{\delta_0}$ be \st\ $\condu_2\leq_{\poP_{\delta_0}}\condu_0$,
$\condu_1$. By \Lemmaof{III:P-msi-1-0},\,\assertof{1}, we have 
$\condu_2\land\conds'\leq_{\poP_\gamma}\condu_2\land\conds$ and
$\condu_2\land\condt'\leq_{\poP_\gamma}\condu_2\land\condt$.

By \Lemmaof{III:P-msi-0-a-1}, it follows that
$\condu_2\land(\conds'\land\condt')\leq_{\poP_\gamma}\condu_2\land(\conds\land\condt)$. 
Thus, 
$\pi_{\delta_0,\gamma}(\pairof{\conds',\condt'})=\conds'\land\condt'\leq_{\poP_\gamma\,|\,\genG_{\delta_0}}
\conds\land\condt=\pi_{\delta_0,\gamma}(\pairof{\conds,\condt})$.\smallskip

To show 
that $\pi_{\delta_0,\gamma}$ also satisfies \xitemof{III:proj-2}, suppose that
$\pairof{\conds,\condt}\in\poS_{\delta_0,\gamma}\times\poT_{\delta_0,\gamma}$ and 
$\condp\in\poP_\gamma\,|\,\genG_{\delta_0}$ are \st\
\begin{xitemize}
\xitem[III:msi-23-0] 
  $\condp\leq_{\poP_\gamma\,|\,\genG_{\delta_0}}\conds\land\condt
  =\tau_{\delta_0,\gamma}(\pairof{\conds,\condt})$.
\end{xitemize}

Let $\condu\in\genG_{\delta_0}$ be \st\
$\condu\land\condp\leq_{\poP_\gamma}\condu\land(\conds\land\condt)$. 

Let $\condp_0$ be a $\gamma\setminus\delta_0$-sequence defined by
\begin{xitemize}
\xitem[III:msi-24] $\condp_0(\xi)=\left\{\,{}
  \begin{array}{@{}ll}
    \utcondq_\xi, &\mbox{if }\xi\in\supp(\condp);\\[\jot]
    \bbone_{\utpoQ_\xi}, &\mbox{otherwise}
  \end{array}
\right.$
\end{xitemize}
for $\xi\in\gamma\setminus\delta_0$, where $\utcondq_\xi$ is a canonical $\poP_\xi$-name of an element of
$\utpoQ_\xi$ \st\
\begin{xitemize}
\xitem[III:msi-25] 
    $\condu\land\condp\restr\xi\forces{\poP_\xi}{\utcondq_\xi\equiv\condp(\xi)}$, and
\xitem[III:msi-25-0] 
    $\condp'\forces{\poP_\xi}{\utcondq_\xi\equiv(\conds\,\land\condt)(\xi)}$,\qquad\qquad
  $\begin{array}[t]{@{}l}
    \mbox{for all }\condp'\in\poP_\xi\\
    \mbox{with }\incmptbl{\poP_\xi}{\condp'}{\condu\land\condp\restr\xi}.
  \end{array}$
\end{xitemize}

Note that, by \xitemof{III:msi-25} and \xitemof{III:msi-25-0}, we have
\begin{xitemize}
\xitem[III:msi-25-1] 
  $\forces{\poP_\xi}{\utcondq_\xi\leq_{\utpoQ_\xi} (\conds\land\condt)(\xi)}$ for all
  $\xi\in\gamma\setminus\delta_0$. 
\end{xitemize}

Let $\conds_0=\condp_0\dhrpr S$ and $\condt_0=\condp_0\dhrpr T$. 

By the construction, it is clear that the following Claim holds, and this  
shows that $\pairof{\conds_0,\condt_0}$ is a witness for \xitemof{III:proj-2}.
\begin{Claim}
  \Label{III:Cl-msi-0}\wassertof{a} 
  $\pi_{\delta_0,\gamma}(\pairof{\conds_0,\condt_0})=\condp_0\leq_{\poP_\gamma\,|\,\genG_{\delta_0}}\condp$,\smallskip

  \wassert{b} $\pairof{\conds_0,\condt_0}\leq_{\poS_{\delta_0,\gamma}\times\poT_{\delta_0,\gamma}}
  \pairof{\conds,\condt}$. 
\end{Claim}
\prfofClaim
\assertof{a}: By \xitemof{III:msi-25}. \assertof{b}: By \xitemof{III:msi-25-1}. 
\qedofClaim\\
\qedofLemma
\qedskip

Let $\utgenG_{\delta_0}$ denote the standard $\poP_{\delta_0}$-name of a
$(\uniV,\poP_{\delta_0})$-generic filter. 

\begin{Lemma}
  \Label{III:P-msi-3-0} Suppose that\/ $\utcondq_0$ and $\utcondq_1$ 
  are $\poP_{\delta_0}$-names of elements of\\ $\poP_\gamma\,|\,\utgenG_{\delta_0}$.
  Then, we have 
  \begin{xitemize}
  \xitem[III:msi-29-a-0] 
    $\forces{\poP_{\delta_0}}{\utcondq_0\leq_{\poP_\gamma\,|\,\utgenG_{\delta_0}}\utcondq_1}$ 
    \ \ $\Leftrightarrow$\ \ $\forces{\poP_{\delta_0}}{\vec{\bbone}_{\delta_0}\land\utcondq_0
      \leq_{\poP_\gamma}\vec{\bbone}_{\delta_0}\land \utcondq_1}$. 
  \end{xitemize}
\end{Lemma}
\prf ``$\Leftarrow$'' is trivial since
$\forces{\poP_\delta}{\vec{\bbone}_{\delta_0}\in\utgenG_{\delta_0}}$. \smallskip

``$\Rightarrow$'': Suppose that the left side of \xitemof{III:msi-29-a-0} holds. This means that 
  $\forces{\poP_{\delta_0}}{\exists p\in\utgenG_{\delta_0}\,
  (p\land\utcondq_0\leq_{\poP_\gamma}p\land\utcondq_1)}$. 
By \Lemmaof{III:P-msi-1-0}, it follows that 
\begin{xitemize}
\xitem[III:msi-29-a-1] 
  $\forces{\poP_{\delta_0}}{\exists p\in\utgenG_{\delta_0}\,\forall p'\leq_{\poP_\gamma}p\ 
  (p'\land\utcondq_0\leq_{\poP_\gamma}p'\land\utcondq_1)}$. 
\end{xitemize}
Suppose, toward a contradiction, that 
\begin{xitemize}
\xitem[III:msi-29-a-2] 
  $\notforces{\poP_{\delta_0}}{\vec{\bbone}_{\delta_0}\land\utcondq_0
      \leq_{\poP_\gamma}\vec{\bbone}_{\delta_0}\land \utcondq_1}$. 
\end{xitemize}
Then, there are $\condp_0\in\poP_{\delta_0}$ and $\delta_0\leq\delta<\gamma$ \st, 
for any $\condp\leq_{\poP_{\delta_0}}\condp_0$, 
\begin{xitemize}
\xitem[III:msi-29-a-3] 
  $\condp\forces{\poP_{\delta_0}}{\vec{\bbone}_{\delta_0}\land\utcondq_0\restr\delta
      \leq_{\poP_\delta}\vec{\bbone}_{\delta_0}\land \utcondq_1\restr\delta}$, but
\end{xitemize}
\begin{xitemize}
\xitem[III:msi-29-a-4] 
  $\condp\forces{\poP_{\delta_0}}{\vec{\bbone}_{\delta_0}\land\utcondq_0\restr\delta
  \forces{\poP_\delta}{\utcondq_0(\delta)\not\leq_{\utpoQ_\delta}\utcondq_1(\delta)}}$. 
\end{xitemize}
By \Lemmaof{III:P-msi-1-0},\,\assertof{1}, it follows that, for any
$\condp\leq_{\poP_{\delta_0}}\condp_0$, 
$\condp\forces{\poP_{\delta_0}}{\check{\condp}\varin\utgenG_{\delta_0}\mbox{ and }
    \check{\condp}\land\utcondq_0\not\leq_{\poP_\gamma}\check{\condp}\land\utcondq_1}$.
This is a contradiction to \xitemof{III:msi-29-a-1} by \Lemmaof{III:P-msi-1-0},\,\assertof{1}. 
\qedofLemma

\begin{Lemma}
  \Label{III:P-msi-4} For $\delta_0<\gamma\leq\kappa$ and $(\uniV,\poP_{\delta_0})$-generic 
  filter $\genG_{\delta_0}$, we have
  $\uniV[\genG_{\delta_0}]\modelof{\poS_{\delta_0,\gamma}\mbox{ is }\LT\nu(\delta_0)\mbox{-closed}}$.
\end{Lemma}
\prf In $\uniV$, let $\utpoS_{\delta_0,\gamma}$ be a $\poP_{\delta_0}$-name of
$\poS_{\delta_0,\gamma}$ and 
$\uth$ be a $\poP_{\delta_0}$-name of a descending $\delta$-sequence 
in $\poS_{\delta_0,\gamma}$ for some 
\begin{xitemize}
\xitem[III:msi-29-0] 
  $\delta<\nu(\delta_0)$. 
\end{xitemize}
By \Lemmaof{III:P-msi-3-0}, 
we have
$\forces{\poP_{\delta_0}}{\vec{\bbone}_{\delta_0}\land\uth(\xi)\leq_{\poP_\gamma}\vec{\bbone}_{\delta_0}\land\uth(\eta)}$
for all $\xi<\eta<\delta$. 

Let
\begin{xitemize}
\xitem[III:msi-30] 
  $D=\setof{\alpha<\gamma}{\condr\forces{\poP_{\delta_0}}{{}
    \begin{array}[t]{@{}l}
      \alpha\in\supp(\uth(\xi))
    \mbox{ for some }\xi<\delta}\\\mbox{for some }\condr\in\poP_{\delta_0}}.
    \end{array}$
\end{xitemize}

Since
$\forces{\poP_{\delta_0}}{(\forall\xi<\delta)\supp(\xi)\subseteq S\setminus\nu(\delta_0)}$, 
we have
\begin{xitemize}
\xitem[III:msi-31] 
  $D\subseteq S\setminus\nu(\delta_0)$.
\end{xitemize}
\begin{Claim}
  \Label{III:Cl-msi-1} For any regular $\delta_0\leq\mu\leq\gamma$, we have
  $\cardof{D\cap\mu}<\mu$. Thus, $D\in\calI_\gamma$. 
\end{Claim}
\prfofClaim
By \xitemof{III:msi-31}, it is enough to show the inequality for all regular cardinal $\mu$ with
$\nu(\delta_0)\leq\mu\leq\gamma$. For such $\mu$, we have
$\forces{\poP_{\delta_0}}{\supp(\uth(\xi))\cap\mu\xmbox{ is a bounded subset of }\mu}$ for 
all $\xi<\delta$. Thus 
\begin{xitemize}
\xitem[III:msi-32] 
  $D_{\mu,\xi,\condr}=\setof{\alpha<\mu}{\condr\forces{\poP_{\delta_0}}{\alpha\in\supp(\uth(\xi))}}$
\end{xitemize}
is a bounded subset of $\mu$ for each $\xi<\delta$ and $\condr\in\poP_{\delta_0}$. 
By \xitemof{III:msi-29-0} and \Lemmaof{III:P-msi-0-0} for $\nu=\nu(\alpha)$, it follows 
that $D\cap\mu=\bigcup_{\xi<\delta,\condr\in\poP_{\delta_0}}D_{\mu,\xi,\condr}$ is a 
bounded subset of $\mu$. 
\qedofClaim\qedskip

Now we define, by induction on $\delta_0\leq i\leq\gamma$, $\poP_{\delta_0}$-names
$\utcondp_i$, $i\in\gamma+1\setminus\delta_0$ \st\ 
\begin{xitemize}
\xitem[III:msi-33]
  $\forces{\poP_{\delta_0}}{\utcondp_i\varin\poP_i\dhrpr(S\setminus\delta_0)}$ 
  for all $i\in\gamma+1\setminus\delta_0$\,;
\xitem[III:msi-36] 
  $\forces{\poP_{\delta_0}}{\supp(\utcondp_i)\subseteq\check{D}}$ for all $i\in\gamma+1\setminus\delta_0$\,;
\xitem[III:msi-34] 
  $\forces{\poP_{\delta_0}}{{}
    \begin{array}[t]{@{}l}
      (\seqof{\utcondp_i}{i\in\gamma+1\setminus\delta_0})^\bullet\mbox{ is an increasing sequence }\\
  \mbox{of sequences}}\,;
    \end{array}
    $\\
    and
\xitem[III:msi-35]   
  $\forces{\poP_{\delta_0}}{
  \begin{array}[t]{@{}l}
    \utcondp_i\mbox{ is a lower bound of }
    \seqof{\uth(\xi)\restr i}{\xi<\delta}\\[-0.6\jot]
      \mbox{\wrt\ }\leq_{\poP_i\,|\,\utgenG_{\delta_0}}\!\!\!\!\!}
  \end{array}$\\[\jot]
  for all $i\in\gamma+1\setminus\delta_0$\,.
\end{xitemize}

For $i=\delta_0$, $\condp_i=\emptyset$ will do.

Suppose now that $i$ is a limit ordinal and $\utcondp_j$, $j<i$ has been defined \st\ 
\begin{xitemize}
\xitemciteb[III:msi-33]{$'$}
  $\forces{\poP_{\delta_0}}{\utcondp_j\varin\poP_j\dhrpr(S\setminus\delta_0)}$ 
  for all $j\in i\setminus\delta_0$\,;
\xitemciteb[III:msi-36]{$'$} 
  $\forces{\poP_{\delta_0}}{\supp(\utcondp_j)\subseteq\check{D}}$ for all $j\in i\setminus\delta_0$\,;
\xitemciteb[III:msi-34]{$'$}
  $\forces{\poP_{\delta_0}}{{}
    \begin{array}[t]{@{}l}
      (\seqof{\utcondp_j}{j\in i\setminus\delta_0})^\bullet\mbox{ is an increasing sequence }\\
      \mbox{of sequences}}\,;
    \end{array}
    $\\
    and
\xitemciteb[III:msi-35]{$'$}   
  $\forces{\poP_{\delta_0}}{
  \begin{array}[t]{@{}l}
    \utcondp_j\mbox{ is a lower bound of }
    \seqof{\uth(\xi)\restr j}{\xi<\delta}\\[-0.6\jot]
      \mbox{\wrt\ }\leq_{\poP_j\,|\,\utgenG_{\delta_0}}\!\!\!\!\!}
  \end{array}$\\[\jot]
  for 
  all $j\in i\setminus\delta_0$.
\end{xitemize}

By \Lemmaof{III:P-msi-3-0}, \xitembof{III:msi-35}{$'$} implies
\begin{xitemize}
\xitem[III:msi-37]
  $\forces{\poP_{\delta_0}}{\vec{\bbone}_{\delta_0}\land\utcondp_j
  \leq_{\poP_j}\vec{\bbone}_{\delta_0}\land\uth(\xi)\restr j}$ for all $\xi<\delta$ and $j<i$.
\end{xitemize}

Let $\utcondp_i$ be the $\poP_{\delta_0}$-name \st\ 
\begin{xitemize}
\xitem[III:msi-39] $\forces{\poP_{\delta_0}}{\utcondp_i\equiv\bigcup(\setof{\utcondp_j}{j<i})^\bullet}$. 
\end{xitemize}

We show that $\utcondp_i$ together with $\utcondp_j$, $j<i$ satisfies
\begin{xitemize}
\xitemciteb[III:msi-33]{$''$}
  $\forces{\poP_{\delta_0}}{\utcondp_j\varin\poP_j\dhrpr(S\setminus\delta_0)}$ 
  for all $j\in i+1\setminus\delta_0$\,;
\xitemciteb[III:msi-36]{$''$} 
  $\forces{\poP_{\delta_0}}{\supp(\utcondp_j)\subseteq\check{D}}$ for all $j\in i+1\setminus\delta_0$\,;
\xitemciteb[III:msi-34]{$''$}
  $\forces{\poP_{\delta_0}}{{}
    \begin{array}[t]{@{}l}
      (\seqof{\utcondp_j}{j\in i+1\setminus\delta_0})^\bullet\mbox{ is an increasing sequence }\\
      \mbox{of sequences}}\,;
    \end{array}
    $\\
    and
\xitemciteb[III:msi-35]{$''$}   
  $\forces{\poP_{\delta_0}}{
  \begin{array}[t]{@{}l}
    \utcondp_j\mbox{ is a lower bound of }
    \seqof{\uth(\xi)\restr j}{\xi<\delta}\\[-0.6\jot]
      \mbox{\wrt\ }\leq_{\poP_j\,|\,\utgenG_{\delta_0}}\!\!\!\!\!}
  \end{array}$\\[\jot]
  for 
  all $j\in i+1\setminus\delta_0$.
\end{xitemize}
\xitembof{III:msi-36}{$''$} follows from \xitembof{III:msi-36}{$'$} and 
\xitemof{III:msi-39}. \xitembof{III:msi-33}{$''$} follows from this. 
\xitembof{III:msi-34}{$''$} is clear by \xitemof{III:msi-39} and 
\xitembof{III:msi-35}{$''$} follows from \xitemof{III:msi-37}. 

Finally, suppose that $\condp_j$, $j\leq i$ has been defined for
$\nu(\delta_0)\leq j<\gamma$ in accordance with \xitemof{III:msi-33} $\sim$ 
\xitemof{III:msi-35}. In particular, we have
\begin{xitemize}
\xitemciteb[III:msi-37]{$'$}
  $\forces{\poP_{\delta_0}}{\vec{\bbone}_{\delta_0}\land\utcondp_i
  \leq_{\poP_i}\vec{\bbone}_{\delta_0}\land\uth(\xi)\restr i}$ for all $\xi<\delta$. 
\end{xitemize}

If $i\not\in S$, then let
$\utcondp_{i+1}=(\utcondp_{i+1}\cup\ssetof{\pairof{i,\bbone_{\utpoQ_i}}})^\bullet_{\poP_{\delta_0}}$. 

If $i\in S$, then we have
\begin{xitemize}
\xitem[III:msi-40] 
  $\forces{\poP_i}{\utpoQ_i\mbox{ is }\LT\nu(\delta_0)\mbox{-closed}}$ 
\end{xitemize}
by \assertof{B},\,\assertof{a} and \assertof{b} in the definition of our mixed support iteration.
By \xitembof{III:msi-37}{$'$} and by the choice of $\uth$, we have
\begin{xitemize}
\xitem[III:msi-41] 
  $\forces{\poP_{\delta_0}}{\vec{\bbone}_{\delta_0}\land\utcondp_i
  \forces{\poP_i}{{}
  \begin{array}[t]{@{}l}
    (\seqof{\uth(\xi)(i)}{\xi<\delta})^\bullet_{\poP_i}\\
    \mbox{ is a descending sequence in }\utpoQ_i}}.
  \end{array}$
\end{xitemize}
By \xitemof{III:msi-40}, there is a $\poP_{\delta_0}$-name $\utcondq$ of $\poP_{\delta_0}$-name 
\st\
\begin{xitemize}
\xitem[III:msi-42] 
  $\forces{\poP_{\delta_0}\!\!}{\vec{\bbone}_{\delta_0}\land\utcondp_i
  \forces{\poP_i\!\!}{{}
    \begin{array}[t]{@{}l}
      \utcondq\mbox{\,is a lower bound of\,}(\seqof{\uth(\xi)(i)}{\xi<\delta})^\bullet_{\poP_i}\!\!}}.
    \end{array}
    $  
\end{xitemize}

Let $\utcondp_{i+1}=(\utcondp_i\cup\ssetof{\pairof{i,\utcondq}})^\bullet_{\poP_{\delta_0}}$. 
Similarly to the previous case, we can show that $\utcondp_{i+1}$ together with
$\utcondp_j$, $j\leq i$ satisfies \xitemof{III:msi-33} 
$\sim$ \xitemof{III:msi-35}.
\qedofLemma

\begin{Lemma}
  \Label{III:P-msi-4-0} Suppose that $\poP$ is a \po, $\utpoQ$ a $\poP$-name of a \po\ with 
  \begin{xitemize}
  \xitem[III:msi-42-a-a] 
    $\forces{\poP}{\utpoQ\mbox{ is ccc}}$, 
  \end{xitemize}
  and\/ $\poS$ is a $\sigma$-closed \po. Then we have
  \begin{xitemize}
  \xitem[III:msi-42-a] 
    $\forces{\poS}{\forces{\check{\poP}}{\check{\utpoQ}\mbox{ is ccc}}}$.
  \end{xitemize}
\end{Lemma}
\prf Suppose that $\utS$ is a $\poS$-name of a $\poP$-name \st
\begin{xitemize}
\xitem[III:msi-42-a-0]
  $\forces{\poS}{\forces{\check{\poP}}{\utS\mbox{ is a subset of }\check{\utpoQ}\mbox{ of cardinality }
    \aleph_1}}$.
\end{xitemize}

We have to show 
\begin{xitemize}
\xitem[III:msi-42-a-0-0] 
  $\forces{\poS}{\forces{\check{\poP}}{\mbox{ there are compatible elements in }\utS}}$.
\end{xitemize}

Let $\utf$ be a $\poS$-name of $\poP$-name \st
\begin{xitemize}
\xitem[III:msi-42-a-1] 
  $\forces{\poS}{\forces{\check{\poP}}{{}
    \begin{array}[t]{@{}l}
      \mapping{\utildef}{\omega_1}{\utS}\mbox{ and}\\\utf
      \mbox{ is an injective enumeration of }\utS}}.
    \end{array}$
\end{xitemize}

Let $\conds\in\poS$ and $\condp\in\poP$ be arbitrary. By $\sigma$-closedness of $\poS$, we 
can find a decreasing sequence $\seqof{\conds_\alpha}{\alpha<\omega_1}$ of elements 
of $\poS$ and a sequence $\utcondq_\alpha$, $\alpha<\omega_1$ of $\poP$-names \st
\begin{xitemize}
\xitem[III:msi-42-a-2] $\conds_0\leq_\poS\conds$, 
\xitem[III:msi-42-a-3] $\conds_\alpha\forces{\poS}{\forces{\check{\poP}}{\smash{\check{\utcondq}}_\alpha\equiv\utf(\alpha)}}$.
\end{xitemize}
By \xitemof{III:msi-42-a-0}, \xitemof{III:msi-42-a-1} and \xitemof{III:msi-42-a-3}, we have 
\begin{xitemize}
\xitem[III:msi-42-a-3-0] 
  $\conds_\alpha\forces{\poS}{\forces{\check{\poP}}{\smash{\check{\utcondq}}_\alpha\varin\check{\utpoQ}}}$.
\end{xitemize}
Since the relation $\cdot\forces{\cdot}{\cdot\varin\cdot}$ is $\Delta_1$, it follows that 
\begin{xitemize}
\xitem[III:msi-42-a-4] $\forces{\poP}{\utcondq_\alpha\varin\utpoQ}$.   
\end{xitemize}

By \xitemof{III:msi-42-a-a}, there are $\condp'\leq_\poP\condp$ and $\alpha_0<\alpha_1<\omega_1$ \st\ 
\begin{xitemize}
\xitem[III:msi-42-a-5] 
  $\condp'\forces{\poP}{\cmptbl{\utpoQ}{\utcondq_{\alpha_0}}{\utcondq_{\alpha_1}}}$. 
\end{xitemize}
By \xitemof{III:msi-42-a-1} and \xitemof{III:msi-42-a-3}, and since
$\seqof{\conds_\alpha}{\alpha<\omega_1}$ is decreasing, 
\begin{xitemize}
\xitem[III:msi-42-a-6] 
  $\conds_{\alpha_1}\forces{\poS}{\condp'\forces{\poP}{\cmptbl{\utpoQ}{\utf(\alpha_0)}{\utf(\alpha_1)}}}$.
\end{xitemize}
Thus
\begin{xitemize}
\xitem[III:msi-42-a-7] 
  $\conds_{\alpha_1}\forces{\poS}{\exists x\leq_{\check\poP}\condp\
  x\forces{\check{\poP}}{\mbox{ there are compatible elements in }\utS}}.$
\end{xitemize}

$\conds_{\alpha_1}\leq_\poS\conds$ by \xitemof{III:msi-42-a-2}. 
Since $\conds$ was arbitrary, if follows that 
\begin{xitemize}
\xitem[III:msi-42-a-8] 
  $\forces{\poS}{\exists x\leq_{\check\poP}\check{\condp}\
  x\forces{\check{\poP}}{\mbox{ there are compatible elements in }\utS}}.$
\end{xitemize}

Now, 
since $\condp$ was arbitrary, \xitemof{III:msi-42-a-0-0} follows.
\qedofLemma

\begin{Lemma}
  \Label{III:P-msi-5}\wassertof{1}
  For $\delta\leq\kappa$, $\iota_\delta$ is an isomorphism from $\poO_{\delta}$ to
  $\poT_{0,\delta}$. $\seqof{i_\beta}{\beta\leq\delta}$ forms a commutative system together 
  with $\seqof{\poO_\beta, i^*_{\beta,\gamma}}{\beta\leq\gamma\leq\delta}$ and
  $\seqof{\poT_{0,\beta},i_{\beta,\gamma}\restr\poT_\beta}{\beta\leq\gamma\leq\delta}$. In 
  particular, $\seqof{\poT_{0,\beta}}{\beta\leq\delta}$ is homomorphic to the sequence of iterands of a finite support 
  iteration of ccc \pos. \smallskip

  \wassert{2}
  For $\delta_0<\gamma\leq\kappa$ and $(\uniV,\poP_{\delta_0})$-generic 
  filter $\genG_{\delta_0}$, we have
  \begin{xitemize}
  \xitem[III:msi-42-0] 
    $\uniV[\genG_{\delta_0}]\modelof{\poT_{\delta_0,\gamma}\mbox{ has the ccc}}$. 
  \end{xitemize}
\end{Lemma}
\prf \assertof{1}: By induction $\delta\leq\kappa$. \smallskip

\assertof{2}: Note first that, by \Lemmaof{III:msi-2},
$\mapping{\pi_{0,\delta_0}}{\poS_{0,\delta_0}\times\poT_{0,\delta_0}}{\poP_{\delta_0}}$ is 
a projection. Let $\genG_\poS\ast\genG_\poT$ be a $(\uniV,\poS_{0,\delta_0}
\times\poT_{0,\delta_0}/\genG_{\delta_0})$-generic filter in the sense of 
\Lemmaof{III:P-proj-0},\,\assertof{2} where we assume that $\genG_\poS$ and $\genG_\poT$ 
are the generic filters over $\poS_{0,\delta_0}$ and $\poT_{0,\delta_0}$ respectively. 
We have $\genG_\poT=\genG_{\delta_0}\downharpoonright T$. Thus
$\poT_{\delta_0,\gamma}=\poT_{0,\gamma}|\genG_\poT$. By the Factor Lemma for finite 
support iteration of ccc \pos, we have
$\uniV[\genG_\poT]\models\poT_{\delta_0,\gamma}\mbox{ is ccc}$.

Since $\poS_{0,\delta_0}$ is $\sigma$-closed by \Lemmaof{III:P-msi-4}, 
$\uniV[\genG_\poS][\genG_\poT]\models\poT_{\delta_0,\gamma}\xmbox{ is ccc}$ by 
\Lemmaof{III:P-msi-4-0}. 
Since $\uniV[\genG_{\delta_0}]$ is an inner model of $\uniV[\genG_\poS][\genG_\poT]$, it 
follows that $\uniV[\genG_{\delta_0}]\models\poT_{\delta_0,\gamma}\xmbox{ is ccc}$.
\qedofLemma
\qedskip

Summarizing what we have proved above, we obtain the following: 
\begin{Prop}
  \Label{III:P-msi-6}Suppose that $\kappa$ is a supercompact cardinal, 
  $\mapping{f}{\kappa}{\kappa}$ a Laver function with $S$, $T\subseteq\kappa$ defined by 
  \xitemof{III:msi-0}, \xitemof{III:msi-1-a}, and let $\mapping{\nu}{\kappa}{\kappa}$ be 
  defined by \xitemof{III:msi-1-0}.

  For the preparatory finite support ccc iteration
  $\seqof{\poO_\alpha,\utpoR_\beta}{\alpha\leq\kappa,\beta<\kappa}$ satisfying 
  \xitemof{III:msi-2}, \xitemof{III:msi-3}, let
  $\seqof{\poP_\alpha,\utpoQ_\beta}{\alpha\leq\kappa,\beta<\kappa}$ be the Easton-type 
  mixed support iteration over
  $\seqof{\poO_\alpha,\utpoR_\beta}{\alpha\leq\kappa,\beta<\kappa}$ as defined in 
  \assertof{A} and \assertof{B} on pages \pageref{III:msi-8-5} $\sim$ 
  \pageref{III:msi-16} with the complete embeddings $i_{\delta,\gamma}$, $\iota_\delta$ and 
  projections $p_{\delta,\gamma}$ for $\delta\leq\gamma\leq\kappa$. 

  For any $(\uniV,\poP_\kappa)$-generic filter $\genG_\kappa$, $\delta_0<\kappa$ and
  $\genG_{\delta_0}=p_{\kappa,\delta_0}\imageof\genG_\kappa$, there are \pos\ $\poS$, 
  $\poT$, a regular sub\po\ $\poQ$ of the completion of $\poS\times\poT$ in $\uniV[\genG_{\delta_0}]$ \st\ 
  \begin{xitemize}
  \xitem[III:msi-44] 
    $\uniV[\genG_{\delta_0}]\modelof{\poS\mbox{ is }\nu(\delta_0)\mbox{-closed and }
    \poT\mbox{ is ccc}}$, and 
  \xitem[III:msi-45] 
  	$\uniV[\genG_{\delta_0}]\modelof{\poQ\sim\poP_\kappa/\genG_{\delta_0}}$. 
  \end{xitemize}
  In particular, 
  there is a $(\uniV[\genG_{\delta_0}],RO(\poS\times\poT))$-generic 
  filter $\tilde{\genH}$ \st, letting $\genH=\tilde{\genH}\cap\poQ$, 
  we have  $\uniV[\genG_{\kappa}]=\uniV[\genG_{\delta_0}][\genH]$. 
\end{Prop}
\prf By \Lemmaof{III:P-msi-3}, we have
$\uniV[\genG_{\delta_0}]\models\poP_\kappa/\genG_{\delta_0}\,\sim\,\poP_\kappa|\genG_{\delta_0}$. 

In $V[\genG_{\delta_0}]$, $\poP_\kappa|\genG_{\delta_0}$ is forcing equivalent to a regular 
sub-\po\ of the completion of $\poS_{\delta_0,\gamma}\times\poT_{\delta_0,\gamma}$ 
by \Lemmaof{III:P-msi-2} (c.f.\ \Lemmaof{III:P-msi-1-a}). By  
\Lemmaof{III:P-msi-4}, $\poS_{\delta_0,\gamma}$ is $\nu(\delta_0)$-closed and by 
\Lemmaof{III:P-msi-5}, $\poT_{\delta_0,\gamma}$ is ccc. 
\qedofProp

\begin{Thm}
  \Label{III:P-msi-7}Suppose that $\kappa$, $f$, $S$, $T$, $\nu$,
  $\seqof{\poO_\alpha,\utpoR_\beta}{\alpha\leq\kappa,\beta<\kappa}$,
  $\seqof{\poP_\alpha,\utpoQ_\beta}{\alpha\leq\kappa,\beta<\kappa}$, $i_{\delta,\gamma}$,
  $\iota_\delta$, 
  $p_{\delta,\gamma}$ for $\delta\leq\gamma\leq\kappa$, and $\genG_\kappa$ are as in \Propof{III:P-msi-6}. 
  \smallskip

  \wassert{0} If\/ $\seqof{\poO_\alpha,\utpoR_\beta}{\alpha\leq\kappa,\beta<\kappa}$ adds 
  $\kappa$ many reals then $\uniV[\genG_\kappa]\models\kappa=\continuum$.\smallskip

  \wassert{1} In $\uniV[\genG_\kappa]$, $\kappa$ is strongly Laver-generically supercompact for 
  $(\gmP,\gmQ_\theta)$ for all $\theta\in\Card$ for the properties of \pos\ $\gmP$ and $\gmQ_\theta$ as in 
  \Propof{III:P-gen-large-3},\,\assertof{2}. \smallskip

  \wassert{1$'$} If $\kappa$ is superhuge, then, in $\uniV[\genG_\kappa]$, $\kappa$ is 
  strongly Laver-generically superhuge for 
  $(\gmP,\gmQ_\theta)$ for all $\theta\in\Card$ for the properties of \pos\ $\gmP$ and $\gmQ_\theta$ as in 
  \Propof{III:P-gen-large-3},\,\assertof{2}. \smallskip

  \wassert{2} If the preparatory iteration
  $\seqof{\poO_\alpha,\utpoR_\beta}{\alpha\leq\kappa,\beta<\kappa}$ is \st\
  $\forces{\poQ_\alpha}{\utpoR_\alpha\equiv\Fn(\omega,2)}$ for all $\alpha\in T$, then 
  in $\uniV[\genG_\kappa]$, $\kappa$ is strongly Laver-generically supercompact for 
  $(\gmP,\gmQ_\theta)$ for all $\theta\in\Card$ for the properties of \pos\ $\gmP$ and $\gmQ_\theta$ as in 
  \Propof{III:P-gen-large-3},\,\assertof{1}. 

  \wassert{2$'$} If the preparatory iteration
  $\seqof{\poO_\alpha,\utpoR_\beta}{\alpha\leq\kappa,\beta<\kappa}$ is \st\
  $\forces{\poQ_\alpha}{\utpoR_\alpha\equiv\Fn(\omega,2)}$ for all $\alpha\in T$ and 
  $\kappa$ is superhuge, then 
  in $\uniV[\genG_\kappa]$, $\kappa$ is strongly Laver-generically superhuge for 
  $(\gmP,\gmQ_\theta)$ for all $\theta\in\Card$ for the properties of \pos\ $\gmP$ and $\gmQ_\theta$ as in 
  \Propof{III:P-gen-large-3},\,\assertof{1}. 
\end{Thm}
\prf \assertof{0}: By \Lemmaof{III:P-msi-0-0}, $\kappa$ is a regular cardinal in 
$\uniV[\genG_\kappa]$ and $\continuum\leq\kappa$ in $\uniV[\genG_\kappa]$. Since 
$\mapping{\iota_\kappa}{\poO_\kappa}{\poP_\kappa}$ is a complete embedding, 
if $\poO_\kappa$ adds $\kappa$ many reals then $\kappa\geq\continuum$.
\smallskip

\assertof{1}: In $V[\genG_\kappa]$, let $\poP=\Col(\kappa,\mu)$ for some cardinal $\mu$ 
and, let $\lambda$ and 
$\theta$ be two other cardinals. \Wolog, we may assume that $\mu^{\LT\kappa}\leq\lambda$. 
Since $f$ is a Laver function there is an elementary embedding
$\elembed{j}{\uniV}{M}$ \st\ 
\begin{xitemize}
\xitem[III:msi-46] 
  $\crit(j)=\kappa$,
\xitem[III:msi-47] 
  $j(\kappa)>\lambda$,
\xitem[III:msi-48] 
  $[M]^\lambda\subseteq M$, and 
\xitem[III:msi-49] 
$j(f)(\kappa)=\pairof{\mu,\theta,\emptyset}$. 
\end{xitemize}

Let $S^*=j(S)$, $\nu^*=j(\nu)$ and, let
$\vec{\poP}^*=j(\seqof{\poP_\alpha,\utpoQ_\beta}{\alpha\leq\kappa,\beta<\kappa})$. 
Since $\vec{\poP}^*$ is $j(\kappa)$-(double) sequence by elementarity, we write
\begin{xitemize}
\xitem[III:msi-50] 
  $\vec{\poP}^*=\seqof{\poP^*_\alpha,\utpoQ^*_\beta}{\alpha\leq j(\kappa),\beta<j(\kappa)}$.
\end{xitemize}

By the elementarity of $j$, \Lemmaof{III:P-msi-0-0}, \xitemof{III:msi-46} and 
\xitemof{III:msi-48}, we have $\poP^*_\alpha=\poP_\alpha$ 
and $\utpoQ^*_\beta=\utpoQ_\beta$ for all $\alpha$, $\beta<\kappa$. 
$\poP^*_\kappa=\poP_\kappa$ by elementarity and \xitemof{III:msi-48}. 
Note that $\kappa\in S^*$ by elementarity. 
So it follows that 
$\utpoQ^*_\kappa=(\Col(\kappa,\mu))^\bullet_{\poP_\kappa}$ by \xitemof{III:msi-10} and \xitemof{III:msi-49}. 
Thus $\utpoQ^*_\kappa[\genG_\kappa]=\poP$. 

Also by \xitemof{III:msi-49}, we have $\nu^*(\kappa)\geq\theta$.
Let $\geng$ be a $(\uniV[\genG_\kappa],\poP)$-generic filter.
In $M[\genG_\kappa][\geng]$, $\poP^*_{j(\kappa)}/\genG_\kappa\ast\geng$ is forcing equivalent to a 
regular sub-\po\ of the completion of a \po\ of the form ``ccc \po\ $\times$ $\theta$-closed 
\po'' by \Propof{III:P-msi-6}. Thus, $\poP^*_{j(\kappa)}/\genG_\kappa\ast\geng\models\gmQ_\theta$. 

Let $\genH^*$ be a
$(\uniV[\genG_\kappa\ast\geng],\poP^*_{j(\kappa)}/\genG_\kappa\ast\geng)$-generic filter.
Then we can find a $(\uniV, \poP^*_{j(\kappa)})$-generic filter $\genH$ \st\
$M[\genH]=M[\genG_\kappa\ast\geng\ast\genH^*]$ and $i^*_{\kappa,j(\kappa)}\imageof\genG_\kappa\subseteq\genH$ for the 
complete embedding $\combed{i^*_{\kappa,j(\kappa)}}{\poP_\kappa}{\poP^*_{j(\kappa)}}$ 
  associated with $\vec{\poP}^*$. 
$j$ can be then lifted to 
\begin{xitemize}
\xitem[III:msi-51] 
  $\elembed{\tilde{j}}{\uniV[\genG_\kappa]}{M[\genH]}$;
  $\uta[\genG]\mapsto j(\uta)[\genH]$.   
\end{xitemize}
It is easy to show that $\tilde{j}$ with $\genH$ satisfies \xitemof{III:gen-large-1-0} $\sim$ 
\xitemof{III:gen-large-6} and \xitemof{III:gen-large-6-0}. The last condition holds by 
\Lemmaof{III:P-gen-large-a-0}. 
This shows that $\kappa$ is strongly Laver-generically supercompact 
for $(\gmP,\gmQ_\theta)$. 
\smallskip

\assertof{2$'$}: is proved similarly to \assertof{2} above. The condition 
\xitembof{III:gen-large-6-0}{$'$} is shown using \Lemmaof{III:P-gen-large-a-1}. 

\assertof{3}, \assertof{3$'$}: can be proved similarly to \assertof{2} and \assertof{2$'$}. 
\qedofThm

\section{Models with strong reflection properties down to $\LT\continuum$ and 
         with even stronger reflection properties but down to $\LE\continuum$}
\Label{III:refl}
As an application of the forcing constructions considered in the previous sections, 
we give two models of large continuum with strong reflection properties around the 
continuum. In one of the models, we have $\HH(\LT\continuum)$ and, in the other, this reflection 
property is negated. Thus, we obtain the independence of $\HH(\LT\continuum)$ from other 
strong reflection principles in the large continuum context. 

In contrast, the most of the other reflection properties 
are situated in a tight web of implications which is (almost) upward directed (see e.g. the 
diagram in the last section of [Fuchino, Sakai and Ottenbreit\cite{III:bib-I}]). 
This suggests that the reflection of non-metrizability is a totally different kind of 
reflection statement from the other reflection principles. 

With an arbitrary preparatory finite support ccc iteration
$\seqof{\poO_\alpha,\utpoR_\beta}{\alpha\leq\kappa,\beta<\kappa}$
we already have the following:

\begin{Thm}
  \Label{III:P-models-0} Let $\kappa$, $\kappa_1$ with $\kappa<\kappa_1$ be two 
  supercompact cardinals and let 
  $\seqof{\poP_\alpha,\utpoQ_\beta}{\alpha\leq\kappa,\beta<\kappa}$ be the Easton-type 
  mixed support iteration over an arbitrary preparatory finite support ccc iterating
  $\seqof{\poO_\alpha,\utpoR_\beta}{\alpha\leq\kappa,\beta<\kappa}$ which adds $\kappa$ 
  many reals. Let $\poP=\poP_\kappa\ast(\Col(\kappa,\kappa_1))^\bullet_{\poP_\kappa}$. 
  Then,  
  in the $\poP$ generic extension over $\uniV$, we have
  \begin{xitemize}
  \xitem[III:models-0] $\continuum=\kappa$; 
  \xitem[III:models-1] $\SDLS^{int}_+(\calL^{\aleph_0}_{stat},\LT\continuum)$ and
    $\GRP^{\LT\continuum}(\LE\continuum)$. 
  \xitem[III:models-2] $\SDLS^{int}_+(\calL^\PKL_{stat},\LT\continuum)$ and
    $\Pkl{\continuum}{\lambda}$ carries a $\sigma$-saturated normal ideal for all
    $\lambda\geq\continuum$. 
  \end{xitemize}
\end{Thm}
\prf \xitemof{III:models-0}: By \Thmof{III:P-msi-7},\,\assertof{0}. 
\smallskip

\xitemof{III:models-1}: By \Thmof{III:P-msi-7},\,\assertof{1}, $\kappa$ is strongly 
Laver-generically supercompact for $(\gmP,\gmQ)$ for properties $\gmP$, $\gmQ$ as in 
\Corof{III:P-gen-large-2}. Thus, by \Corof{III:P-gen-large-2},\,\assertof{2}, we have 
$\SDLS^{int}_+(\calL^{\aleph_0}_{stat},\LT\continuum)$ and
$\GRP^{\LT\continuum}(\LE\continuum)$ in the generic extension. \smallskip
\smallskip

\xitemof{III:models-2}: By \Thmof{III:P-msi-7},\,\assertof{1}, 
$\kappa$ is strongly 
Laver-generically supercompact for $(\gmP,\gmQ_\theta)$ for any $\theta$ for $\gmP$ and
$\gmQ_\theta$ as in \Propof{III:P-gen-large-3}, \assertof{2}. 
Thus \xitemof{III:models-2} holds by \Propof{III:P-gen-large-3},\,\assertof{2} and \assertof{3}. 
\qedofThm
\qedskip

The following strengthening of the forcing axiom $\MA(\calP)$ for a class $\calP$ of \pos\ 
was studied in [Fuchino, Ottenbreit and Sakai\cite{III:bib-II}]:

For a \po\ $\poP$, $\poP$-name $\utilde{S}$ of a set of subsets of $\On$ and a 
filter $\genG$ on $\poP$, let
\begin{xitemize}
\xitem[III:models-3] 
  $\utilde{S}(\genG)=\setof{b}{{}
  \begin{array}[t]{@{}l}
    b=\setof{\alpha\in\On}{\condp\forces{\poP}{\check{\alpha}\varin\utilde{s}}
    \mbox{ for a }\condp\in\genG}\\\mbox{for a }\poP\mbox{-name }\utilde{s}\mbox{ \st\ }\\
    \forces{\poP}{\utilde{s}\varin\utilde{S}\mbox{ and }\sup(\utilde{s})\equiv\sup(b)}}.
  \end{array}$
\end{xitemize}

Note that if $\genG$ is a $(\uniV,\poP)$-generic filter, then
$\utilde{S}(\genG)=\utilde{S}[\genG]$.
\ifextended{\small\color{darkelectricblue}
  [ ...]
  }\fi

For uncountable cardinals $\mu$ and $\kappa>\aleph_1$, let
$\MA^{++\mu}(\calP,\LT\kappa)$ be the strengthening 
of $\MA^{+\mu}(\calP,\LT\kappa)$ defined by:
\begin{xitemize}
\item[$\MA^{++\mu}(\calP,\LT\kappa)$: ]\it For any $\poP\in\calP$, any family 
  $\calD$ of dense subsets of\/ $\poP$ with $\cardof{\calD}<\kappa$ and any family 
  $\calS$ of\/ $\poP$-names \st\ $\cardof{\calS}\leq\mu$ and
  $\forces{\poP}{\utilde{S}\xmbox{ is a stationary subset of }\Pkl{\eta_{\mbox{\scriptsize\(\scriptstyle\utS\)}}}{\theta_{\utilde{S}}}}$ 
  for some $\omega<\eta_{\utilde{S}}\leq\theta_{\utilde{S}}<\continuum$ with $\eta_{\utilde{S}}$ regular, for all
  $\utilde{S}\in\calS$, there is a $\calD$-generic filter $\genG$ over $\poP$ \st\
  $\utilde{S}(\genG)$ is stationary in $\Pkl{\eta_{\mbox{\scriptsize\(\scriptstyle\utS\)}}}{\theta_{\utilde{S}}}$ for all
  $\utilde{S}\in\calS$. 
\end{xitemize}

In case of $\mu=\omega_1$, the principle $\MA^{++\omega_1}(\calP,\LT\kappa)$ is 
equivalent to the usual $\MA^{+\omega_1}(\calP,\LT\kappa)$.

\begin{Prop}
  \Label{III:P-models-0-0}
  Suppose that $\kappa$ is a supercompact cardinal and $f$ a Laver function on $\kappa$. 
  Let $S$ and $T$ be defined by \xitemof{III:msi-0} and \xitemof{III:msi-1-a}.\smallskip

  \wassert{1} Suppose that the preparatory finite support ccc iteration 
  $\seqof{\poO_\alpha,\utpoR_\beta}{\alpha\leq\kappa,\beta<\kappa}$ is defined by:
  \begin{xitemize}
  \xitem[III:models-4] $\utpoR_\beta=\left\{{}
    \begin{array}{@{}ll}
      (P)^\surd_{\poO_\beta}, &\mbox{if }\beta=\alpha+1\mbox{ for an }\alpha\in S\mbox{ and}\\
      &f(\alpha)=\pairof{\mu,\theta,P}\mbox{ for cardinals }\mu,\,\theta\mbox{ and a \po\ }P\\
      &\mbox{\st\ }\forces{\poO_\beta}{(P)^\surd_{\poO_\beta}\mbox{ is a ccc \po}};\\[\jot]
      \ssetof{\bbone_{\utpoR_\beta}}, &\mbox{otherwise.}
    \end{array}\right.$
  \end{xitemize}
  Then, for the Easton-type mixed support iteration
  $\seqof{\poP_\alpha,\utpoQ_\beta}{\alpha\leq\kappa,\beta<\kappa}$ over
  $\seqof{\poO_\alpha,\utpoR_\beta}{\alpha\leq\kappa,\beta<\kappa}$ and
  $\poP=\poP_\kappa\ast(\Col(\kappa,\mu))^\bullet_{\poP_\kappa}$ for a regular
  $\mu>\kappa$, 
  \begin{xitemize}
  \xitem[III:models-5] 
    $\forces{\poP}{{}
    \begin{array}[t]{@{}l}
      \MA^{++\eta}(\calP,\LT\continuum)\mbox{ for all cardinal }\eta<\continuum\\
    \mbox{holds for }\calP=\setof{P}{P\mbox{ is a ccc \po\ and }P\in\uniV}},\end{array}$
  \end{xitemize}
  where ``$P$ is ccc \po'' is meant ``ccc \po\ in the $\poP$-generic extension'' while 
  $\uniV$ denotes here the ground model before extending generically by $\poP$. \smallskip

  \wassert{2} Suppose that
  $\seqof{\poO_\alpha,\utpoR_\beta}{\alpha\leq\kappa,\beta<\kappa}$ is a preparatory finite 
  support iteration \st\ each $\utpoQ_\alpha$ for $\alpha\in T$ is a $\poO_\alpha$-name of 
  the Cohen \po\ $\Fn(\omega_2)$. Then, 
 for the Easton-type mixed support iteration
  $\seqof{\poP_\alpha,\utpoQ_\beta}{\alpha\leq\kappa,\beta<\kappa}$ over
  $\seqof{\poO_\alpha,\utpoR_\beta}{\alpha\leq\kappa,\beta<\kappa}$ and
  $\poP=\poP_\kappa\ast(\Col(\kappa,\mu))^\bullet_{\poP_\kappa}$ for a regular
  $\mu>\kappa$, 
  \begin{xitemize}
  \xitem[III:models-6] 
    $\forces{\poP}{{}
    \begin{array}[t]{@{}l}
      \MA^{++\eta}(\calP,\LT\continuum)\mbox{ for all cardinal }\eta<\continuum\mbox{ holds for }\\
      \calP=\setof{P}{P\mbox{ is forcing equivalent to }\Fn(\mu,2)
      \mbox{ for some }\mu}},\end{array}$
  \end{xitemize}
\end{Prop}
\prf
\assertof{1}: Suppose that $\genG_\kappa$ is a $(\uniV,\poP_\kappa)$-generic filter and 
$\geng$ a $(\uniV[\genG_\kappa],\Col(\kappa,\mu)^{\uniV[\genG_\kappa]})$-generic filter.
Let $\eta$, $\nu<\kappa$, 
$\poR\in V$, and $\calD$, $\calS$ $\in V[\genG_\kappa\ast\geng]$ be \st\ 
\begin{xitemize}
\xitem[III:models-7] 
  $\uniV[\genG_\kappa][\geng]\modelof{{}
  \begin{array}[t]{@{}l}
    \poR\mbox{ is a ccc \po,}\\
    \calD\mbox{ is a family of dense subsets of }\poR\mbox{ with }\cardof{\calD}=\nu, \\
    \mbox{and }
    \calS\mbox{ is a family of }\poR\mbox{-names with }\cardof{\calS}=\eta\mbox{ such }\\
    \mbox{that each element }\utS\mbox{ of }\calS\mbox{ is a }\poR\mbox{-name 
      of a station-}\\[-1\jot] 
    \mbox{nary subset of }\Pkl{\eta_{\mbox{\scriptsize\(\scriptstyle\utS\)}}}{\theta_{\utS}}\mbox{ for some }
    \omega<\eta_{\utS}\leq\theta_{\utS}<\continuum\\
    \mbox{with }\eta_{\utS}\mbox{ regular}}.
  \end{array}$
\end{xitemize}

Let $\cardof{\poR}=\lambda$. \Wolog, we may assume that the underlying set of $\poR$ is
$\lambda$. Thus $\poR=\pairof{\lambda\leq_\poR}$. Let $\theta$ be sufficiently large and 
let $\elembed{j}{V}{M}$ be \st\ $\crit(j)=\kappa$,   $j(\kappa)>\theta$,
\begin{xitemize}
\xitem[III:models-8] 
  $[M]^\theta\subseteq M$, and 
\xitem[III:models-9] 
  $j(f)(\kappa)=\pairof{\mu,\theta,\poR}$. 
\end{xitemize}

Let $S^*=j(S)$, $\nu^*=j(\nu)$ and, let
$\vec{\poP}^*=j(\seqof{\poP_\alpha,\utpoQ_\beta}{\alpha\leq\kappa,\beta<\kappa})$. 
As before, we write
\begin{xitemize}
\xitem[III:models-10] 
  $\vec{\poP}^*=\seqof{\poP^*_\alpha,\utpoQ^*_\beta}{\alpha\leq j(\kappa),\beta<j(\kappa)}$.
\end{xitemize}
We have $\poP_\alpha=\poP^*_\alpha$ for $\alpha\leq\kappa$.

Let $\genG_\kappa$ be a $(\uniV,\poP_\kappa)$-generic filter. Then 
$\utpoQ^*_\kappa[\genG_\kappa]=\Col(\kappa,\mu)^{\uniV[\genG_\kappa]}$ and
$\nu^*(\kappa)\geq\theta$ by 
\xitemof{III:models-9}. Let $\geng$ be a
$(\uniV[\genG_\kappa],\,\Col(\kappa,\mu)^{\uniV[\genG_\kappa]})$-generic filter. 
By \xitemof{III:models-9} and \xitemof{III:models-7},
$\utpoQ_{\kappa+1}[\genG_\kappa\ast\geng]\sim\poR$.

Thus, in $M$, $\poP^*_{j(\kappa)}$ is factored as
\begin{xitemize}
\xitem[III:models-11] 
  $\poP^*_{j(\kappa)}$\ \ $\sim$\ \ $\poP_\kappa\ast\Col(\kappa,\mu)^\bullet_{\poP_\kappa}\ast\utpoR\ast\utpoR_1$
\end{xitemize}
where $\utpoR$ corresponds to the ground model \po\ $\poR$. We have 
\begin{xitemize}
\xitem[III:models-12] 
  $\forces{\poP_\kappa\ast\Col(\kappa,\mu)^\bullet_{\poP_\kappa}\ast\utpoR}{{}
  \begin{array}[t]{@{}l}
    \utpoR_1\mbox{ is a regular sub-\po\ of the completion of}\\
    \mbox{a \po\ of the form `ccc \po\,}\times\nu^*(\kappa)\mbox{-closed}\\
    \mbox{\po'}\ }
  \end{array}$
\end{xitemize}
by \xitemof{III:P-msi-2}, \xitemof{III:P-msi-4}, and \xitemof{III:P-msi-5}. 

Now, let $\genh$ be $(M[\genG\ast\geng],\poR)$-generic filter and $\genH$ be\vspace{-1\jot}
$(M[\genG\ast\geng\ast\genh],\,\utpoR_1[\genG\ast\geng\ast\genh])$-generic filter.
$\elembed{j}{\uniV}{M}$ is then lifted to 
\begin{xitemize}
\xitem[III:models-13] 
  $\elembed{j^*}{\uniV[\genG_\kappa]}{M[\genG_\kappa\ast\geng\ast\genh\ast\genH]}$;
  $\uta[\genG_\kappa]\mapsto j(\uta)[\genG_\kappa\ast\geng\ast\genh\ast\genH]$.
\end{xitemize}

We have
\begin{xitemize}
\xitem[III:models-13-0] 
  $M\modelof{\poC^*=j^*(\Col(\kappa,\mu)^{\uniV[\genG_\kappa]})\mbox{ is }\LT j(\kappa)\mbox{-directed closed}}$
\end{xitemize}
by elementarity and $(\mu^{<\kappa})^+<\theta<j(\kappa)$. Thus, we can find an
$(M[\genG_\kappa\ast\geng\ast\genh\ast\genH],\,\poC^*)$-generic filter $\geng^*$ \st\
$j^*\imageof\geng\subseteq\geng^*$.
$j^*$ is then further lifted to
\begin{xitemize}
\xitem[III:models-14] 
  $\elembed{j^{**}}{\left(\uniV[\genG_\kappa]\right)[\geng]}{
  \left(M[\genG_\kappa\ast\geng\ast\genh\ast\genH]\right)[\geng^*]}$;
  $\uta[\geng]\mapsto j^*(\uta)[\geng^*]$
\end{xitemize}
for $\Col(\kappa,\mu)^{\uniV[\genG_\kappa]}$-names $\uta$ in $\uniV[\genG_\kappa]$.

In $M[\genG_\kappa\ast\geng\ast\genh]$, $\genh$ is a filter on $\poR$ with intersects with 
each element of $\calD$ and each element $\utS$ of $\calS$ is interpreted as a stationary 
subset of $\Pkl{\eta_{\mbox{\scriptsize\(\scriptstyle\utS\)}}}{\theta_{\utS}}$ by the genericity of $\genh$ and 
since $\calD$, $\calS\in M$ by the closedness \xitemof{III:models-8} of $M$. These 
interpretations of elements of $\calS$ remain stationary in
$M[\genG_\kappa\ast\geng\ast\genh\ast\genH\ast\geng^*]$ by \xitemof{III:models-12}. and 
\xitemof{III:models-13-0}
(see \Lemmaof{III:L-preserv-a}).

Since $\calD$ and $\calS$ have cardinality $<\kappa$, we have
\begin{xitemize}
\xitem[III:models-15] $j^{**}(\calD)=\setof{j^{**}(D)}{D\in\calD}$, and\\
  $j^{**}(\calS)=\setof{j^{**}(\utS)}{\utS\in\calS}$. 
\end{xitemize}
It follows that $\genh$ generates a filter on $j^{**}(\poR)$ which intersects each element of 
$j^{**}(\calD)$ and interprets each element $\utS$ of $j^{**}(\calS)$ as a superset of the 
corresponding interpretation of the element of $\calS$ by $\genh$ is a stationary subset of
$\Pkl{\eta_{\mbox{\scriptsize\(\scriptstyle\utS\)}}}$ in $M[\genG_\kappa\ast\geng\ast\genh\ast\genH\ast\geng^*]$. 

Thus, we have
\begin{xitemize}
\xitem[III:models-16] 
  $M[\genG_\kappa\ast\geng\ast\genh\ast\genH\ast\geng^*]\modelof{{}
  \begin{array}[t]{@{}l}
    \mbox{there is a }j^{**}(\calD)\mbox{-generic filter on }j^{**}(\poR)\\[\jot]
  \mbox{which interprets each element }\utS\mbox{ of }j^{**}(\utS)\\[-0.5\jot]
  \mbox{as a stationary subset of }\Pkl{\eta_{\mbox{\scriptsize\(\scriptstyle\utS\)}}}{\theta_{\utS}}}.
  \end{array}$
\end{xitemize}

By the elementarity of $j^{**}$, it follows that 
\begin{xitemize}
\xitem[III:models-17] 
  $\uniV[\genG_\kappa\ast\geng]\modelof{{}
  \begin{array}[t]{@{}l}
    \mbox{there is a }\calD\mbox{-generic filter on }\poR\\[\jot]
    \mbox{which interprets each element }\utS\mbox{ of }\utS\\[-0.5\jot]
    \mbox{as a stationary subset of }\Pkl{\eta_{\mbox{\scriptsize\(\scriptstyle\utS\)}}}{\theta_{\utS}}}.
  \end{array}$
\end{xitemize}

\assertof{2}: can be proved similarly to \assertof{1}. 
\qedofProp

\begin{Thm}
  \Label{III:P-models-1} \wassertof{1} Suppose that the existence of two supercompact cardinals is 
  consistent. Then the following combination of the principles is also 
  consistent:
  \begin{xitemize}
  \xitemcite[III:models-1] 
    $\SDLS^{int}_+(\calL^{\aleph_0}_{stat},\LT\continuum)$, $\GRP^{\LT\continuum}(\LE\continuum)$; 
  \xitemcite[III:models-2] 
    $\SDLS^{int}_+(\calL^\PKL_\stat,\LT\continuum)$;
  \xitem[III:models-21] 
    $\MA^{++\eta}(\calP,\LT\continuum)$ for
    $\calP=\setof{\poP}{\poP\sim\Fn(\lambda,2)\mbox{ for some }\lambda}$\\for all
    $\eta<\continuum$; and 
  \xitem[III:models-22] 
    $\HH(\LT\continuum)$. 
  \end{xitemize}

  \wassert{2} If there is a superhuge cardinal and a supercompact cardinal above it, then 
  the combination of the principles \xitemof{III:models-1} $\sim$ \xitemof{III:models-22} 
  above together with
  \begin{xitemize}
  \xitem[III:models-22-0] 
    $\Pkl{\continuum}{\lambda}$ carries a $\sigma$-saturated normal ideal for all
    $\lambda\geq2^{\aleph_0}$
  \end{xitemize}
  is consistent. 
\end{Thm}
\prf 
\assertof{1}: For two supercompact cardinals $\kappa<\kappa_1$, let
$\seqof{\poO_\alpha,\utpoR_\beta}{\alpha\leq\kappa,\beta<\kappa}$ and 
$\seqof{\poP_\alpha,\utpoQ_\beta}{\alpha\leq\kappa,\beta<\kappa}$ 
be as in \Propof{III:P-models-0-0},\,\assertof{2}.

Then, the generic extension of $\uniV$ by
$\poP=\poP_\kappa\ast(\Col(\kappa,\kappa_1))^\bullet_{\poP_\kappa}$ is as desired:  
$\forces{\poP}{\mbox{\xitemof{III:models-0}, \xitemof{III:models-1}, \xitemof{III:models-2}}}$
by \Thmof{III:P-models-0}. $\forces{\poP}{\mbox{\xitemof{III:models-21}}}$ by 
\Propof{III:P-models-0-0},\assertof{2} and $\forces{\poP}{\mbox{\xitemof{III:models-22}}}$ 
by \Thmof{III:P-msi-7},\,\assertof{2} and \Propof{III:P-gen-large-3}. 
\smallskip

\assertof{2}: Let $\poP=\poP_\kappa\ast(\Col(\kappa,\kappa_1))^\bullet_{\poP_\kappa}$ be as 
in \assertof{1} but for a superhuge $\kappa$ and a supercompact 
$\kappa_1$ above $\kappa$. Then
$\forces{\poP}{\xmbox{\xitemof{III:models-1} }\sim\xmbox{ \xitemof{III:models-22}}}$ as 
in \assertof{1} and $\forces{\poP}{\xmbox{\xitemof{III:models-22-0}}}$  by 
\Thmof{III:P-msi-7}\,,\assertof{2$'$} and \Propof{III:P-gen-large-4}. 
\qedofThm

\begin{Thm}
  \Label{III:P-models-2} \wassertof{1} Suppose that the existence of two supercompact cardinals is 
  consistent. Then the following combination of principles is also 
  consistent:
  \begin{xitemize}
  \xitemcite[III:models-1] 
    $\SDLS^{int}_+(\calL^{\aleph_0}_{stat},\LT\continuum)$, $\GRP^{\LT\continuum}(\LE\continuum)$, 
  \xitemcite[III:models-2] 
    $\SDLS^{int}_+(\calL^\PKL_\stat,\LT\continuum)$;
  \xitem[III:models-23] 
    There is an inner model\/ $M$ of\/ $\uniV$ \st\ $(\continuum)^\uniV=(\continuum)^M$ and $V$ is 
    reached from $M$ by the forcing with a regular sub-\po\ of the completion of the 
    product of ccc and $\LT\continuum$-closed \pos, and  
    $\MA^{++\eta}(\calP,\LT\continuum)$ for
    $\calP=\setof{\poP}{\poP\mbox{ is a ccc \po\ }\poP\in M}$\\for all
    $\eta<\continuum$, and 
  \xitem[III:models-24] 
  	$\neg\,\HH(\LT\continuum)$. 
  \end{xitemize}
  \wassert{2} If there is a superhuge cardinal and a supercompact cardinal above it, then 
  the combination of the principles \xitemof{III:models-1}, \xitemof{III:models-2}, 
  \xitemof{III:models-23} and \xitemof{III:models-24} 
  above together with
  \begin{xitemize}
  \xitemcite[III:models-22-0] 
    $\Pkl{\continuum}{\lambda}$ carries a $\sigma$-saturated normal ideal for all
    $\lambda\geq2^{\aleph_0}$
  \end{xitemize}
  is consistent. 

\end{Thm}
\prf \assertof{1}: Let $\seqof{\poO_\alpha,\utpoR_\beta}{\alpha\leq\kappa,\beta<\kappa}$ be the following 
modification of the preparatory ccc finite support iteration \xitemof{III:models-4} in 
\Propof{III:P-models-0-0},\,\assertof{1}:
\begin{xitemize}
\xitem[III:models-25] $\utpoR_\beta=\left\{{}
  \begin{array}{@{}ll}
    (P)^\surd_{\poO_\beta},\ &\mbox{if }\beta=\alpha+1\mbox{ for an }\alpha\in S\mbox{ and}\\
    &f(\alpha)=\pairof{\mu,\theta,P}\mbox{ for cardinals }\mu,\,\theta\mbox{ and a \po\ }P\\
    &\mbox{\st\ }\forces{\poO_\beta}{(P)^\surd_{\poO_\beta}\mbox{ is a ccc \po}};\\[\jot]
    \mbox{\rlap{Hechler real forcing over \(\poO_\beta\),\qquad if \(\beta\in T\) but}}\\
    &\beta\mbox{ is not a successor of an element of }S\,;\\[\jot]
    \ssetof{\bbone_{\utpoR_\beta}}, &\mbox{otherwise.}
  \end{array}\right.$
\end{xitemize}

Let $\seqof{\poP_\alpha,\utpoQ_\beta}{\alpha\leq\kappa,\beta<\kappa}$ be the Easton-type 
mixed support iteration over  $\seqof{\poO_\alpha,\utpoR_\beta}{\alpha\leq\kappa,\beta<\kappa}$ and
$\poP=\poP_\kappa\ast(\Col(\kappa,\kappa_1))^\bullet_{\poP_\kappa}$.

Then, the generic extension of $\uniV$ by
$\poP=\poP_\kappa\ast(\Col(\kappa,\kappa_1))^\bullet_{\poP_\kappa}$ is as desired:  
$\forces{\poP}{\mbox{\xitemof{III:models-0}, \xitemof{III:models-1}, \xitemof{III:models-2}}}$
by \Thmof{III:P-models-0}.

$\forces{\poP}{\mbox{\xitemof{III:models-23}}}$ follows from (the proof of) 
\Propof{III:P-models-0-0},\,\assertof{1}. Note that the proof of 
\Propof{III:P-models-0-0},\,\assertof{1} does not rely on the value of $\poR_\beta$ for
$\beta\in T$ which is not a successor of the element of $S$.

Now the Hechler part of the preparatory iteration introduces 
an $\leq^*$-increasing sequence $\vec{h}=\seqof{f_\alpha}{\alpha<\continuum}$ of functions 
of length $\continuum$ (in the generic extension by $\poP_\kappa$)
$\calF=\setof{f_\alpha}{\alpha<\continuum}$ is still unbounded in the generic extension by
$\poP_\kappa$ by the genericity of $f_\alpha$'s: $f_\alpha$'s may no more Hechler reals 
above corresponding intermediate models in $\uniV^{\poP_\kappa}$ but each of them adds a 
Cohen real as its coordinatewise summand (see [Truss\cite{III:bib-truss}]). 
$\calF$ remains 
unbounded in $\poP$-generic extension $\uniV[\genG_\kappa\ast\genH]$ since no new reals are added by
$\Col(\kappa,\kappa_1)$. Thus,  
in $\uniV[\genG]$, the first 
countable topological 
space $X_\calF$ constructed in \sectionof{III:hamburger} 
is non-metrizable but all 
subspaces of $X_\calF$ of size $\LT2^{\aleph_0}$ are metrizable. Thus
$\forces{\poP}{\mbox{\xitemof{III:models-24}}}$. \smallskip

\assertof{2}: 
Let $\poP=\poP_\kappa\ast(\Col(\kappa,\kappa_1))^\bullet_{\poP_\kappa}$ be as 
in \assertof{1} with superhuge $\kappa$ and supercompact 
$\kappa_1$ above $\kappa$. Then
$\forces{\poP}{\xmbox{\xitemof{III:models-1},\xmbox{\xitemof{III:models-2},\xitemof{III:models-23}, 
      \xitemof{III:models-24}}}}$
as 
in \assertof{1} and $\forces{\poP}{\xmbox{\xitemof{III:models-22-0}}}$  by 
\Thmof{III:P-msi-7}\,,\assertof{2$'$} and \Propof{III:P-gen-large-4}. 
\qedofThm
\qedskip

We end up with mentioning some remaining open problems. As noted in 
\sectionof{III:hamburger}, 
Hamburger's Problem i.e.\ the consistency of $\HH(\LT\aleph_2)$ is still widely open. 
Galvin's Conjecture is also a persistingly open problem which can be discussed in our context 
(see e.g. [Todorcevic \cite{III:bib-stevo}]).

Both of the following two problems, which might be more at hand, are related to the last theorem in this section:
\begin{Problem}
Can we have the full $\MA^{++\mu}(ccc,\LT\continuum)$ together with all other strong reflection 
properties in some modification of the model of \Thmof{III:P-models-2}?
\end{Problem}

\begin{Problem}
  What is the $\Refl_\HH$ in the model of \Thmof{III:P-models-2}? Can we make it
  $(\continuum)^+$ or $\infty$ by some modification of the construction in the proof?
\end{Problem}

\phantomsection
\addcontentsline{toc}{section}{References}


\begin{thebibliography}{99}
\Label{III:ref}




\bibitem{III:bib-bm} Joan Bagaria and Menachem Magidor, On $\omega_1$-strongly compact 
  cardinals, Vol.79, (1),  (2014), 266--278. 
\bibitem{III:bib-bin} R.H.\,Bing, Metrization of topological spaces, Canadian Journal of Mathematics Vol.3 (1951), 175--186.
\bibitem{III:bib-cox} Sean Cox, The diagonal reflection principle, 
Proceedings of the American Mathematical Society, Vol.140, No.8 (2012), 
2893-2902. 
\bibitem{III:bib-cummings} James Cummings, Iterated Forcing and Elementary Embeddings, in:
  (Matthew Foreman and Akihiro Kanamori, eds.) Handbook of set Theory, Vol.2,  
  (2009), 775--884. 
\bibitem{III:bib-vdouwen} Eric K.\,van Douwen, The integers and topology, in: K.\,Kunen and 
  J.\,Vaughan (eds.), Handbook of Set-Theoretic Topology, Elsevier, (1984). 
\bibitem{III:bib-dow} Alan Dow, An introduction to applications of elementary submodels 
  to topology, Topology Proceedings 13, No.1 (1988), 17--72. 
\bibitem{III:bib-dtw2} A.\,Dow, F.D.Tall, and W.A.R.\,Weiss,  New proofs of 
  the consistency of the normal Moore space conjecture I, II, Topology and its 
  Applications, 37 (1990), 33--51, 115--129. 
\bibitem{III:bib-erice} Saka\'e Fuchino, Istvan Juh\'asz, Lajos Soukup, Zolt\'an Szentmikl\'ossy, 
   and Toshimichi Usuba, Fodor-type Reflection Principle and
  reflection of metrizability and meta-Lindel\"ofness, Topology and its Applications, Vol.157, 8 (2010), 1415--1429. 
\bibitem{III:bib-I} Saka\'e Fuchino, Andr\'e Ottenbereit 
  Maschio Rodrigues, and Hiroshi Sakai, Strong downward L\"owenheim-Skolem 
  theorems for stationary logics, I,\\to appear in Archive for Mathematical Logic. \\
  Extended version of the paper:\ \  \scalebox{0.8}[1]{\url{https://fuchino.ddo.jp/papers/SDLS-x.pdf}}
\bibitem{III:bib-II} \bysame{Fuchino, Saka\'e etal}: Strong downward L\"owenheim-Skolem 
  theorems for stationary logics, II --- 
  reflection down to the continuum, 
  to appear in Archive for Mathematical Logic. \\
  Extended version of the paper: \scalebox{0.75}[1]{\url{https://fuchino.ddo.jp/papers/SDLS-II-x.pdf}}
\bibitem{III:bib-refl-conti} Saka\'e Fuchino and Andr\'e Ottenbreit Maschio Rodrigues, 
Reflection principles, generic large cardinals, and the Continuum Problem, to 
appear in the Proceedings of the Symposium on Advances in Mathematical Logic 
2018.\ \ 
Extended version of the paper:\ \ \\
    \scalebox{0.67}[1]{\url{https://fuchino.ddo.jp/papers/refl\_principles\_gen\_large\_cardinals\_continuum\_problem-x.pdf}}
\bibitem{III:bib-hajnal-juhasz} A.\,Hajnal and I.\,Juh\'asz, On spaces in which every 
  small subspace is metrizable, Bulletin de lAcad\'emie Polonaise des Sciences, 
  S\'erie des sciences math, astr.\ et phys, Vol.24, No9, (1976), 727--731.
\bibitem{III:bib-millennium-book} Thomas Jech, Theory, The Third Millennium
	Edition, Springer (2001/2006).
\bibitem{III:bib-kanamori} Akihiro Kanamori, {\em The Higher Infinite}, Second 
  Edition, Springer Monographs in Mathematics, Springer-Verlag, (2003/2009). 
\bibitem{III:bib-krueger1} John Krueger, A general Mitchell style iteration, Mathematical Logic 
  Quarterly, 54(6), (2008), 641--651.
\bibitem{III:bib-krueger2} \bysame{John Krueger}, Some applications of mixed support iterations, 
  Annals of Pure and Applied Logic, 158(1-2), (2009), 40--57.
\bibitem{III:bib-kunen} K.\ Kunen, {\em Set Theory, An Introduction to Independence Proofs}, 
North-Holland (1980).
\bibitem{III:bib-andre} Andr\'e Ottenbreit Maschio Rodrigues, 
Some reflection principles in generic extensions by mixed support iteration 
\ifarxived
\mbox{(\ \hspace{30.6ex}}\\
\mbox{\hspace{34.2ex}\ )}, PhD thesis,
\else
（ミックスサポート反復強制によるジェネリック拡大におけるいくつかの反映原理）, PhD thesis, 
\fi
Kobe University (2020). 
\bibitem{III:bib-sakai} Hiroshi Sakai, Posets of cardinality $\kappa$ which destroy stationary subsets of
$\Pkl{}{}$, in preparation.
\bibitem{III:bib-stevo} Stevo Todorcevic, Combinatorial dichotomies in set theory, 
The Bulletin of Symbolic Logic, Vol.17, No.1, (2011), 1--72. 
\bibitem{III:bib-truss} John Truss,  Sets having calibre $\aleph_1$, in: R.Gandy, M.Hyland 
  (Eds.), Logic Colloquium 76 , North-Holland PUblishing Company (1977), 595--612.
\end{thebibliography}
\end{document}